\documentclass[11pt,twoside]{article}
\usepackage{amsmath,amsfonts,amssymb,cases,mathdots,color,colordvi,url,tikz}
\usepackage{hyperref}
\usepackage{bm}
\usepackage{bbm}
\usepackage{pgfpages}
\usepackage{tabularx}
\usepackage{graphics}
\usepackage{subfigure}
\usepackage{epsfig}
\usepackage[figuresright]{rotating}
\usepackage{rotating}
\usepackage{algorithm}
\usepackage{algorithmic}
\usepackage{physics}
\usepackage{epstopdf}
\usepackage{parskip}
\usepackage{footmisc}
\usepackage{tablefootnote}
\usepackage[justification=centering]{caption}

\newtheorem{proposition}{Proposition}
\newtheorem{theorem}[proposition]{Theorem}
\newtheorem{lemma}[proposition]{Lemma}

\newtheorem{definition}[proposition]{Definition}

\newtheorem{remark}{Remark}
\newtheorem{assumption}{Assumption}
\newtheorem{example}{Example}

\newcommand{\be}{\begin{equation}}
	\newcommand{\ee}{\end{equation}}
\newcommand{\ba}{\begin{eqnarray}}
	\newcommand{\ea}{\end{eqnarray}}
\newcommand{\bas}{\begin{eqnarray*}}
	\newcommand{\eas}{\end{eqnarray*}}

\graphicspath{ {./Figures/} }

\def\H{{\mathcal{H}}}

\def\mbR{{\mathbb R}}

\def\S{{\mathcal S}}
\def\C{{\mathcal C}}
\def\F{{\mathcal F}}

\def\J{{\mathcal J}}

\def\N{{\mathcal N}}
\def\O{{\mathcal O}}

\def\Y{{\mathcal Y}}

\def\bfp{{\bf p}}
\def\bfq{{\bf q}}
\def\bfr{{\bf r}}
\def\bfs{{\bf s}}

\def\bfa{{\bf a}}
\def\bfb{{\bf b}}
\def\bfc{{\bf c}}
\def\bfd{{\bf d}}
\def\bfe{{\bf e}}

\def\bfv{{\bf v}}
\def\bfx{{\bf x}}
\def\bfy{{\bf y}}
\def\bfz{{\bf z}}
\def\bfw{{\bf w}}
\def\bfu{{\bf u}}

\def\bfone{{\bf 1}}
\def\Prox{\mbox{Prox}}
\def\OG{{\overline{\Gamma}}}

\def\wz{{\widetilde{\zeta}}}
\def\obz{\overline{\bfz}}
\def\bfdt{\boldsymbol{\delta}}
\def\proxg{ {\rm Prox}_{\tau g(\cdot)} }

\def\bp{ \textbf{Proof.} }
\def\ep{ \hfill $\Box$ }
\def\bfta{{\bf {\boldsymbol \tau}}}
\def\oT{ \overline{T} }
\def\wz{ \widetilde{\bfz} }
\def\cone{(\textbf{C1})}
\def\ctwo{(\textbf{C2})}
\def\mbN{ \mathbb{N} }
\def\wq{ \widetilde{\bfq} }
\def\oj{{\overline{j}}}
\def\ul{{\underline{l}}}
\def\ol{{\overline{l}}}
\def\oS{\overline{\S}}
\def\ode{\overline{\delta}}
\def\bfla{\boldsymbol{\lambda}}
\def\whp{\widehat{\partial}}

\def\bfmu{\boldsymbol{\mu}}
\def\indf{\boldsymbol{\mathcal{I}}}

\def\VPO{ \texttt{VPO} }
\def\VDO{ \texttt{VDO} }
\def\TIME{ \texttt{TIME} }
\def\AUC{ \texttt{AUC} }

\title{Composite Optimization with Indicator Functions: \\ Stationary Duality and a Semismooth Newton Method}

\author{Penghe Zhang\thanks{Department of Data Science and Artificial Intelligence, The Hong Kong Polytechnic University, Hong Kong SAR, China, E-mail: {penghe.zhang@polyu.edu.hk} },
\ \ Naihua Xiu\thanks{School of Mathematics and Statistics, Beijing Jiaotong University, Beijing 100044, China, E-mail: {nhxiu@bjtu.edu.cn} } \ \ and \ \
Hou-Duo Qi\thanks{Department of Data Science and Artificial Intelligence, and Department of Applied Mathematics, The Hong Kong Polytechnic University, Hong Kong SAR, China, E-mail: {houduo.qi@polyu.edu.hk} }
	}

\date{}

\begin{document}
	
	\maketitle	

\begin{abstract}
	Indicator functions of taking values of zero or one are essential to numerous applications in machine learning and statistics.
	The corresponding primal optimization model has been researched in several recent works. However, its dual problem is a more challenging topic that has not been well addressed. One possible reason is that the Fenchel conjugate of any indicator function is finite only at the origin. This work aims to explore the dual optimization for the sum of a strongly convex function and a composite term with indicator functions on positive intervals. 
	For the first time, a dual problem is constructed by extending the classic conjugate subgradient property to the indicator function. This extension further helps us establish the equivalence between the primal and dual solutions. The dual problem turns out to be a sparse optimization with a $\ell_0$ regularizer and a nonnegative constraint. The proximal operator of the sparse regularizer is used to identify a dual subspace to implement gradient and/or semismooth Newton iteration with low computational complexity. This gives rise to a dual Newton-type method with both global convergence and local superlinear (or quadratic) convergence rate under mild conditions. Finally, when applied to AUC maximization and sparse multi-label classification, our dual Newton method demonstrates satisfactory performance on computational speed and accuracy.
	
		\vspace{3mm}
	
\end{abstract}

\begin{sloppypar}
	\section{Introduction} 
	In this work, we aim to solve the following nonconvex composite optimization
	\begin{equation} \label{01cop}
		\min_{\bfx \in \mbR^n} f( \bfx ) + \indf(A\bfx + \bfb),
	\end{equation}
	where $f:\mathbb{R}^n \to \mbR$ is a $\sigma_f$-strongly convex function (possibly nonsmooth), and $A \in \mbR^{m \times n}$ and $\bfb \in \mbR^m$ are given data. For $\bfu = (u_1,\cdots,u_m)^\top$, we define $\indf (\bfu) = \lambda \sum_{i = 1}^{m} \bfone_{(0,\infty)} (u_i) $, where $\lambda$ is a positive constant and $\bfone_{(0,\infty)}$ is the indicator function\footnote{This paper adopts the terminology ``indicator function" from the statistical learning theory 
			\cite[Page 19]{vapnik1999nature},
			differing from the usage in optimization where the function takes 0 on a given set and $\infty$ elsewhere. It is also called the characteristic function \cite{zadeh1965fuzzy} in other fields.} 
	on the positive interval (also called zero-one loss), defined as
	\begin{align*}
		\bfone_{(0,\infty)} ( t ) = \left\{ \begin{aligned}
			& 1, \quad \mbox{   if } t > 0, \\
			&0, \quad \mbox{   if } t \leq 0.
		\end{aligned} \right.
	\end{align*}
	Problem \eqref{01cop} arises from extensive applications, including multi-label classification \cite{zhang2013review}, area under the receiver operating characteristic curve (AUC) \cite{gao2013one,hanley1982meaning}, maximum ranking correlation \cite{han1987non,kendall1938new}, to mention a few. The strongly convex function $f$ often serves as a regularizer in these tasks, including squared $\ell_2$ norm for margin maximization (e.g. \cite{CV95,shivaswamy2010maximum}) or (group) elastic net for feature selection (e.g. \cite{zou2005regularization,munch2021adaptive,xie2023group}). 
	Problem \eqref{01cop}, in particular, covers the composite $\ell_0$-norm:
	$
	\| B\bfx + \bfb \|_0 
	= \indf(B\bfx + \bfb) + \indf(-B\bfx- \bfb).
	$
	The matrix $A$ becomes $[B; -B]$ (Matlab notation) and this implies that
	the rows of $A$ are dependent.
	Therefore, it is essential we do not assume that $A$ has a full-row rank 
	property. And we do not assume that $f$ is smooth.

	When $f$ is smooth, several progresses have been made in solving the primal nonconvex composition problem involving the indicator function or a lower semi-continuous (l.s.c) function. These methods include, but are not limited to the Newton method \cite{zhou2021quadratic}, augmented Lagrangian-based method \cite{bot2019proximal,boct2020proximal,li2015global,wang2018con}, lifted reformulation \cite{cui2023minimization,hananalysis}, and complementarity reformulation \cite{feng2013complementarity,kanzow2025sparse}. However, the composition of a nonconvex, discontinuous indicator function and a large-scale matrix brings significant difficulty in both convergence and computation. For example, certain regularity conditions (such as a full row-rank matrix $A$) are often required for the convergence analysis of the aforementioned algorithms. 
	Moreover, it still remains a difficult task
	to globalize a primal Newton method while retaining local quadratic convergence.
	Finally, we should stress that the function $f$ in \eqref{01cop} is not necessarily smooth. This together with the combinatorial nature of $\indf$ 
	poses significant barriers to designing a globally convergent algorithm for \eqref{01cop}.

	This paper offers an alternative methodology, strongly motivated by the classical 
	duality theory of convex optimization.
	We will propose a dual reformulation of \eqref{01cop} with a conducive structure for numerical computation. 
	In the following, we first explain the motivation that leads to a framework for
	constructing the dual problem.
	We then briefly summarize the main contributions of the paper.

	\subsection{Motivation from convex optimization}
	
	Let us recall the Fenchel duality in convex optimization (see \cite[Chap.~31]{rockafellar1970convex}
	and \cite{beck2014fast,chen1994proximal}):
	\begin{align}  \label{convex-dual}
		\inf_{\bfx \in \mbR^n}~ 
		\Big\{ f(\bfx) + \varphi (A\bfx + \bfb) \Big\} 
		=	\sup_{\bfz \in \mbR^m} 
		\left\{ -f^*(-A^\top \bfz) + \langle \bfb, \bfz \rangle -\varphi^*(\bfz) \right\}.
	\end{align}
	where both $f$ and $\varphi$ are closed, proper, and convex, and
	$f^*$ and $\varphi^*$ are their respective conjugate functions.
	For simplicity of discussion, we omit the relative interior conditions for the duality
	to hold.
	The dual problem is often advantageous to solve due to a few features.
	Firstly, the matrix $A$ has been separated from the possibly nonsmooth function $\varphi$ while composited with $f^*$. 
	Secondly, $f^*$ has Lipschitzian continuous gradient when $f$ is strongly convex. 
	Thirdly, in many applications, $\varphi^*$ has a closed-form proximal operator that would
	allow proximal gradient method to be developed, see \cite{beck2014fast}.
	An impressive progress has been made on applying the dual approach to LASSO-type problems and highly efficient semismooth Newton algorithms are developed, see 
	\cite{li2018highly,zhang2020efficient} and the references therein. Recently, \cite{gfrerer2021semismooth} proposed the semismoothness* of a set-valued mapping, which is an improvement of classic semismoothness. For minimization of a prox-regular function, \cite{mordukhovich2021generalized} designed a generalized Newton method based on the second-order subdifferential \cite{mordukhovich1992sensitivity} and proved its local superlinear convergence under semismoothness* and tilt stability assumptions. This method is further developed for solving composite optimization \cite{duy2023generalized,khanh2024coderivative,khanh2024globally} and difference programming \cite{aragon2024coderivative} under some assumptions. 
		A more comprehensive introduction to this method can be found in the monograph \cite{mordukhovich2024second}.
	These works motivated us to investigate whether a dual approach is viable for
	\eqref{01cop} and whether a generalized Newton method can be constructed.
	
	However, if we blindly apply the dual formulation \eqref{convex-dual} to \eqref{01cop}
	with $\varphi = \indf$,
	we would not gain much. 
	This is because the conjugate function $\indf^*$ only takes finite value at $\{0\}$:
	\begin{equation*}
		\indf^*(\bfz) := \sup_{\bfu \in \mbR^m} 
		\Big\{ \langle \bfz, \bfu \rangle - \indf(\bfu) \Big\}
		= \left\{ \begin{aligned}
			& 0, &&\mbox{if } \bfz = 0, \\
			& \infty, &&\mbox{otherwise}.
		\end{aligned} \right. 
	\end{equation*}
	This means that the dual problem would have only one feasible point $\bfz =0$ and everything
	we care about in duality breaks down.

\begin{figure}[h]
		\subfigure{
			\begin{minipage}[t]{0.5\linewidth}
				\centering
				\includegraphics[width=6in]{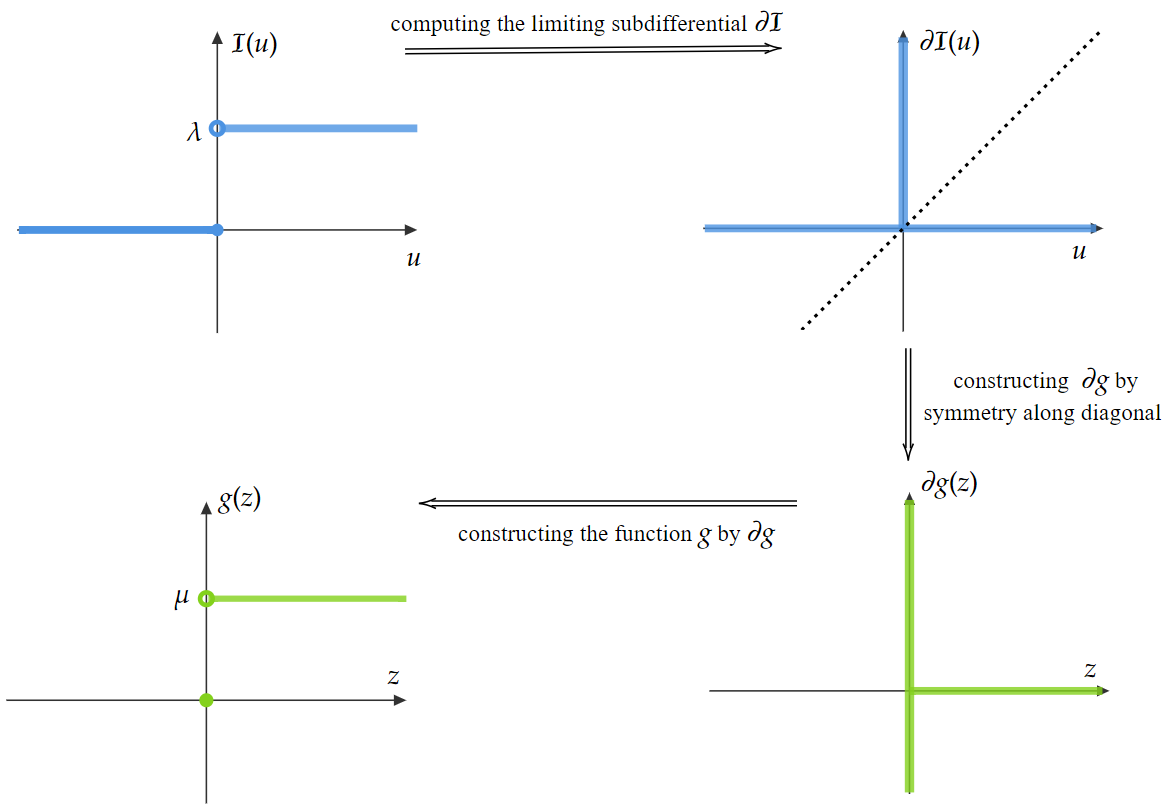}
			\end{minipage}%
		}
		\caption{Constructing a function $g$ such that $z \in \partial \indf(u)$ if and only if $u \in \partial g(z)$ for any given $z,u \in \mbR$.}
		\label{Fig-g}
	\end{figure}

	Now let us take a step back to recall one essential property for the conjugate 
	function pair $(\varphi, \varphi^*)$ when $\varphi$ is proper, closed, and convex:
	\begin{align*}
		\bfz \in \partial \varphi(\bfu)~~ \Longleftrightarrow ~~ \bfu \in  \partial \varphi^*(\bfz),
	\end{align*}
	where $\partial \varphi(\cdot)$ is the subdifferential of $\varphi(\cdot)$ 
	in convex analysis \cite{rockafellar1970convex}.
	We call this property the stationary duality property. 
	We would like to find a function $g(\cdot)$ that plays the role of $\indf^*(\cdot)$ and
	satisfies the stationary duality property:
	\begin{align} \label{inverse-I-g}
		\bfz \in \partial \indf(\bfu)~~ \Longleftrightarrow ~~ \bfu \in  \partial g(\bfz), 
	\end{align}
	where $\partial \indf(\cdot)$ is the limiting subdifferential in variational analysis 
	\cite{mordukhovich2018variational,RockWets98}.
	In the meantime, we require the function $g(\bfz)$ to inherit some useful information
	from the primal Problem \eqref{01cop}.
	For the one-dimensional case, the construction of such function $g(z)$ 
	is demonstrated in Fig.~\ref{Fig-g}.
	
	\subsection{Stationary duality and main contributions}
	
	Extension of the one-dimensional case to the multi-dimensional case leads to the following
	function:
	\be \label{g-function}
	g(\bfz) := \mu \| \bfz \|_0 + \bfdt_+(\bfz),
	\ee
	where $\mu>0$ is a parameter and $\bfdt_+(\cdot)$ 
	is defined as follows:
	\begin{align*}
		\bfdt_+(\bfz) = \left\{ \begin{aligned}
			& 0, &&\mbox{ if } \bfz \geq 0, \\
			& +\infty, &&\mbox{ otherwise}.
		\end{aligned} \right.
	\end{align*}
	Using the characterizations in \eqref{partial01} and \eqref{subdif-g},
	it is easy to verify that the stationary dual relationship \eqref{inverse-I-g} holds for
	$g(\cdot)$.
	This gives rise to the stationary dual problem:
	\begin{equation} \label{01dco-max} 
		\max_{\bfz \in \mathbb{R}^m}~  
		- h(\bfz) -  g(\bfz) \qquad \mbox{with} 
		\quad
		h(\bfz) := f^*(-A^\top \bfz) - \langle \bfz, \bfb \rangle .
	\end{equation}
	Although we cannot expect the perfect duality in convex optimization, we
	are able to establish the following important results.
	
	\begin{itemize} 
		
		\item[(i)] There is one-to-one correspondence between local minimizers of the primal problem
		\eqref{01cop} and that of the dual problem \eqref{01dco-max}. See Thm.~\ref{foc} for the detailed
		formula for computing the corresponding solution from the other.
		This is achieved through the study of proximal-type (``P-" for short) stationary point of \eqref{01dco-max} and its 
		correspondence to the KKT (Karush-Kuhn-Tucker) point of a local convex reformulation of \eqref{01cop} at a local minimizer. 
		This result is the basis for further numerical development.

		\item[(ii)] Since $f$ is strongly convex, its conjugate $f^*$ has Lipschitzian continuous gradient and well-defined generalized Jacobian.
		More importantly, the nonsmooth function $g$ 
		is a sparse regularizer with a closed-form proximal operator, which helps us identify a dual subspace to perform low-complexity iteration. In this reduced subspace, we can first implement a gradient step and then accelerate it by a semismooth Newton step. In this process, easy-to-check conditions are proposed to ensure the dual objective function in \eqref{01dco-max} ascent. These steps form the basic framework of our subspace gradient semismooth Newton (SGSN) method, see Alg.~\ref{SGSN}.
		
		\item[(iii)] In addition to the low per-iteration complexity, our SGSN has global convergence with local superliner (or quadratic) rate. If the generated sequence is bounded and the dual objective function in \eqref{01dco-max} satisfies Polyak-{\L}ojasiewicz-Kurdyka (P{\L}K) property, then the sequence converges to a {\rm P}-stationary point of \eqref{01dco-max}
		(Thm.~\ref{Thm-Global}). Thereby, a local minimizer of the primal Problem \eqref{01cop} is found. This result does not require the regularity condition that the matrix $A$ has full-row rank. If we further assume that the gradient function $\nabla f^*$ is (strongly) semismooth and a second-order sufficient condition holds, then the sequence converges to a {\rm P}-stationary point of \eqref{01dco-max} with local superlinear (or quadratic) rate. To the best of our knowledge, SGSN is the first algorithm that enjoys both global convergence and local superlinear rate for Problem \eqref{01cop} (Thm.~\ref{Thm-Local}). 
		Finally, we demonstrate the use of SGSN on two important applications: AUC maximization and the sparse multi-label classification. The experiments show that SGSN has satisfactory performance on computational speed and solution accuracy when benchmarking against several competitive algorithms.
				
	\end{itemize}

	\subsection{Organization}
	The paper is organized as follows. Some frequently used symbols and notations are
	introduced in Section \ref{pre}. The one-to-one correspondence between solutions of primal and dual problems is established in Section \ref{pd-solution}. 
	To find a dual solution, we describe the framework of SGSN and prove its global convergence with a local superlinear (or quadratic) rate in Section \ref{sec_sgsn}. The numerical experiments on AUC maximization and multi-label classification are conducted in Section \ref{numerical-experiments}. 
	We conclude the paper in Section \ref{conclusion} with some discussion on possible 
	future research.
	
	\section{Preliminaries} \label{pre}
	
	In this part, we first describe the notation used in the paper. 
	We then recall the definition of the limiting subdifferential and apply it to the function $g(\cdot)$ 
	in \eqref{g-function}.
	We also characterize its proximal operator. 
	Finally, we introduce some function classes. 
	
	\subsection{Notation} \label{Subsection-Notation}
	
	The boldfaced lowercase letter $\bfx \in \mathbb{R}^n$ denotes a column vector of size $n$ and $\bfx^\top$ is its transpose.
	Let $x_i$ or $[\bfx]_i$ denote the $i$th element of $\bfx$.
	We let $\mbR^n_+$ denote the non-negative orthant in $\mathbb{R}^n$.
	The norm $\| \cdot\|$ denotes the $\ell_2$ norm of $\bfx$ or the Frobenius norm for
	a matrix $A$, and thus we always have
	$
	\| A\bfx \| \le \| A \| \| \bfx\|.
	$
	A neighborhood of $\bfx^* \in \mathbb{R}^n$ is often denoted as $\mathcal{N}(\bfx^*)$,
	while $\mathcal{N}(\bfx^*, \delta) := \{ \bfx \in \mathbb{R}^n \ | \ \| \bfx - \bfx^* \| < \delta \}$ when the radius $\delta > 0$ is given. 
	We denote 
	$
	I
	$ 
	as the identity matrix of appropriate dimension. 
	Let $\bfe$ represent a vector of all 1's with dimension implied by the context.
	We denote ceiling and floor functions as $\lceil \cdot \rceil$ and $\lfloor \cdot \rfloor$ respectively.
	
	Let $[m]$ denote the set of indices $\{1, \ldots, m\}$.
	For a subset $\Gamma \subset [m]$,  $\OG$ consists of those indices of $[m]$ not in $\Gamma$ and $|\Gamma|$ denotes the number of elements in $\Gamma$ (cardinality of $\Gamma$).
	For vector $\bfz \in \mathbb{R}^m$ (resp. matrix $A \in \mbR^{m \times n}$), $\bfz_{\Gamma}$ (resp. $A_{\Gamma:}$) denotes the subvector of
	$\bfz$ indexed by $\Gamma$ (resp. the submatrix of $A$ with rows indexed by $\Gamma$). Particularly, for a differentiable function $h$, we denote $\nabla_{\Gamma} h(\bfz):= (\nabla h(\bfz))_{\Gamma}$. For a symmetric matrix $H \in \mbR^{m \times m}$, $H_{\Gamma}$ is the submatrix with both rows and columns indexed by $\Gamma$. Meanwhile, the largest and smallest eigenvalue of $H$ are denoted by $\bfla_{\max}$ and $\bfla_{\min}$ respectively. 
	Given vector $\bfz \in \mathbb{R}^m$, we denote its support set by 
	$\mathcal{S}(\bfz):= \{ i\in [m]: z_i \neq 0  \}$, where ``$:=$'' means ``define''.
	For a set $\Omega \subset \mbR^m$, ${\rm co}\{ \Omega \}$ is the convex hull of
	$\Omega$.
	
	\subsection{Limiting subdifferential} 
	
	The major reference for this part is \cite{mordukhovich2006variational,mordukhovich2018variational,RockWets98}.
	Let us consider a proper, lower semi-continuous (l.s.c.) and bounded-below function $\varphi: \mathbb{R}^m  \to (-\infty, +\infty]$. The regular subdifferential of $\varphi$ 
	at $\obz \in {\rm dom}\,\varphi := \{ \bfx \in \mathbb{R}^m: \varphi(\bfx) < \infty \}  $ 
	is defined by
	\begin{align*}
		& \widehat{\partial} \varphi (\obz) := \left\{  \bfv \in \mbR^m : \liminf_{\bfz \to \overline{\bfz} \atop \bfz \neq \overline{\bfz}} \frac{\varphi(\bfz) - \varphi(\overline{\bfz}) - \langle \bfv, \bfz - \obz \rangle }{\| \bfz - \obz \|} \geq 0  \right\} .
	\end{align*}
	For $\obz \notin \hbox{dom}\,\varphi$, we set $\partial \varphi (\obz) = \varnothing$. The limiting or Mordukhovich subdifferential of $\varphi$ at $\obz \in \mbR^m $ is defined by 
	\begin{align*}
		\partial \varphi (\obz):= \{ \bfv \in \mbR^m: \exists \bfz^k \to \obz,~ \varphi(\bfz^k) \to \varphi (\obz) \mbox{ and }  \widehat{\partial} \varphi (\bfz^k)  \ni  \bfv^k \to  \bfv \} .
	\end{align*}
	It was initially proposed in \cite{mordukhovich1976maximum} in an equivalent formulation. Particularly, the limiting subdifferential of $\indf(\cdot)$ 
	takes the following form \cite{zhou2021quadratic}:
	\be \label{partial01}
	\partial  \indf(\bfu)  =
	\left\{
	\bfv \in \mathbb{R}^m \; \left|
	\ v_i \left\{
	\begin{array}{ll}
		\ge 0 , \ & \mbox{if} \ u_i = 0 \\
		=0,     \ & \mbox{if} \ u_i \not=0,
	\end{array} , \ \ i \in [m]
	\right.
	\right .
	\right\} .
	\ee 
	For given $\tau > 0$, 
	the proximal operator of $\varphi$ is defined by
	\begin{eqnarray*}
		\Prox_{\tau \varphi(\cdot)}(\bfu) &:=& \arg\min_{\bfz \in \mathbb{R}^n}\left\{
		\varphi(\bfz) 
		+ \frac{1}{2\tau} \| \bfz - \bfu\|^2
		\right\}.
	\end{eqnarray*}
	
	The following result reveals the relationship of $\partial \indf$, $\partial g$ and $\proxg(\cdot)$, and they will be utilized to establish the equivalence of solutions between the primal and dual problems.
	Although it is straightforward, a proof is included for easy reference.
	
	\begin{lemma} \label{property-g}
		Given points $\bfz \in \mbR^m_+$ and $\bfu \in \mbR^m$, the following results hold.
		
		(i) The limiting subdifferential of $g$ at $\bfz$ can be computed by 
		\begin{align} \label{subdif-g}
			\partial g(\bfz) = \left\{ \bfv \in \mbR^m: v_i \in \left\{ \begin{aligned}
				& \mbR , &&\mbox{if } z_i = 0 \\
				& \{ 0 \},    &&\mbox{if } z_i > 0
			\end{aligned} \right.  \right\} .
		\end{align}
		
		(ii) The proximal operator of $g$ with $\tau > 0$ can be represented as
		\begin{equation}\label{proxg}
			\Big[ 	{{\rm Prox}_{\tau g (\cdot)}} ( \bfu )  \Big]_i 
			= \left\{  \begin{aligned}
				& 0, && u_i < \sqrt{2\tau\mu}, \\
				&\{ 0, u_i \}, && u_i = \sqrt{2\tau\mu}, \\
				&u_i, && u_i > \sqrt{2\tau\mu},
			\end{aligned} \right. \quad i=1, \ldots, m.
		\end{equation}
		
		(iii) $\bfz \in \partial \indf(\bfu)$ if and only if $\bfu \in \partial g(\bfz)$.
		
		(iv) $\bfz \in \proxg ( \bfz + \tau \bfu )$ holds with some $\tau > 0$ if and only if $\bfu$ and $\bfz$ take the following value
		\begin{equation} \label{uz-range}
			(u_i,z_i) \in \Omega_1 \cup \Omega_2, \ \ \ \ \forall i \in [m], 
		\end{equation}
		where $\Omega_1:= \{ 0 \} \times [\sqrt{2\mu\tau}, \infty)$ and $\Omega_2:= (-\infty, \sqrt{2\mu/\tau}] \times \{0\}$. 
		
		(v) If $\bfz \in \proxg ( \bfz + \tau \bfu )$ with some $\tau > 0$, then $\bfu \in \partial g(\bfz)$. Conversely, if $\bfu \in \partial g(\bfz)$, then $\bfz \in \proxg ( \bfz + \tau \bfu )$ holds with $\tau \in (0, \min\{ \bfta_1(\bfz), \bfta_2(\bfu) \})$, where 
		\begin{align*}
			& \bfta_1(\bfz) = \left\{  \begin{aligned}
				& +\infty, && \bfz = 0, \\
				&\min\left\{ \frac{z_i^2}{ 2\mu } : \ z_i > 0 \right\}, && {\rm otherwise}, 
			\end{aligned} \right. \\  
			&\bfta_2(\bfu) := \left\{  
			\begin{aligned}
				& +\infty,  && \bfu \leq 0, \\
				&\min\left\{ \frac{2\mu}{ u_i^2 }:\ u_i > 0 \right\}, &&{\rm otherwise}.
			\end{aligned} \right.
		\end{align*}
		
	(vi) If $\bfz \in \proxg ( \bfz + \tau \bfu )$ with some $\tau > 0$, then for any $\tau' \in (0,\tau)$, ${\rm Prox}_{\tau' g(\cdot)} ( \bfz + \tau' \bfu )$ is single-valued and it holds that $\bfz = {\rm Prox}_{\tau' g(\cdot)} ( \bfz + \tau' \bfu )$.
	\end{lemma}
	\bp
	(i) Let us consider the regular and limiting subdifferential of a l.s.c. function $g_1: \mathbb{R} \to \mbR\cup\{ \infty \}$, defined by
	\begin{align*}
		g_1(t):= \left\{ \begin{aligned}
			& \mu, && \mbox{if } t > 0, \\
			& 0, && \mbox{if } t = 0,  \\
			& \infty, && \mbox{otherwise.}
		\end{aligned} \right.
	\end{align*}
	Obviously, $g(\bfz) = \sum_{i=1}^m g_1(z_i)$ holds and it follows from \cite[Prop.~10.5]{RockWets98} that
	\begin{align}\label{subdif-sep0}
		\widehat{\partial} g(\bfz) = \widehat{\partial} g_1(z_1) \times \cdots \times \widehat{\partial} g_1(z_m)\quad \mbox{and}\quad \partial g(\bfz) = \partial g_1(z_1) \times \cdots \times \partial g_1(z_m).
	\end{align}
	Given $t\geq 0$, we claim that 
	\begin{align} \label{subdif-sep}
		\widehat{\partial} g_1(t) = \left\{ v \in \mbR: v \in \left\{ \begin{aligned}
			& \mbR , &&\mbox{if } t = 0, \\
			& \{ 0 \},    &&\mbox{if } t > 0 
		\end{aligned} \right.  \right\} .
	\end{align}
	Indeed, when $t > 0$, $g_1$ is differentiable and then $	\widehat{\partial} g_1(t) = g'_1(t) = 0$ from \cite[Exercise 8.8]{RockWets98}. If $t = 0$, then for any sequence $t_k \downarrow 0$ and $t_k \neq 0$, we have $g_1(t_k) - g_1(0) \geq \mu$. Thus, given any $v \in \mbR$, we have $(g_1(t_k) - g_1(0) - \langle v, t_k \rangle)/|t_k| \to +\infty$. Based on these facts, \eqref{subdif-sep} can be derived by the definition of regular subdifferential.
	We can further prove $\partial g_1(t) = \widehat{\partial} g_1(t)$ by the definition of limiting subdifferential. 
	Combining this result with \eqref{subdif-sep0}, we obtain \eqref{subdif-g}.
	
	(ii) By the definition of the proximal operator and the separable property of $\| \cdot \|^2$ and $\| \cdot \|_0$, we only need to solve the following optimization problem for each $i \in [m]$
	\begin{align} \label{prox-op}
		\Big[ 	{\rm Prox}_{\tau g (\cdot)} ( \bfu )  \Big]_i  = \arg \min_{z_i \geq 0} \frac{1}{2} (z_i - u_i)^2 + \tau \mu \| z_i \|_0.
	\end{align}
	When $z_i = 0$, the objective function value is $\texttt{obj}_1 = u_i^2/2$. 
	Now we consider a constrained quadratic optimization
	\begin{align*}  
		\min_{z_i \geq 0}\ \frac{1}{2} (z_i - u_i)^2 + \tau \mu .
	\end{align*}
	Its optimal solution is $z_i = \max\{ 0,u_i \}$ and the objective function value is $\texttt{obj}_2 = (\min\{ u_i,0 \})^2 /2 + \tau\mu$. Finally, we just need to compare $\texttt{obj}_1$ and $\texttt{obj}_2$ to determine the solution to \eqref{prox-op}. 
	
	Case I: if $\texttt{obj}_1 < \texttt{obj}_2$, this means $u_i < \sqrt{2\tau\mu}$ and $[ {\rm Prox_{\tau g (\cdot)}} ( \bfu ) ]_i = 0$.
	
	Case II: $\texttt{obj}_1 > \texttt{obj}_2$ is equivalent to $u_i > \sqrt{2\tau\mu}$ and $[ 	{\rm Prox_{\tau g (\cdot)}} ( \bfu ) ]_i = \max\{ 0,u_i \} = u_i$.
	
	Case III: if $\texttt{obj}_1 = \texttt{obj}_2$, this leads to $u_i = \sqrt{2\tau\mu}$ and $[ {\rm Prox_{\tau g (\cdot)}} ( \bfu ) ]_i = \{0,u_i\}$.
	
	We arrive at the conclusion of (ii).
	
	(iii) This result can be directly obtained from \eqref{partial01} and \eqref{subdif-g}. 
	
	(iv) Suppose $\bfz \in \proxg ( \bfz + \tau \bfu )$ holds, we can use \eqref{proxg} to derive the range of $(\bfu,\bfz)$. If $z_i + \tau u_i < \sqrt{2\tau\mu}$, then $z_i = 0$ and $u_i < \sqrt{2\mu/\tau}$. If $z_i + \tau u_i > \sqrt{2\tau\mu}$, then $z_i = z_i + \tau u_i$, leading to $u_i = 0$ and $z_i > \sqrt{2\mu\tau}$. For $z_i + \tau u_i = \sqrt{2\tau\mu}$, we have $z_i \in \{ 0, z_i + \tau u_i \}$, yielding $(u_i,z_i) = (0,\sqrt{2\tau\mu})$ or $(u_i,z_i) = (\sqrt{2\mu/\tau},0)$.
	
	Now let us prove the converse counterpart. If $(u_i,z_i) \in \Omega_1$, then $z_i + \tau u_i \geq \sqrt{2\tau\mu}$ and thus $z_i = z_i + \tau u_i \in [\proxg ( \bfz + \tau \bfu )]_i$. If $(u_i,z_i) \in \Omega_2$, then $z_i + \tau u_i \leq \sqrt{2\tau\mu}$ and thus $z_i = 0 \in [\proxg ( \bfz + \tau \bfu )]_i$. Overall, we complete the proof of (iii).
	
	(v) The first part of this assertion directly follows from (i) and (iv). To prove the converse part, it suffices to prove that $\bfu \in \partial g( \bfz)$ implies \eqref{uz-range} 
	with $\tau \in (0, \min\{ \bfta_1(\bfz), \bfta_2(\bfu) \})$. 
	Indeed, for any $i_0$ such that $z_{i_0} > 0$, we have $\tau < \bfta_1(\bfz) < z_{i_0}^2/(2\mu)$, implying $z_{i_0} > \sqrt{2\tau\mu}$. Meanwhile, for any $i_0$ such that $u_{i_0} > 0$, we have $\tau < \bfta_2(\bfu) < (2\mu)/u_{i_0}^2$, leading to $u_{i_0} < \sqrt{2\mu/\tau}$. These results together with $\bfu \in g(\bfz)$ yields $(u_i,z_i) \in \Omega_1 \cup \Omega_2$ for $i\in [m]$, and we obtain \eqref{uz-range}.
	
	(vi) It suffices to prove $\bfz = {\rm Prox}_{\tau' g(\cdot)} ( \bfz + \tau' \bfu )$. By (iv), we have $(u_i,z_i) \in \Omega_1 \cup \Omega_2$ for $i \in [m]$. Then let us consider the following two cases.
		
		Case I: If $(u_i, z_i) \in \Omega_1$, then $u_i=0$ and $z_i + \tau' u_i \geq \sqrt{2\mu\tau} > \sqrt{2\mu\tau'}$. Using \eqref{proxg}, we can derive $[{{\rm Prox}_{\tau' g (\cdot)}} ( \bfz + \tau' \bfu )]_i = z_i + \tau' u_i = z_i$.
		
		Case II: If $(u_i, z_i) \in \Omega_2$, then 
		$z_i=0$ and $u_i \le \sqrt{2\mu/\tau}$. Hence,
		$z_i + \tau' u_i  \leq \tau' \sqrt{2\mu/\tau} < \sqrt{2\mu\tau}$. Using \eqref{proxg}, we can derive $[{{\rm Prox}_{\tau' g (\cdot)}} ( \bfz + \tau' \bfu )]_i = 0 = z_i$.
		
		Therefore, we obtain the desired conclusion.
	\ep

	\subsection{Function classes}
	
	This paper involves a few classes of functions. We briefly introduce them below.
	
	\textbf{(A) Strongly convex functions, $L$-smooth functions and conjugate functions.}
	We say a function $\psi:\mathbb{R}^n \to (-\infty, +\infty]$ is $\sigma_{\psi}$-strongly convex for $\sigma_{\psi} > 0$, if dom$\,\psi$ is convex and the following inequality holds for any $\bfx,\bfy \in \mbox{dom }$ and $s \in [0, 1]$,
	\begin{align*}
		\psi(s\bfx + (1-s)\bfy) \leq s\psi(\bfx) + (1-s)\psi(\bfy) - \frac{\sigma_\psi}{2} s(1-s) \| \bfx - \bfy \|^2.
	\end{align*}
	We say $\psi$ is $\ell_\psi$-smooth for $\ell_\psi > 0$ if it is continuously differentiable and 
	\begin{equation*}
		\| \nabla \psi(\bfx) - \nabla \psi(\bfy) \| \leq \ell_\psi \| \bfx - \bfy \|, ~~\forall \bfx,\bfy \in \mbR^n.
	\end{equation*}
	\begin{lemma}[Descent Lemma in \cite{beck2017first}] \label{descent-lemma}
		Let $\psi: \mathbb{R}^n \to (- \infty, \infty]$ be a $\ell_\psi$-smooth function. Then we have
		\[
		\psi(\bfx) \leq \psi(\bfy) + \langle \nabla \psi(\bfy), \bfx - \bfy \rangle + \frac{\ell_\psi}{2} \| \bfx - \bfy \|^2, \quad \forall \ \bfx, \ \bfy \in \mathbb{R}^n.
		\]
	\end{lemma}
	
	For a function $\varphi:\mathbb{R}^m \to (-\infty, +\infty]$, its conjugate function $\varphi^*: \mathbb{R}^m  \to (-\infty, +\infty]$ is defined by
	\begin{align*}
		\varphi^*( \bfv ):= \sup_{ \bfz \in \mbR^m } \langle \bfv, \bfz  \rangle - \varphi(\bfz) .
	\end{align*}
	
	Some properties of conjugate function from \cite{beck2017first} are summarized as follows.
	
	\begin{lemma}\label{conjugate-property} Let $\varphi$ be a proper, l.s.c. and convex function, then we have
		
		(i) $\varphi^*$ is proper, l.s.c. and convex.
		
		(ii) $\bfu \in \partial \varphi(\bfz)$ if and only if $\bfz \in \partial \varphi^*(\bfu)$.
		
		(iii) If $\varphi$ is $\sigma_\varphi$-strongly convex, then $\varphi^*$ is $\frac{1}{\sigma_\varphi}$-smooth.
		
		(iv) If $\varphi$ is $\ell_\varphi$-smooth, then $\varphi^*$ is $\frac{1}{\ell_\varphi}$-strongly convex. 
	\end{lemma}
	
	\textbf{(B) P{\L}K functions.}
	They are also widely known as Kurdyka-{\L}ojasiewicz (K\L) functions defined
		by K{\L} inequality and popularized in
		\cite{attouch2010proximal,bolte2007clarke}.
	A recent paper \cite{bento2025convergence},
		which attributes the origin of K{\L} inequality to Polyak \cite{polyak1963gradient}, suggested the term Polyak-{\L}ojasiewicz-Kurdyka (P{\L}K) inequality to reflect
		the order of the time the inequality and its variants were first proposed. 
		We refer to the introduction part of \cite{bento2025convergence} for a short historical account of P{\L}K inequality. 
	
	Let $s \in (0,+\infty]$, we denote $\Psi_s$ as the class of all concave and continuous functions: $\psi:[0,s) \to \mbR_+$ which satisfies the following conditions
	\begin{itemize}
		\item[(i)] $\psi(0) = 0$;
		\item[(ii)] $\psi$ is continuously differentiable on $(0,s)$ and continuous at 0;
		\item[(iii)] $\psi'(t) > 0$ for all $t \in (0,s)$.
	\end{itemize}
	Next let us define Polyak-{\L}ojasiewicz-Kurdyka (P{\L}K) property. A proper and l.s.c. function $\varphi: \mbR^m \to (-\infty,+\infty]$ is said to have the P{\L}K property at $\obz \in \hbox{dom}\, \partial \varphi:= \{ \bfz \in \mbR^m: \partial \varphi(\bfz) \neq \varnothing \}$ if there exists $s \in (0,+\infty]$, a neighborhood $\mathcal{N}(\obz)$ of $\obz$ and a function $\psi \in \Psi_s$, such that the following inequality holds
	\begin{align} \label{KL-ineq}
		\psi'( \varphi(\bfz) - \varphi(\obz) ) \hbox{dist}( 0,\partial \varphi(\bfz) ) \geq 1,~~\forall \bfz \in \mathcal{N}(\obz) \cap [ \varphi(\obz) < \varphi(\bfz) < \varphi(\obz) + s ]
	\end{align}
	where for $\overline{c}_1,\overline{c}_2 > 0$, $[ \overline{c}_1 < \varphi(\bfz) < \overline{c}_2 ]: = \{ \bfz \in \mbR^m: \overline{c}_1 < \varphi(\bfz) < \overline{c}_2 \}$. If $\varphi$ satisfies the P{\L}K property at every point of $\hbox{dom } \partial \varphi$, then we call $\varphi$ a P{\L}K function.
	
	The P{\L}K property reveals that the values of $\varphi$ can be reparameterized to sharpness \cite{polyak1979sharp,polyak1987introduction} and this is crucial for global convergence analysis of algorithms in nonconvex settings. One attractive aspect of P{\L}K functions is that they are ubiquitous in applications. The common classes of P{\L}K functions include real subanalytic, semi-algebraic, convex functions with growth conditions and so on. For more illustrating examples about P{\L}K properties, we recommend \cite{attouch2010proximal} and \cite{bolte2014proximal}.
	
	\textbf{(C) Semismooth functions.} 
	Let $\Phi: \mbR^m \to \mbR^n$ be a locally Lipschitzian function. Suppose that $D_\Phi$ is the set containing all the differentiable points of $\Phi$ and $J\Phi(\bfz)$ is the Jacobian at a differentiable point $\bfz$. Then the generalized Jacobian of $\Phi$ is defined by
	\begin{align*}
		\partial_C \Phi (\obz) := {\rm co}\Big\{ 
		V ~|~\mbox{there exists a sequence $\bfz^k \in D_\Phi$ such that} \
		V = \lim_{\bfz^k \to \obz} J\Phi(\bfz^k) \Big\}.
	\end{align*} 
	Furthermore, $\Phi$ is said to be semismooth \cite{cui2021modern,facchinei2003finite,qi1993nonsmooth} if for any $\bfz \in \mbR^m$, $\Phi$ is directionally differentiable and it holds 
	\begin{align} \label{def-semismoothness}
		\Phi (\bfz + \bfd) - \Phi(\bfz) - H \bfd = o(\| \bfd \|) \quad 
		\forall \ H \in \partial_C \Phi(\bfz + \bfd), \ \bfd \in \mbR^m .
	\end{align}
	For a continuously differentiable function $\psi: \mbR^m \mapsto (-\infty, +\infty]$,
	we say it is $\mbox{SC}^1$ function if the gradient function $\nabla \psi(\cdot)$ is semismooth.

	\section{Optimality of Primal and Dual Problems} \label{pd-solution}
	
	In this section, we study the optimality conditions of the primal and dual problems. 
	The ultimate result is 
	the equivalence of local minimizers between the two problems in the sense that
	one can be obtained from the other. 
	This one-to-one correspondence can be seen as the duality theorem for the problem \eqref{01cop}. 
	For convenience, we consider the following equivalent reformulation of \eqref{01cop}
	\begin{equation} \label{01cop-c}
		\min_{\bfx, \bfu} f( \bfx ) + \indf(\bfu) , \quad s.t. \quad  
		(\bfx, \bfu) \in \F := \left\{ (\bfx, \bfu) \ | \  A\bfx + \bfb = \bfu \right\}.
	\end{equation}
	We interchangeably refer to them as the primal problem 
	depending on whether the variable $\bfu$ is needed or not. 
	
	\subsection{Equivalence of local minimizers and KKT points of the primal problem}
	
	We note that the strong convexity of $f$ ensures that the objective of \eqref{01cop} is l.s.c and coercive, then there must exist an optimal solution for the primal problem 
	\eqref{01cop} and \eqref{01cop-c}. 
	The KKT point of \eqref{01cop-c} with a continuously differentiable $f$ has been introduced in \cite{zhou2021quadratic}. 
	By utilizing limiting subdifferential, we can extend this definition to the case when $f$ is nonsmooth.
	
	\begin{definition}
		A point $\bfw^*:=(\bfx^*,\bfu^*) \in \mbR^{n+m}$ is a KKT point of \eqref{01cop-c} 
		if there exists a multiplier $\bfz^*$ such that the following system holds:
		\begin{equation} \label{01cop-KKT}
			\left\{
			\begin{array}{rl}
				\partial f(\bfx^*) + A^\top \bfz^* & \ni 0, \\
				- \bfz^* + \partial \indf(\bfu^*)  & \ni 0, \\
				A\bfx^* + \bfb - \bfu^*            & = 0.
			\end{array} 
			\right .
		\end{equation}
		We call $(\bfw^*,\bfz^*)$ a KKT pair of \eqref{01cop-c}.

	\end{definition}

	The convexity of $f$ and the piecewise linear structure of the composite term $\indf(A \bfx + \bfb)$
	allow us to characterize a local minimizer in terms of KKT point.
	
	\begin{proposition} \label{primal-opt}
		For problem \eqref{01cop-c}, $\bfw^*:=(\bfx^*,\bfu^*) \in \mbR^{n+m}$ is a local minimizer if and only if it is a KKT point.
	\end{proposition}
	
	\bp
	First, for a given reference point $\bfw^* =(\bfx^*, \bfu^*)$, 
	define the following linearly constrained convex programming:
	\begin{align} \label{cp}
		\min_{\bfx \in \mbR^n} f(\bfx) \quad \mbox{s.t.} \quad
		(\bfx, \bfu) \in \F^* := \left\{ (\bfx, \bfu)\ | \   A\bfx + \bfb - \bfu = 0, \ \ \bfu_{\J^*_{-}} \leq 0 \right\},
	\end{align}
	where $\J^*_{-}:= \{ i \in [m]: u^*_i \leq 0 \}$. 
	A minimizer $(\bfx,\bfu)$ of \eqref{cp} must satisfy the following system together with the multiplier $(\bfz, \bfz_{\J^*_-})$
	\begin{equation}   \label{cp-kkt}
		\left\{  \begin{aligned}
			& \partial f(\bfx) + A^\top \bfz \ni 0,\\
			& \bfz_{\overline{\J}^*_-} = 0, \  \bfz_{\J^*_-} \geq 0, \
			\bfu_{\J^*_-} \leq 0, \ 
			\langle \bfz_{\J^*_-}, \; \bfu_{{\J^*_-}}  \rangle = 0, \\
			& A \bfx + \bfb - \bfu = 0 .
		\end{aligned}  \right.
	\end{equation}
	Comparing \eqref{01cop-KKT} and \eqref{cp-kkt}, and utilizing 
	representation \eqref{partial01}, we see that
	a KKT point $\bfw^*$ of \eqref{01cop-c} must be a global minimizer of the convex problem \eqref{cp}. 
	With this equivalence in hand, we prove the claim in the proposition.

	\textbf{Necessity:}
	Let $ \bfw^*= (\bfx^*,\bfu^*)$ be a local minimizer of \eqref{01cop-c}, then there exists $\delta >0$ such that
	\begin{align*}
		f(\bfx) + \indf(\bfu) \geq f(\bfx^*) + \indf(\bfu^*), ~~~ \forall \bfw \in \N(\bfw^*,\delta) \cap \F.
	\end{align*}
	By the definition of $\J^*_-$, we have $\indf(\bfu) \leq \indf(\bfu^*)$ for any $\bfu \in \F^*$, and then we can obtain
	\begin{align*}
		f(\bfx)  \geq f(\bfx^*), ~~~ \forall \bfw \in \N(\bfw^*,\delta) \cap \F^*.
	\end{align*}
	This means $\bfw^*$ is also a local minimizer of \eqref{cp}. 
	Since Problem \eqref{cp} is convex, $\bfw^*$ is also a global minimizer of \eqref{cp}, which in
	turn is equivalent to a KKT point of \eqref{01cop-c}.

	\textbf{Sufficiency:}
	Let $ \bfw^*= (\bfx^*,\bfu^*)$ be a KKT point of \eqref{01cop-c}.
	We have proved at the beginning that $\bfw^*$ is also a global minimizer of \eqref{cp}.
	Thus we have
	\begin{align} \label{f>f*}
		f(\bfx) \geq f(\bfx^*), ~~\forall\; \bfw = (\bfx, \bfu) \in \F^* .
	\end{align}
	Now let us consider a sufficiently small radius $\delta^*$ such that for any $\bfw \in \N(\bfw^*,\delta^*)$, we have
	\begin{align} \label{relation}
		\{  i \in [m]: u_i > 0 \} \supseteq \overline{\J}^*_- 
		\quad \mbox{and} \quad
		f(\bfx) - f(\bfx^*) > - \lambda/2, 
	\end{align}
	where the second inequality follows from the lower semicontinuity of $f$. Then for any $\bfw \in \N(\bfw^*, \delta^*) \cap \F^* $, \eqref{relation} leads to $\indf(\bfu) \geq \indf(\bfu^*)$. This together with \eqref{f>f*} implies
	\begin{align*}
		f(\bfx) + \indf(\bfu) \geq f(\bfx^*) +  \indf(\bfu^*), ~~~ \forall \bfw \in \N(\bfw^*,\delta^*) \cap \F^*.
	\end{align*}
	If $\bfw \in \N(\bfw^*, \delta^*) \cap (\F \backslash \F^*)$, then there exists $i_0 \in \J^*_-$ such that $u_{i_0} > 0$. Combining this with \eqref{relation} leads to $\indf(\bfu) \geq \indf(\bfu^*) + \lambda $. Taking the inequality in \eqref{relation} into account, we have
	\begin{align*}
		f(\bfx) + \indf(\bfu) \geq f(\bfx^*) +  \indf(\bfu^*) + \lambda/2,  ~~~ \forall \bfw \in \N(\bfw^*,\delta^*) \cap \F^*.
	\end{align*}
	Hence, we proved that $\bfw^*$ is also a local minimizer of \eqref{01cop-c}.
	\ep

	\subsection{Equivalence of local minimizers and KKT points of the dual problem}
	
	We rewrite the dual problem \eqref{01dco-max} as a minimization problem:
	\begin{align} \label{01dco-uc}
		\min_{\bfz \in \mbR^m} F(\bfz):= h(\bfz) + g(\bfz),
	\end{align}
	Since $f$ is $\sigma_f$-strongly convex, it follows from Lem.~\ref{conjugate-property} that $h$ is convex and $\ell_h$-smooth with $\ell_h = \| A \|^2/\sigma_f$.
	Utilizing the definition of $\partial g$ and $\proxg(\cdot)$, we introduce the KKT and P-stationary points of \eqref{01dco-uc}. These two types of stationary points have also been given in \cite{kanzow2021globalized,mordukhovich2023globally} for minimizing the sum of two functions.
	
	\begin{definition} \label{Def-KKT}
		Let us consider a reference point $\bfz^* \in \mbR^m$.
		
		(i) $\bfz^*$ is said to be a KKT point of \eqref{01dco-uc} if 
		\begin{equation} \label{01dco-kkt}
			0 \in \nabla h(\bfz^*) + \partial g(\bfz^*) .
		\end{equation}
		
		(ii) $\bfz^*$ is called a P-stationary point of \eqref{01dco-uc} if there is $\tau > 0$ such that
		\begin{equation} \label{01dco-pstat}
			\bfz^* \in \proxg( \bfz^* - \tau \nabla h(\bfz^*) ) .
		\end{equation}
	\end{definition}
	
	In comparison to the definition of KKT point, a P-stationary point involves 
	a positive parameter $\tau>0$.
	For a given $\bfz^*$, let $\bfu^* = -\nabla h( \bfz^* )$.
	The conditions \eqref{01dco-kkt} and \eqref{01dco-pstat} can be rewritten respectively as
	\begin{align*}
		\bfu^* \in \partial g(\bfz^*)\quad \mbox{and} \quad 	\bfz^* \in \proxg( \bfz^* + \tau \bfu^* ).
	\end{align*}
	It then follows from Lem.~\ref{property-g}(v) that a KKT point of \eqref{01dco-uc} is also its
	P-stationary point as stated in the result below.
	
	\begin{lemma} \label{Lemma-KKT-P-Dual}
		If $\bfz^*$ is a P-stationary point of \eqref{01dco-uc}, then it is also a KKT point of \eqref{01dco-uc}. Conversely, if $\bfz^*$ is a KKT point of \eqref{01dco-uc}, then it is a P-stationary point with $\tau\in (0,\min\{\bfta_1(\bfz^*),\bfta_2(-\nabla h(\bfz^*))\})$, where $\bfta_1$ and $\bfta_2$ have been defined in Lem.~\ref{property-g}.
	\end{lemma}
	
	Moreover, every P-stationary point (KKT point) of \eqref{01dco-uc} is also its local minimizer, as proved
	below.
	
	\begin{proposition} \label{Prop-LocalMinimizer-P-Dual}
		For problem \eqref{01dco-uc}, $\bfz^* $ is a local minimizer if and only if it is a P-stationary point. 
	\end{proposition} 
	
	\bp
	If $\bfz^*$ is a local minimizer of \eqref{01dco-uc}, then using \cite[Thm.~10.1, Exercise 10.10]{RockWets98}, we can derive \eqref{01dco-kkt}, which means $\bfz^*$ is a P-stationary point of \eqref{01dco-uc} as we have claimed in Lem.~\ref{Lemma-KKT-P-Dual}. 
	
	Now let us prove the sufficiency. 
	Suppose $\bfz^*$ is a P-stationary point of \eqref{01dco-uc}.
	We must have $\bfz^* \in \mbR^m_+$ due to the function $g(\bfz)$.
	We choose $\widetilde{\delta}^*$ sufficiently small such that
	\begin{align} \label{relation1}
		\mathcal{S}(\bfz) \supseteq \mathcal{S}(\bfz^*):= \S^* \quad \mbox{and} \quad
		| h(\bfz) - h(\bfz^*) | < \mu/2, ~~\forall\; \bfz \in \N(\bfz^*, \widetilde{\delta}^*),
	\end{align}
	where we used the inclusion property of the support-set function $\mathcal{S}(\cdot)$ near $\bfz^*$
	and the inequality used the continuity of $h$. 
	Let $\Theta^*: = \{ \bfz \in \mathbb{R}^m: \S(\bfz) = \S(\bfz^*) \}$.
	We consider two cases.
	
	Case I: $\bfz \in \N(\bfz^*, \widetilde{\delta}^*) \cap \mbR^m_+ \cap \Theta^*$.
	For this case, $g(\bfz) = g(\bfz^*)$. Denote $\bfu^* = - \nabla h(\bfz^*)$, by the convexity of $h$, we obtain the following inequality:
	\begin{align*}
		h(\bfz) - h(\bfz^*) \geq -\langle \bfu^*, \bfz - \bfz^*  \rangle = -\langle \bfu^*_{\S^*}, (\bfz - \bfz^*)_{\S^*}  \rangle = 0, 
	\end{align*}
	where the last two equations follows from $\bfz \in \Theta^*$ and Lem.~\ref{property-g} (iv) respectively. 
	
	Case II: $\bfz \in \N(\bfz^*, \widetilde{\delta}^*) \cap \mbR^m_+ \backslash \Theta^*$.
	Then \eqref{relation1} implies that $\S(\bfz)$ contains at least one more element than 
	$\S(\bfz^*)$. This in turn implies $g(\bfz) \geq g(\bfz^*) + \mu$. 
	It follows from \eqref{relation1} that $h(\bfz) \ge h(\bfz^*) - \mu/2$. 
	The two inequalities just obtained means
	$h(\bfz) + g(\bfz) \geq h(\bfz^*) + g(\bfz^*) + \mu/2$.
	
	In summary, we have proved the following:
	\be \label{Bound-F}
	\left\{
	\begin{array}{ll}
		g(\bfz) = g(\bfz^*), \ h(\bfz) \ge h(\bfz^*), \  & \mbox{for}\ \bfz \in \N(\bfz^*, \widetilde{\delta}^*) \cap \mbR^m_+ \cap \Theta^* , \\
		h(\bfz) + g(\bfz) \geq h(\bfz^*) + g(\bfz^*) + \mu/2, \ & \mbox{for}\
		\bfz \in \N(\bfz^*, \widetilde{\delta}^*) \cap \mbR^m_+ \backslash \Theta^* ,\\
		F(\bfz) = +\infty, \  & \mbox{for}\ \bfz \not\in \mbR^m_+.
	\end{array} 
	\right .
	\ee 
	Consequently, we have $F(\bfz) \ge F(\bfz^*)$ for $\bfz \in \N(\bfz^*, \widetilde{\delta}^*)$
	and $\bfz^*$ is a local minimizer of \eqref{01dco-uc}.
	\ep
	
	We finish this subsection with a result on the relationship between global minimizers and {\rm P}-stationary points of \eqref{01dco-uc}.
		
		\begin{proposition} \label{Prop-global-P-stationary}
			For problem \eqref{01dco-uc}, the following assertions hold:
			
			(i) If $\bfz^*$ is a global minimizer, then it is a {\rm P}-stationary point with $\tau \in (0,1/\ell_h)$.
			
			(ii) If $h$ is $\sigma_h$-strongly convex and $\bfz^*$ is a {\rm P}-stationary point with some $\tau > 1/\sigma_h$, then it is the unique global minimizer.
		\end{proposition}
		
		\bp
		(i) directly follows from \cite[Prop.~4.4]{khanh2024coderivative}. 
		Now we prove (ii). By the definitions of {\rm P}-stationary point and proximal point, we have
		\begin{align*}
			\frac{1}{2\tau} \| \bfz - (\bfz^* - \tau \nabla h(\bfz^*)) \|^2 + g(\bfz) \geq \frac{1}{2\tau} \| \tau \nabla h(\bfz^*) \|^2 + g(\bfz^*).
		\end{align*}
		After simplification, we obtain
		\begin{align*}
			g(\bfz ) - g(\bfz^*) \geq -\frac{1}{2\tau} \| \bfz - \bfz^* \|^2 - \langle \nabla h(\bfz^*), \bfz - \bfz^* \rangle.
		\end{align*}
		The $\sigma_h$-strong convexity of $h$ implies
		\begin{align*}
			h(\bfz) - h(\bfz^*) \geq \langle \nabla h(\bfz^*), \bfz - \bfz^* \rangle + \frac{\sigma_h}{2} \| \bfz - \bfz^* \|^2.
		\end{align*}
		Adding the two inequalities above leads to
		\begin{align*}
			F(\bfz) - F(\bfz^*) \geq \Big(\frac{\sigma_h}{2} - \frac{1}{2\tau} \Big) \| \bfz - \bfz^* \|^2 > 0, \ \mbox{whenever} \ \bfz \neq \bfz^*.
		\end{align*}
		which means that $\bfz^*$ is the unique global minimizer of \eqref{01dco-uc}.
		\ep
	
	\subsection{Equivalence of local minimizers of primal and dual problems} 
	
	The equivalence results in the two preceding subsections lead to the main 
	result in this section. 
	
	\begin{theorem} [Equivalence Characterization of Local Minimizer] \label{foc}
		If $\bfw^* = (\bfx^*, \bfu^*)$ is a local minimizer of the primal problem \eqref{01cop-c} with associated multiplier $\bfz^*$, then $\bfz^*$ is a local minimizer of the dual problem \eqref{01dco-uc}.  
		Conversely, if $\bfz^*$ is a local minimizer of the dual problem \eqref{01dco-uc}, 
		then $\bfw^* := (\nabla f^*( -A^\top \bfz^* ),-\nabla h(\bfz^*))$ is a local minimizer of
		the primal problem \eqref{01cop-c}.
	\end{theorem}
	
	\bp
	In Prop.~\ref{primal-opt} for the primal problem \eqref{01cop-c}, we established the equivalence of
	its local minimizer and KKT point.
	Combining Lem.~\ref{Lemma-KKT-P-Dual} and Prop.~\ref{Prop-LocalMinimizer-P-Dual}
	for the dual problem \eqref{01dco-uc},  we established the equivalence of its local minimizer and its
	KKT point. We only need to prove the equivalence of the KKT point for both primal and the dual 
	problems.
	
	Suppose $\bfz^*$ is a  KKT point of the dual problem \eqref{01dco-uc}.
	Let $\bfx^* := \nabla f^*( -A^\top \bfz^* )$ and $\bfu^* := -\nabla h (\bfz^*)$. 
	Then by using Lem.~\ref{conjugate-property}(ii) and $h(\bfz) = f^*(-A^\top \bfz) - \langle \bfb, \bfz \rangle$, we know that the KKT condition \eqref{01dco-kkt} for the dual problem is equivalent to  
	\be \label{01cop-kkt1}
	\left\{
	\begin{array}{l}
		\partial f(\bfx^*) + A^\top \bfz^* \ni 0, \\
		- \bfu^* + \partial g(\bfz^*) \ni 0, \\
		A\bfx^* + \bfb - \bfu^* = 0 .
	\end{array} 
	\right .
	\ee 
	Comparing \eqref{01cop-kkt1} with \eqref{01cop-KKT}, we just need to prove their second lines are equivalent. This directly follows from Lem.~\ref{property-g}(iii).
	Consequently, $\bfw^*=(\bfx^*, \bfu^*)$ is a KKT point of the primal problem.
	This establishes the equivalence result.
	\ep

	At the end of this section, we introduce a second-order sufficient condition (SOSC) for $\eqref{01dco-uc}$. Just like the case of classic smooth nonlinear programming (see e.g. \cite{nocedal2006numerical}), we will show that our SOSC is useful to derive a quadratic growth condition for a P-stationary point and local superlinear convergence rate for a second-order algorithm. 
	
	Given a P-stationary point $\bfz^*$, denote the support set $\S^*:= \S(\bfz^*)$. One can verify that it is also a KKT point with multiplier $(\nabla_{\S^*} h(\bfz^*), \nabla_{\oS^*} h(\bfz^*))$ of the following convex programming
	\begin{align}\label{cp-d}
		\min_{\bfz \in \mbR^m} h(\bfz),  \hbox{   s.t.  } \bfz_{\S^*} \geq 0, \hbox{ } \bfz_{\oS^*} = 0.
	\end{align}
	The critical cone of the above problem at 
	$\bfz^*$ is $\C:=\{ \bfd \in \mbR^m: \bfd_{\oS^*} = 0 \}$ (it happens to be a subspace). Since $h$ is $\ell_h$-smooth, for any $\bfz \in \mbR^m$,  $\partial^2 h(\bfz):= \partial_C( \nabla h(\bfz) ) \neq \varnothing.$ Then we can give the following definition based on the second-order sufficient condition for the convex problem \eqref{cp-d} (see e.g. \cite{facchinei2003finite}).

	\begin{definition} \label{def-sosc}
		Given a P-stationary point $\bfz^* \in \mbR^m$ of \eqref{01dco-uc}, we say the second-order sufficient condition (SOSC) holds at $\bfz^*$ if 
		\begin{align} \label{sosc}
			\bfd^\top H^*_{\S^*} \bfd > 0,~~ \forall H^* \in \partial^2 h (\bfz^*) \hbox{ and } \bfd \in \mbR^{| \S^* |} \backslash \{ 0 \}.
		\end{align}
	\end{definition}
	
	Similar to the results for $\mbox{SC}^1$ functions (see \cite[Prop.~7.4.12]{facchinei2003finite}), we can prove the quadratic-growth condition for \eqref{01dco-uc} at a P-stationary point $\bfz^*$ satisfying SOSC.
	
	\begin{proposition}[Quadratic Growth Condition]\label{Prop-QuadraticGrowth}
		If $\nabla h$ is semismooth and $\bfz^* \in \mbR^m$ is a P-stationary point of \eqref{01dco-uc} satisfying SOSC, then there exist $\overline{\delta}^*, c^* > 0$ such that 
		the following quadratic growth at $\bfz^*$ holds:
		\begin{align}
			F(\bfz) \geq F(\bfz^*) + c^* \| \bfz - \bfz^* \|^2, \qquad
			\forall \; \bfz \in \N(\bfz^*,\overline{\delta}^*) \cap \mathbb{R}^m_+ .
		\end{align} 
	\end{proposition}
	\bp
	As we have mentioned before Def. \ref{def-sosc}, $\bfz^*$ is also a KKT point of \eqref{cp-d}. Considering that \eqref{sosc} holds, it follows from \cite[Prop.~7.4.12]{facchinei2003finite} that there exists $\overline{\delta}_1^*$ and $c^* > 0$ such that
	\begin{align} \label{quad-gro-cp}
		h(\bfz) \geq h(\bfz^*) + c^*\| \bfz - \bfz^* \|^2,~ \forall \bfz \in \N( \bfz^*, \overline{\delta}_1^* ) \cap \{ \bfz \in \mbR^m: \bfz_{\S^*} \geq 0, \bfz_{\oS^*} = 0  \}.
	\end{align}
	Then we can take $\ode^* =\min\{ \widetilde{\delta}^*, \ode^*_1, \sqrt{\mu/(2c^*)} \}$, where $\widetilde{\delta}^*$ has been defined in the proof of Prop.~\ref{Prop-LocalMinimizer-P-Dual}. 
	Then \eqref{relation1} holds.
	Furthermore, for any $\bfz \in \N( \bfz^*, \overline{\delta}^* )$, we have $\{ \bfz \in \mbR^m: \bfz_{\S^*} \geq 0, \hbox{ } \bfz_{\oS^*} = 0  \} =  \mbR^m_+ \cap \Theta^*$, where $\Theta^* = \{ \bfz \in \mathbb{R}^m: \S(\bfz) = \S(\bfz^*) \}$. 
	Combination of \eqref{quad-gro-cp} and \eqref{Bound-F} leads to 
	\begin{align*}
		&	F(\bfz) - F(\bfz^*) \geq c^* \| \bfz - \bfz^* \|^2,~~\forall \bfz \in \N(\bfz^*, \ode^*) \cap \mbR^m_+ \cap \Theta^* , \\
		&F(\bfz) - F(\bfz^*) \geq \mu/2 \geq c^* \| \bfz - \bfz^* \|^2,~~\forall \bfz \in \N(\bfz^*, \ode^*) \cap \mbR^m_+ \backslash \Theta^*.
	\end{align*} 
	We arrived at the claimed result.
	\ep

	\section{A Subspace Gradient Semismooth Newton Method} \label{sec_sgsn}
	Since we have established the equivalence of solutions between the primal problem \eqref{01cop-c} and the dual problem\eqref{01dco-uc}, 
	we now develop an efficient method for the dual problem.
	It is essentially a subspace method, which first identifies a subspace and then within this subspace applies a 
	gradient step followed by a semismooth Newton step if certain conditions are met.
	We call the resulting method a subspace gradient semismooth Newton (SGSN) method.
	The key challenge is to establish its global convergence as well as its local superlinear 
	convergence rate 
	under mild assumptions. We now describe the framework for SGSN followed by its convergence
	analysis.
	
	\subsection{Steps of SGSN}
	
	Suppose $\bfz^k$ is the current iterate. 
	First, we perform the classical proximal gradient step
	\begin{align} \label{pg-step1}
		\bfv^k \in \proxg ( \bfz^k - \tau \nabla h(\bfz^k) ).
	\end{align}
	The aim of this step is to ensure global convergence and to 
	identify the subspace $\mbR^{|T_k|}$ indexed by 
	\[ 
	T_{k}:= \{ i \in [m]: [ \bfz^k - \tau \nabla h(\bfz^k) ]_i > \sqrt{2\tau\mu}  \}.
	\] 
	From \eqref{proxg}, 
	once such $\bfv^k$ is given by
	\begin{align} \label{pg-step}
		\bfv^{k}_{T_k} = [ \bfz^k - \tau \nabla h(\bfz^k) ]_{T_k} 
		\quad  \mbox{and }  \quad \bfv^{k}_{\oT_k} = 0.
	\end{align}
	In other words, $\bfv^k$ takes a gradient descent on $\mbR^{|T_k|}$ while setting the
	complementary part to zero.
	
	The second step is to compute a Newton direction along the subspace identified in the first step.
	Specifically, we minimize a second-order Taylor approximation of $h( \bfv^k + \bfd )$ with $\bfd_{\oT_k} = 0$:
	\begin{align*}
		\bfd^k := \arg\min_{ \bfd \in \mbR^m }~ \langle \nabla h( \bfv^k ), \bfd \rangle + \frac{1}{2} \bfd^\top ( H^k + \gamma_k I)  \bfd, \hbox{   \ \ s.t. } \bfd_{\oT_k} = 0.
	\end{align*}
	where $\gamma_k>0$ is a regularization parameter, $H^k \in \whp^2 h(\bfv^k)$, and the set $\whp^2 h (\bfz)$ is defined by
	\begin{align} \label{hat-partial}
		\whp^2 h (\bfz): = \{ A Q A^\top \ | \  Q \in \partial^2 f^*(-A^\top \bfz)   \} .
	\end{align}
	Since $f^*$ is $(1/\sigma_f)$-smooth, $\whp^2 h (\bfz)$ is always nonempty. 
	The reason to use $\whp^2 h$ instead of $\partial^2 h$ is 
	because $h$ includes the composition term $f^*\circ(-A)$ and it is usually difficult to 
	characterize the exact form of $\partial^2 h$. 
	The relationship between the two sets is given in \cite[Thm.~2.6.6]{clarke1990optimization}:  
	\begin{align} \label{chain-rule}
		\partial^2 h(\bfz) \overline{\bfd} \subseteq \whp^2 h(\bfz) \overline{\bfd}, \qquad
		\forall \; \overline{\bfd} \in \mbR^m.
	\end{align}
	Since $H^k$ is positive semidefinite, $( H^k + \gamma_k I)$ is positive definite and
	hence  $\bfd^k$ is well-defined.
	
	The final step is to measure whether the Newton direction is ``sufficiently'' good. Let
	\[
	\wz^{k+1} := \bfv^k + \alpha_k \bfd^k,
	\]
	where $\alpha_k >0$ can be explicitly computed.
	We accept $\wz^{k+1}$ as our next iterate if it satisfies the following two conditions:
	\begin{align}
		& F(\bfv^k) - F(\wz^{k+1}) \geq c_1 \|  \wz^{k+1} - \bfv^k \|^2  \tag{\textbf{C1}} \label{c1} \\
		& \| \nabla_{T_k} h(\wz^{k+1}) \| \leq c_2 \| \wz^{k+1} - \bfv^k \|, \tag{\textbf{C2}} \label{c2}
	\end{align}
	where $c_1,c_2 >0$.
	Otherwise, we take $\bfv^k$ as our next iterate.
	Those steps are summarized in Alg.~\ref{SGSN}.
	
	\begin{algorithm}[htbp] 
		\caption{SGSN: Subspace Gradient Semismooth Newton Method} \label{SGSN}
		\begin{algorithmic}
			
			\STATE{Initialization: Take $\tau \in (0,1/\ell_h)$, $c_1, c_2, \gamma > 0$, and then start with $\bfz^0 \in  \mathbb{R}^m$.}
			\FOR{$k=0,1,\cdots$}
			
			\STATE{\textbf{1. Proximal gradient: } }	
			Compute $T_k$, and then update the gradient iterate $\bfv^{k}$ by \eqref{pg-step}.

			\STATE{\textbf{2. Semismooth Newton:} } Denote $\gamma_{k}:= \gamma \| \nabla_{T_k} h(\bfv^{k}) \|$, and choose $H^{k} \in \whp^2 h(\bfv^{k})$. Compute Newton direction $\bfd^{k}_{T_k}$ by solving
			\begin{align} \label{Nt-eq}
				( H^{k}_{T_k} + \gamma_{k} I )\bfd^{k}_{T_k} = - \nabla_{T_k} h(\bfv^{k}) .
			\end{align}
			Update Newton iterate by
			\begin{align} \label{newton-step}
				\wz^{k+1}_{T_k} = \bfv^{k}_{T_k} + \alpha_{k} \bfd^{k}_{T_k} \hbox{ and } \wz^{k+1}_{\oT_k} = 0.
			\end{align}
			where $\alpha_k = \min\{ 1, \beta_k \}$ with
			\begin{align} \label{step-beta}
				\beta_k := \left\{ \begin{aligned}
					& 1, &&\mbox{if } \bfd^k_{T_k} \geq 0 \\
					& \min\{ -v^k_i/d^k_i : i \in T_k \mbox{ and } d^k_i < 0 \}, &&\mbox{otherwise}.
				\end{aligned} \right.
			\end{align}
			
			\STATE{\textbf{3. Update step:} }	Update $\bfz^k$ either by the semismooth Newton step or
			the proximal gradient step as follows:
			\be \label{Newton-Condition}
			\bfz^{k+1} = \left\{
			\begin{array}{ll}
				\wz^{k+1} , & \ \hbox{if } \cone \hbox{ and } \ctwo \hbox{ hold} \\ [1ex]
				\bfv^{k}, & \ \mbox{otherwise.}
			\end{array} 
			\right .
			\ee 	
			
			\ENDFOR
		\end{algorithmic}
	\end{algorithm}
	
	\begin{remark} \label{rem-sgsn}
		(i) From \eqref{pg-step}, one can observe that $T_k$ exactly includes all nonzero entries of $\bfv^k$. Together with \eqref{newton-step}, we have
		\begin{align} \label{Tk-support}
			T_k = \S(\bfv_k) \supseteq \S(\wz^{k+1}) .
		\end{align}
		%
		
		(ii) Denote $\bfx^k := \nabla f^*( -A^\top \bfz^k )$, then $-\nabla h(\bfz^k) = A \bfx^k + \bfb$ and $\bfv^{k+1}_{T_k} = [\bfz^k + \tau ( A \bfx^k + \bfb )]_{T_k}$. 
		Following this setting,
		the computational cost of the proximal gradient step is $\O(|T_k|n)$. 
		Regarding the Newton step, let us consider a simple case when $f$ is twice continuously differentiable and so is $f^*$. Then the Hessian is $H^k = A Q^k A^\top$ with $Q^k:= \nabla^2 f^*(-A^\top \bfz^k)$ and Newton equation \eqref{Nt-eq} will be 
		\begin{align*} 
			(  \gamma_{k} I + A_{T_k:} Q^k A^\top_{T_k:} )\bfd^{k}_{T_k} = - \nabla_{T_k} h(\bfv^{k}).
		\end{align*} 
		In many applications, including the AUC maximization and multi-label classification, $f$ is often a separable regularizer. Thereby, $Q^k$ is diagonal. Then the above equation can be efficiently solved by linear conjugate gradient method and the computational cost is $\O(n|T_k|^2)$. 
		In the case $n < |T_k|$ and the inverse of $Q^k$ is cheap to compute, 
		we can make use of the Sherman-Morrison-Woodbury formula so that the overall computational cost is in the order of $\O(n^2|T_k|)$.

		(iii) At the $k$-th iteration, if it happens that $\bfv^k = \bfz^k$, then $\nabla_{T_k} h(\bfv^k) = 0$ and $\bfv^k_{\oT_k} = 0$, which further implies $0 \in \nabla h(\bfv^k) + \partial g(\bfv^k)$ and $\widetilde{\bfz}^{k+1} = \bfv^k$. 
		This means that $\bfz^k$ is a KKT point (see \eqref{01dco-kkt} in Def.~\ref{Def-KKT}). 
		Consequently, after the $k$-th iteration, the iterate would not change  and is
		a P-stationary point of \eqref{01dco-uc} by Lem.~\ref{Lemma-KKT-P-Dual}. 
	\end{remark}

	\subsection{Global Convergence} \label{global-convergence}
	
	In this subsection, our main goal is to prove the sequence $\{ \bfz^k \}_{k \in \mbN}$ generated by SGSN converges to a P-stationary point of \eqref{01dco-uc}.
	We first state some basic properties of the generated sequence.
	%
	%
	
	\begin{proposition} \label{Prop-SGSN-property}
		Let $\{ (\bfz^k, \bfv^k) \}_{k \in \mbN }$ be the sequence generated by SGSN. 
		The following holds.
		
		\begin{itemize}
			\item[(i)] The sequence of objective function value satisfies sufficient descent 
			\begin{align} \label{suffcient-des}
				F(\bfz^k) - F(\bfv^{k}) \geq \eta_0 \| \bfv^{k} - \bfz^k \|^2,~~F(\bfv^k) - F(\bfz^{k+1}) \geq c_1 \| \bfv^{k} - \bfz^{k+1} \|^2,
			\end{align}
			where $\eta_0 = 1/(2\tau) - \ell_h/2 > 0$.

			\item[(ii)] There exist $\bfp^k \in \partial F(\bfv^k)$ such that $\| \bfp^{k} \| \leq (1/\tau + \ell_h ) \| \bfv^k - \bfz^k \|$. If Newton iterate $\wz^{k+1}$ is accepted, then there exists $\wq^{k+1} \in \partial F(\wz^{k+1})$ such that $\| \wq^{k+1} \| \leq \eta_1  \| \bfz^{k+1} - \bfv^k \|$,
			where $\eta_1:= \max\{ c_2, 1/\tau + \ell_h \}$. 
			
			
			\item[(iii)] 
			Suppose the function sequence 
			$\{ F(\bfz^k)\}_{k \in \mbN }$ is bounded from below,
			then there exists $C \in \mbR$ such that 
			\begin{equation}\label{lim-F}
				\lim_{k \to \infty} F(\bfz^k) = \lim_{k \to \infty} F(\bfv^k) = C.
			\end{equation}
			Furthermore, it holds
			\begin{align} \label{successive-change}
				\lim_{k \to \infty} \| \bfv^{k} - \bfz^k \| = 0 \quad \hbox{and}\quad \lim_{k \to \infty} \| \bfv^{k} - \bfz^{k+1} \| = 0 .
			\end{align}
			
			\item[(iv)] Let $\bfz^*$ be an accumulation point of $\{ \bfz^k \}_{k \in \mbN}$, then $\bfz^*$ is a P-stationary point of \eqref{01dco-uc} and $F(\bfz^*) = C$.
			The corresponding conclusions for sequence $\{ \bfv^k \}_{k \in \mbN}$ also hold. 
			
		\end{itemize} 
	\end{proposition}
	
	\bp
	(i) By \eqref{pg-step1} and the definition of the proximal operator, we have
	\begin{align*}
		\| \bfv^k - ( \bfz^k - \tau\nabla h (\bfz^k) ) \|^2/2 + \tau g(\bfv^k) \leq  \| \bfz^k - ( \bfz^k - \tau\nabla h (\bfz^k) ) \|^2/2 + \tau g(\bfz^k),
	\end{align*}
	which implies
	\begin{align*}
		\| \bfv^k - \bfz^k \|^2/(2\tau) + \langle \nabla h(\bfz^k), \bfv^k - \bfz^k \rangle + g(\bfv^k) \leq g(\bfz^k) .
	\end{align*}
	Since $h$ is $\ell_h$-smooth, using Lem.~\ref{descent-lemma}, we  obtain
	\begin{align*}
		h(\bfv^k) \leq h(\bfz^k) + \langle \nabla h(\bfz^k), \bfv^k - \bfz^k \rangle + (\ell_h/2) \| \bfv^k - \bfz^k \|^2 .
	\end{align*}
	Adding the above two inequalities leads to $F(\bfz^k) - F(\bfv^{k}) \geq \eta_0 \| \bfv^{k} - \bfz^k \|^2$. If $\bfz^{k+1} = \wz^{k+1}$, then \eqref{suffcient-des} follows from \eqref{c1}, and \eqref{suffcient-des} automatically holds otherwise.
	
	(ii) From \eqref{pg-step1} and the definition of the proximal operator, we have
	\begin{align*}
		\bfv^k \in \arg\min_{\bfz} \| \bfz - (\bfz^k -  \tau\nabla h(\bfz^k)) \|^2/2 
		+ \tau g(\bfz) .
	\end{align*}
	By \cite[Exercise 8.8]{RockWets98}, we have 
	\begin{align*}
		0 \in \bfv^k - (\bfz^k -  \tau\nabla h(\bfz^k)) + \tau \partial  g(\bfv^k).
	\end{align*}
	Denoting $\bfp^k:= -(\bfv^k - \bfz^k)/\tau + \nabla h(\bfv^k) - \nabla h(\bfz^k)$, then $\bfp^k \in \nabla h(\bfv^k) + \partial g(\bfv^k) = \partial F(\bfv^k)$. We can also estimate $\| \bfp^{k} \| \leq (1/\tau + \ell_h ) \| \bfv^k - \bfz^k \|$ by using $\ell_h$-smoothness of $h$. 
	
	If $\bfz^{k+1} = \wz^{k+1}$, let us take $\wq^{k+1} = (\nabla_{T_k} h(\wz^{k+1});0_{\oT_{k}})$. Noticing that $\wz^{k+1}_{\oT_k} = 0$ from \eqref{newton-step}, we have $( 0_{T_k}; -\nabla_{\oT_k} h(\wz^{k+1}) ) \in \partial g (\wz^{k+1})$. This means that $\wq^{k+1} \in  \nabla h(\wz^{k+1}) + \partial g(\wz^{k+1}) = \partial F(\wz^{k+1})$. In this case, the upper bounded of $\| \wq^{k+1} \|$ follows from \eqref{c2}. Overall, we finish the proof of (ii).
	
	(iii) 
	Since $\{ F(\bfz^k) \}_{k \in \mbN}$ is assumed to be bounded below, \eqref{suffcient-des} implies $ F(\bfz^{k+1}) \leq F(\bfv^k) \leq F(\bfz^k) $, and thus $\{ F(\bfv^k) \}_{k \in \mbN}$ is also bounded below. Then it follows from the monotone convergence theorem that \eqref{lim-F} holds. Using \eqref{suffcient-des} and \eqref{lim-F}, we can derive \eqref{successive-change}. 
	
	(iv) 
	Since $\bfz^*$ is an accumulation point, there exists a subsequence $\{ \bfz^{k_j} \}_{j \in \mbN}$ such that $\lim_{j \to \infty} \bfz^{k_j} = \bfz^*$. Then we have $\lim_{j \to \infty} \bfz^{k_j} - \tau \nabla h(\bfz^{k_j}) = \bfz^* - \tau \nabla h(\bfz^*)$, as well as $\lim_{j \to \infty} \bfv^{k_j} = \bfz^*$ from \eqref{successive-change}. It follows from 
	the Proximal Behavior theorem \cite[Thm.~1.25]{RockWets98} that $\bfz^*$ is a P-stationary point of \eqref{01dco-uc}. 
	
	Now let us prove $F(\bfz^*) = C$. From \eqref{pg-step1} and the definition of proximal operator, we have
	\begin{align*}
		\| \bfv^{k_j} - (\bfz^{k_j} - \tau \nabla h(\bfz^{k_j})) \|^2/2 + \tau g(\bfv^{k_j}) \leq \| \bfz^* - (\bfz^{k_j} - \tau \nabla h(\bfz^{k_j})) \|^2/2 + \tau g(\bfz^*) .
	\end{align*}
	Meanwhile, the second inequality of \eqref{suffcient-des} implies
	\begin{align*}
		h(\bfz^{k_j+1}) + g(\bfz^{k_j+1}) + c_1 \| \bfv^{k_j} - \bfz^{k_j+1} \|^2 \leq h(\bfv^{k_j}) + g(\bfv^{k_j}).
	\end{align*}
	Considering that $\lim_{j \to \infty} \bfz^{k_j} = \lim_{j \to \infty} \bfz^{k_j+1} =\lim_{j \to \infty} \bfv^{k_j} = \bfz^*$ holds, taking limit superior on above two inequalities implies
	$\limsup_{j \to \infty} g (\bfz^{k_j}) \leq \limsup_{j \to \infty} g (\bfv^{k_j}) \leq g(\bfz^*)$. This together with lower semi-continuity of $g$ means $\lim_{j \to \infty} g (\bfz^{k_j}) = g(\bfz^*)$. Considering the continuity of $h$, we have $F(\bfz^*) = \lim_{j \to \infty} F(\bfz^{k_j}) = C$. By a similar procedure, we can prove the counterpart conclusion for sequence $\{ \bfv^k \}_{k \in \mbN}$.
	\ep
	
	\begin{remark} \label{Remark-Sequence-Boundedness-Assumption}
		The assumption of the boundedness of the function sequence $\{ F(\bfz^k)\}_{k \in \mbN }$
		and the existence of an accumulation point of the sequence $\{\bfz^k\}_{k \in \mbN }$ can be 
		ensured by the following assumption:
		\begin{assumption}\label{seq-bound}
			The sequence $\{\bfz^k\}_{k \in \mbN}$ generated by SGSN is bounded.
		\end{assumption}
		This sequence boundedness assumption is commonly used for convergence analysis of algorithms in nonconvex settings (see e.g. \cite{bolte2014proximal,bolte2018nonconvex,chen2017augmented,hallak2023adaptive}).
		Under Assumption~\ref{seq-bound}, the generated sequence $\{ \bfv^k\}_{k \in \mbN}$ is also
		bounded due to \eqref{successive-change}.
	\end{remark}

	Prop.~\ref{Prop-SGSN-property} presents some basic properties about the sequence
	$\{ (\bfz^k, \bfv^k)\}_{k \in \mbN }$ in separate cases. 
	When Condition \eqref{c1} or \eqref{c2} is not satisfied at the $k$-th iteration, we would have
	$\bfz^{k+1} = \bfv^k$ (repeating points). 
	We now remove those repeating points from $\{ (\bfz^k, \bfv^k)\}_{k \in \mbN }$
	and relabel the sequence by $\bfy^j$:
	\begin{align*}
		\{ \bfy^j \}_{j \in \mbN}:= \{ \bfz^0, \bfv^0, \cdots, \bfz^k, \bfv^k, \cdots   \} \backslash \{ \bfz^{k+1}: \bfz^{k+1} = \bfv^{k}, k \in \mathbb{N}  \}.
	\end{align*}
	We can see that  $\{ \bfz^k \}_{k \in \mbN}$ is a subsequence of $\{ \bfy^j \}_{j \in \mbN}$, because if $\bfz^{k+1} = \bfv^k$, then $\bfv^k \in \{ \bfy^j \}_{j \in \mbN}$, otherwise  $\bfz^{k+1} \in \{ \bfy^j \}_{j \in \mbN}$. We will show that $\{ \bfy^j \}_{j \in \mbN}$ is convergent and so is $ \{ \bfz^k \}_{k \in \mbN}$.
	The benefit of removing the repeated points and relabeling 
	allow us to restate the properties established in Prop.~\ref{Prop-SGSN-property} in a unified fashion:
	\begin{itemize}
		\item[$\bullet$] Denote $\eta_2:= \min\{ \eta_0, c_1 \}$, then we have
		\begin{align} \label{suf-des-y}
			F(\bfy^j) - F(\bfy^{j+1}) \geq \eta_2 \| \bfy^{j+1} - \bfy^j \|^2 .
		\end{align}
		\item[$\bullet$] There exists $\bfq^j \in \partial F(\bfy^j)$ such that 
		\begin{align} \label{subgradient-ub}
			\| \bfq^j \| \leq \eta_1 \| \bfy^{j+1} - \bfy^j \| .
		\end{align} 
		\item[$\bullet$] The Ostrowski's condition holds
		\begin{align} \label{os-con}
			\lim_{j \to \infty} \| \bfy^{j+1} - \bfy^j \| = 0 .
		\end{align}
	\end{itemize}
	
	These crucial properties enable us to utilize the P{\L}K property for convergence analysis of $\{ \bfy^j \}_{j \in \mbN}$.
	
	\begin{assumption} \label{semi-alge}
		$F(\cdot)$ is a P{\L}K function.
	\end{assumption}
	
	Particularly, the semi-algebraic function, which remains stable under basic operations, is an important class of the P{\L}K function (see e.g. \cite{bolte2014proximal}). We can give the following sufficient condition for Assumption \ref{semi-alge}. 
	
	\begin{lemma} \label{F-kl}
		If $f$ is semi-algebraic, then $F$ is semi-algebraic, and thus it is a P{\L}K function.
	\end{lemma}
	\bp
	Since $f$ is semi-algebraic, it follows from \cite[Example 2]{bolte2014proximal} that $\widetilde{f}(\bfv,\bfz):= \langle \bfv, \bfz \rangle - f(\bfz)$ is semi-algebraic, and thus $f^*(\bfv):= \sup_{\bfv \in \mathbb{R}^m} \widetilde{f}(\bfv,\bfz)$ is also semi-algebraic. As $h$ is the composition of $f^*$ and a linear mapping, it is semi-algebraic. Finally, taking semi-algebraicity of $\bfdt_+$ and $\| \cdot \|_0$ (see \cite[Example 4]{bolte2014proximal}) into account, $F(\cdot) = h(\cdot) + \bfdt_+(\cdot) + \mu \| \cdot \|_0 $ is semi-algebraic, which implies its P{\L}K property.
	\ep
	
	\begin{lemma} \label{accumulation set} Suppose Assumption \ref{seq-bound} holds and let $\Y$ be the set consisting of all the accumulation points of $\{ \bfy^j \}_{j \in \mbN}$. Then we have
		
		(i) $\lim_{j \to \infty} \hbox{\rm dist} (\bfy^j, \Y) = 0$.	
		
		(ii) $\Y$ is a nonempty, compact and connected set. 
		
		(iii) Each element $\bfy^* \in \Y$ is a P-stationary point of \eqref{01dco-uc} with $F(\bfy^*) = \lim_{j \to \infty} F(\bfy^j) = C$, where the constant $C$ has been defined in 
		Prop.~\ref{Prop-SGSN-property} (ii).
		
		(iv) If Assumption \ref{semi-alge} holds, then there exist $\epsilon > 0$, $s > 0$ and $\psi \in \Psi_s$ such that for any $\bfz \in \mathbb{B}^*:= \{ \bfz \in \mbR^m: {\rm dist}(\bfz, \Y) < \epsilon \} \cap [ C < F(\bfz) < C + s ]$, the following inequality holds
		\begin{align} \label{uniform-kl}
			\psi'( F(\bfz) - C )~{\rm dist}(0, \partial F(\bfz)) \geq 1.
		\end{align}  	
	\end{lemma}
	
	\bp
	Since Ostrowski's condition \eqref{os-con} holds, (i) and (ii) follow from \cite[Lem. 5]{bolte2014proximal}. (iii) follows from Prop.~\ref{Prop-SGSN-property}. Under Assumption \ref{semi-alge}, $F$ is a P{\L}K function with a constant value on $\Y$. Combining this with the compactness of $\Y$, we can derive (iv) by \cite[Lem. 6]{bolte2014proximal}.
	\ep
	
	Lem.~\ref{accumulation set}(iv) is also known as uniform P{\L}K property. This is a more favorable property than those introduced in Section \ref{pre}, because the data $\epsilon$, $s$ and $\psi$ is the same for all the points in $\mathbb{B}^*$, making the formula \eqref{uniform-kl} more applicable for convergence analysis.
	
	\begin{theorem} \label{Thm-Global}
		Suppose that Assumptions \ref{seq-bound} and \ref{semi-alge} hold. Then the sequence $\{ \bfz^k \}_{k \in \mbN}$ generated by SGSN converges to a P-stationary point of \eqref{01dco-uc}.
	\end{theorem}
	\bp
	Since $\{ \bfz^k \}_{k \in \mbN}$ is a subsequence of $\{ \bfy^j \}_{j \in \mbN}$, we can arrive at the desired conclusion by proving that $\{ \bfy^j \}_{j \in \mbN}$ converges to a P-stationary point of \eqref{01dco-uc}. First, let us consider a trivial case that there exists $\overline{j} \in \mbN$ satisfying $F(\bfy^\oj) = C$. Taking $\lim_{j \to \infty} F(\bfy^j) = C$ and the sufficient descent \eqref{suf-des-y} into account, we have $\bfy^j = \bfy^\oj$ for all $j \geq \oj$. As we have discussed in Remark~\ref{rem-sgsn}(iii), $\bfy^\oj$ actually is a P-stationary point of \eqref{01dco-uc}.
	
	Now let us consider the case when $F(\bfy^j) > C$ for any $j \in \mathbb{N}$. Since $\lim_{j \to \infty} \hbox{\rm dist} (\bfy^j, \Y) = 0$ and $\lim_{j \to \infty} F(\bfy^j) = C$, it follows from Lem.~\ref{accumulation set}(iv) that when $j$ is sufficiently large, we have
	\begin{align*}
		\psi'( F(\bfy^j) - C )~{\rm dist}(0, \partial F(\bfy^j)) \geq 1.
	\end{align*}  
	This together with \eqref{subgradient-ub} leads to 
	\begin{align} \label{psi'-lb}
		\psi'( F(\bfy^j) - C ) \geq 1/( \eta_1 \| \bfy^j - \bfy^{j-1} \| ) .
	\end{align}
	Given $l_1,l_2 \in \mbN$, denote $\Delta_{l_1,l_2}:=\psi( F(\bfy^{l_1}) - C ) - \psi( F(\bfy^{l_2}) - C )$, then by the concavity of $\psi$, \eqref{psi'-lb} and \eqref{suf-des-y}, we can derive
	\begin{align*}
		\Delta_{j,j+1} \geq \psi'( F(\bfy^j) - C ) ( F(\bfy^j) - F(\bfy^{j+1}) ) \geq \frac{\eta_2}{\eta_1} \frac{\| \bfy^{j+1} - \bfy^j \|^2}{\| \bfy^j - \bfy^{j-1} \|},
	\end{align*}
	which leads to $\| \bfy^{j+1} - \bfy^j \| \leq \sqrt{(\eta_1/\eta_2) \Delta_{j,j+1} \| \bfy^j - \bfy^{j-1} \|}$. Using the inequality of arithmetic and geometric means, we have
	\begin{align*}
		2\| \bfy^{j+1} - \bfy^j \| \leq \| \bfy^j - \bfy^{j-1} \| 
		+ (\eta_1/\eta_2) \Delta_{j,j+1} .
	\end{align*}
	Now let us sum up the above inequality from sufficiently large $\underline{l}$ to $\overline{l}$ to get
	\begin{align*}
		2 \sum_{j =\ul}^{\ol} \| \bfy^{j+1} - \bfy^j \| \leq & \sum_{j =\ul}^{\ol} \| \bfy^j - \bfy^{j-1} \| + \sum_{j =\ul}^{\ol} \frac{\eta_1}{\eta_2} \Delta_{j,j+1} \\
		= & \sum_{j =\ul}^{\ol} \| \bfy^{j+1} - \bfy^j \| + \| \bfy^\ul - \bfy^{\ul-1} \| + \frac{\eta_1}{\eta_2} \Delta_{\ul,\ol+1}.
	\end{align*}
	This further leads to
	\begin{align*}
		\sum_{j =\ul}^{\ol} \| \bfy^{j+1} - \bfy^j \| \leq \| \bfy^\ul - \bfy^{\ul-1} \| + (\eta_1/\eta_2) (\psi( F(\bfy^{\ul}) - C ) - \psi( F(\bfy^{\ol+1}) - C )).
	\end{align*}
	Taking $\ol \to \infty$, we can derive
	\begin{align}
		\sum_{j=1}^{\infty} \| \bfy^{j+1} - \bfy^j \| < \infty .
	\end{align} 
	This indicates that $\{ \bfy^j \}_{j \in \mbN}$ is a Cauchy sequence. It must converge to a P-stationary point of \eqref{01dco-uc}, and so does $\{ \bfz^k \}_{k \in \mbN}$. We have finished the proof. 
	\ep
	
	\subsection{Local Superlinear Convergence Rate} \label{locl-superlinear}
	
	In this Subsection, we will derive the local superlinear convergence rate of SGSN under the following assumption.
	
	\begin{assumption} \label{semismooth+sosc}
		$\nabla f^* (\cdot)$ is a semismooth function, and there exists an accumulation point $\bfz^*$ of the sequence $\{ \bfz^k \}_{k \in \mbN}$ satisfying  
		\begin{align} \label{sosc1}
			\bfd^\top H^*_{\S^*} \bfd > 0,~~ \forall H^* \in \whp^2 h (\bfz^*) \hbox{ and } \bfd \in \mbR^{| \S^* |} \backslash \{ 0 \}.
		\end{align}
	\end{assumption}
	By the relationship \eqref{chain-rule}, we can derive that \eqref{sosc1} is stronger than \eqref{sosc}.
	Assumption \ref{semismooth+sosc} is mainly used to ensure the nonsingularity for each element of the generalized Jacobian in the neighborhood of the reference point $\bfz^*$. 
	We have the following result.
	
	\begin{lemma} \label{bounded}
		If Assumption \ref{semismooth+sosc} holds, there exist radius $\epsilon^* > 0$ and $\eta_3,\eta_4 > 0$ such that 
		\begin{align} \label{bound}
			\eta_3 \leq \inf_{ H \in \whp^2 h(\bfz), \atop \bfz \in \N(\bfz^*,\epsilon^*) }  \bfla_{\min} (H_{\S^*}) \leq	\sup_{ H \in \whp^2 h(\bfz), \atop \bfz \in \N(\bfz^*,\epsilon^*) }  \bfla_{\max} (H_{\S^*}) \leq \eta_4 .
		\end{align}
	\end{lemma}
	\bp
	Under Assumption \ref{semismooth+sosc}, for any $H^* \in \whp^2 h (\bfz^*)$, $H^*_{\S^*}$ is positive definite. According to \cite[Prop.~7.1.4]{facchinei2003finite}, $\partial^2 f^*(-A^\top \bfz^*)$ is a nonempty and compact set. Then the same property is true for $\whp h (\bfz^*)$ by the definition \eqref{hat-partial}. Therefore, there exist $\eta_3, \eta_4 > 0$ such that 
	\begin{align} \label{bound-H*}
		\inf_{ H^* \in \whp^2 h(\bfz^*)}  \bfla_{\min} (H^*_{\S^*}) \geq 2\eta_3,~~~ \sup_{ H^* \in \whp^2 h(\bfz^*)}  \bfla_{\max} (H^*_{\S^*}) \leq \frac{\eta_4}{2} .
	\end{align}
	By the upper semi-continuity of $\partial^2 f^*$ (see \cite[Prop.~7.1.4]{facchinei2003finite}), there exists $\epsilon^* > 0$ such that 
	\begin{align*}
		\partial^2 f^* (-A^\top \bfz)  \subseteq \partial^2 f^*(-A^\top \bfz^*) + \N ( 0, \delta_1^*), ~\forall \bfz \in \N (\bfz^*, \epsilon^*),
	\end{align*}
	where $\delta_1^*:= \min\{ \eta_3, \eta_4/2 \}/ \| A_{\S^*:} \|^2. $ Then for any $\bfz \in \N (\bfz^*, \epsilon^*)$ and $H = A Q A^\top \in \whp h( \bfz )$, we can find $H^* = A Q^* A^\top \in \whp h( \bfz^* )$ with $\| Q - Q^* \| \leq \delta^*_1$. Using the Weyl's perturbation theorem (see \cite[Corollary III.2.6]{bhatia2013matrix}), we have 
	\begin{align*}
		& \max\{ | \bfla_{\min} (H_{\S^*}) -  \bfla_{\min} (H^*_{\S^*}) |, 
		| \bfla_{\max} (H_{\S^*}) -  \bfla_{\max} (H^*_{\S^*}) | \} \leq \| H_{\S^*} - H^*_{\S^*} \| \\
		& \leq \| A_{\S^* :} (Q - Q^*) (A_{\S^* :})^\top \| \leq \min\{ \eta_3, \eta_4/2 \} .
	\end{align*}
	Combining this with \eqref{bound-H*}, we can derive \eqref{bound}.
	\ep

	\begin{remark}
		(i) \eqref{chain-rule} implies that for any $\overline{\bfd} \in \mbR^m$, $ \overline{\bfd}^\top(\partial^2 h(\bfz)) \overline{\bfd} \subseteq \overline{\bfd}^\top (\whp^2 h(\bfz)) \overline{\bfd}$. Then by the Rayleigh-Ritz theorem (see \cite[Chapter 5]{lutkepohl1997handbook}), we can also obtain
		\begin{align} \label{bound1}
			0 < \eta_3 \leq \inf_{ H \in \partial^2 h(\bfz), \atop \bfz \in \N(\bfz^*,\epsilon^*) }  \bfla_{\min} (H_{\S^*}) \leq	\sup_{ H \in \partial^2 h(\bfz), \atop \bfz \in \N(\bfz^*,\epsilon^*) }  \bfla_{\max} (H_{\S^*}) \leq \eta_4 .
		\end{align}	
		(ii) Semismoothness of $\nabla f^*$ implies semismoothness of $\nabla h$. This can be verified by \cite[Thm.~7.4.3]{facchinei2003finite}. Therefore, Assumption \ref{semismooth+sosc} ensures the SOSC in Prop.~\ref{Prop-QuadraticGrowth}.
		
		(iii) As we have shown in Prop.~\ref{Prop-SGSN-property},
		if Assumption \ref{seq-bound} holds, $\bfz^*$ will be a P-stationary point and $F$ is constant on all the accumulation points of $\{ \bfz^k \}_{k \in \mbN}$. If we further take Assumption \ref{semismooth+sosc} into consideration, SOSC holds and 
		Prop.~\ref{Prop-QuadraticGrowth} ensure that each accumulation point of $\{ \bfz^k \}_{k \in \mbN}$ is isolated. Since $\lim_{k \to \infty } \| \bfz^{k+1} - \bfz^k \| = 0$ (derived by \eqref{successive-change}), it follows from \cite{kanzow1999qp} that $\lim_{k \to \infty} \bfz^k = \bfz^*$. Then from \eqref{successive-change}, we also have $\lim_{k \to \infty} \bfv^k = \bfz^*$.
	\end{remark}
	
	To derive the local superlinear convergence of SGSN, we need to overcome several obstacles. First, SGSN performs Newton iterate on changing subspace $\mbR^{|T_k|}$. The technique for the convergence rate of classic Newton method is unavailable unless $\S^*$ is exactly identified by $T_k$. Moreover, the full-Newton step is only implemented when \eqref{c1}, \eqref{c2} and $\alpha_k = 1$ are met. 
	In the subsequent lemma, we show that these conditions hold after finite iterations.
	
	\begin{lemma} \label{pre-superlinear}
		Suppose that Assumptions \ref{seq-bound} and \ref{semismooth+sosc} hold. Then there exists a sufficiently large $k^* \in \mbN$ such that for $k \geq k^*$
		
		(i) The support set of $\bfz^*$ is identified by $T_k$, i.e.
		\begin{align}
			T_k = \S(\bfv_k) = \S^* .
		\end{align} 
		
		(ii) The step length $\alpha_k$ always equals to 1.
		
		(iii) If $c_1 \leq \eta_3/4$ and $c_2 \geq \ell_h + \eta_4 + 1$, then the Newton iterate $\wz^{k+1}$ is always accepted.
	\end{lemma}
	
	\bp
	(i) As we just discussed, $\lim_{k \to \infty} \bfv^k = \bfz^*$ holds. Then when $k$ is sufficiently large, we have $\S(\bfv^k) \supseteq \S^*$. Moreover, according to 
	Prop.~\ref{Prop-SGSN-property}, $\lim_{k\to \infty} F(\bfv^k) = F(\bfz^*)$ holds, which further leads to $\lim_{k \to \infty} \| \bfv^k \|_0 = \| \bfz^* \|_0$. 
	Noticing that the values of $\| \cdot \|_0$ are from the set $\{0, 1, \ldots, m\}$, 
	therefore $\| \bfv^k \|_0 = \| \bfz^* \|_0$ must be true when $k$ is large enough. Combining this with $\S(\bfv^k) \supseteq \S^*$, we can derive assertion (i).
	
	(ii) Now let us consider $k$ is large enough such that $T_k = \S^*$ holds. From Lem.~\ref{bounded}, $H^k_{T_k}$ is nonsingular and we can estimate
	\begin{align*}
		\| \bfd^k_{T_k} \| \leq \| (H^k_{T_k} + \gamma_k I)^{-1} \nabla_{T_k} h(\bfv^k) \| \leq \| \nabla_{T_k} h(\bfv^k) \|/\eta_3.
	\end{align*}
	Since $\nabla_{\S^*} h(\bfz^*) = 0$ holds from \eqref{01dco-kkt}, we have $\lim_{k \to \infty} \| \bfd^k_{T_k} \| = 0$. Meanwhile, since for any $i \in T_k = \S^*$, we have $\lim_{k \to \infty} v^k_i = z^*_i > 0$. Then according to \eqref{step-beta}, when $k$ is large enough $\beta_k \geq 1$, and thus $\alpha_k = 1$.
	
	(iii) We need to show that \eqref{c1} and \eqref{c2} hold when $k$ is large enough. Let us begin with some preliminary estimations. Using $\nabla_{T_k} h(\bfz^*) = 0$ and Lipschitz continuity of $h$, we can estimate
	\begin{align} \label{ga-bound}
		\gamma_k = \| \nabla_{T_k} h(\bfv^k) \| = \| \nabla_{T_k} h(\bfv^k) - \nabla_{T_k} h(\bfz^*) \| \leq  \ell_h  \| \bfv^k - \bfz^* \| .
	\end{align}
	Let us denote $\Phi := \nabla f^*$ and $H^k = AQ^kA^\top$ with $Q^k \in \partial^2 f^* (-A^\top \bfv^k) = \partial_C \Phi(-A^\top \bfv^k)$. Since $\alpha_k = 1$ when $k$ is sufficiently large, from \eqref{Nt-eq} and \eqref{newton-step}, we have
	\begin{align}
		&\| \wz^{k+1} - \bfz^* \| \notag \\
		= & \| (\wz^{k+1} - \bfz^*)_{T_k} \| = \|  (\bfv^k - \bfz^*)_{T_k} - ( H^k_{T_k} + \gamma_k I )^{-1} \nabla_{T_k} h(\bfv^k) \| \notag \\
		\leq & \frac{1}{\eta_3} \| ( H^k_{T_k} + \gamma_k I ) (\bfv^k - \bfz^*)_{T_k} - (\nabla h (\bfv^k) - \nabla h(\bfz^*))_{T_k} \| \notag \\
		=& \frac{1}{\eta_3}\| H^k_{T_k} (\bfv^k - \bfz^*)_{T_k} - (\nabla h (\bfv^k) - \nabla h(\bfz^*))_{T_k} \| +  o(\| \bfv^k - \bfz^* \|) \notag \\
		=& \frac{1}{\eta_3} \| A_{T_k:} Q^k A^\top_{T_k:} (\bfv^k - \bfz^*)_{T_k} + A_{T_k:} ( \Phi (- A^\top \bfv^k ) - \Phi (-A^\top \bfz^*) ) \| +  o(\| \bfv^k - \bfz^* \|) \notag \\
		\leq &( \|A\|/\eta_3) \| \Phi (-A^\top \bfv^k) - \Phi (-A^\top \bfz^*) - Q^k(-A^\top(\bfv^k - \bfz^*)) \|  + o(\| \bfv^k - \bfz^* \| ) \notag \\
		= &o(\| \bfv^k - \bfz^* \|) , \label{superlinear}
	\end{align}
	where the second and last equations follow from \eqref{ga-bound} and the semismoothness property \eqref{def-semismoothness} respectively. The above relationship also implies $\lim_{k \to \infty} \| \wz^{k+1} - \bfv^k \| = 0$. From \eqref{Nt-eq} and \eqref{ga-bound}, we can also estimate
	\begin{align*}
		\| \nabla_{T_k} h ( \bfv^k ) \| = & \| (H^{k+1/2}_{T_k} + \gamma_k I) (\wz^{k+1} - \bfv^k)_{T_k} \| \\
		\leq & \eta_4 \| \wz^{k+1} - \bfv^k \| + o( \| \wz^{k+1} - \bfv^k \| ) .
	\end{align*}
	Now we are ready to prove \eqref{c2}. When $k$ is sufficiently large, it follows from $o( \| \wz^{k+1} - \bfv^k \| ) \leq \| \wz^{k+1} - \bfv^k \|$ and Lipschitz continuity of $h$ that 
	\begin{align*}
		\| \nabla_{T_k} h(\wz^{k+1}) \| \leq &	\| \nabla_{T_k} h ( \bfv^k ) \| + \| \nabla_{T_k} h ( \bfv^k ) - \nabla_{T_k} h(\wz^{k+1}) \| \\
		\leq & (\ell_h + \eta_4 + 1) \| \wz^{k+1} - \bfv^k \| .
	\end{align*}
	Finally, let us prove \eqref{c1}. According to \cite[Prop.~7.4.10]{facchinei2003finite} and \eqref{superlinear}, we have
	\begin{align}
		h(\wz^{k+1}) =& h(\bfz^*) + \langle \nabla h(\bfz^*), \wz^{k+1} - \bfz^* \rangle + \frac{1}{2}(\wz^{k+1} - \bfz^*)^\top \widetilde{H}^{k+1} (\wz^{k+1} - \bfz^*) \notag \\
		&+ o(\| \wz^{k+1} - \bfz^* \|^2) \notag \\
		= &   h(\bfz^*) + \langle \nabla h(\bfz^*), \wz^{k+1} - \bfz^* \rangle + o(\| \bfv^k - \bfz^* \|^2) \label{taylor1} \\
		h(\bfv^k) =& h(\bfz^*) + \langle \nabla h(\bfz^*), \bfv^{k} - \bfz^* \rangle + \frac{1}{2}(\bfv^{k} - \bfz^*)^\top P^k (\bfv^{k} - \bfz^*) \notag \\ &+ o(\| \bfv^{k} - \bfz^* \|^2), \label{taylor2}
	\end{align}
	for any $\widetilde{H}^{k+1} \in \partial^2 h(\wz^{k+1})$ and $P^k \in \partial^2 h (\bfv^k)$.
	Moreover, using $T^k = \S^*$ and $\nabla_{S^*} h(\bfz^*) = 0$, we have $\langle \nabla h(\bfz^*), \wz^{k+1} - \bfz^* \rangle = \langle \nabla_{T_k} h(\bfz^*), (\wz^{k+1} - \bfz^*)_{T_k} \rangle = 0$ and $\langle \nabla h(\bfz^*), \bfv^{k} - \bfz^* \rangle = \langle \nabla_{T_k} h(\bfz^*), (\bfv^{k} - \bfz^*)_{T_k} \rangle = 0$. Then subtracting \eqref{taylor1} and \eqref{taylor2} leads to
	\begin{align*}
		h(\bfv^k) - h (\wz^{k+1}) \geq \frac{\eta_3}{2}  \| \bfv^{k} - \bfz^* \|^2 + o( \| \bfv^{k} - \bfz^* \|^2 ) .
	\end{align*}
	From $\| \wz^{k+1} - \bfv^k \| \leq \| \wz^{k+1} - \bfz^* \| +\| \bfv^k - \bfz^* \|$ and \eqref{superlinear}, we have $\| \wz^{k+1} - \bfv^k \|^2 \leq \| \bfv^{k} - \bfz^* \|^2 + o(\| \bfv^k - \bfz^* \|^2)$, and thus we can verify
	\begin{align*}
		h(\bfv^k) - h (\wz^{k+1}) \geq & \frac{\eta_3}{4}  \| \bfv^{k} - \bfz^* \|^2 + \frac{\eta_3}{4}  \| \bfv^{k} - \wz^{k+1} \|^2 + o( \| \bfv^{k} - \bfz^* \|^2 ) \\
		\geq & \frac{\eta_3}{4}  \| \bfv^{k} - \wz^{k+1} \|^2 .
	\end{align*}
	Since $\wz^{k+1}_{\oT_k} = 0$ and $T_k = \S(\bfv^k)$, we have $g(\wz^{k+1}) \leq g(\bfv^k)$. Finally, we can prove \eqref{c1} by 
	\begin{align*}
		F(\bfv^k) - F(\wz^{k+1}) \geq  h(\bfv^k) - h(\wz^{k+1}) \geq \frac{\eta_3}{4}  \| \bfv^{k} - \wz^{k+1} \|^2 .
	\end{align*}
	Overall, we have completed the proof.
	\ep
	
	Lem.~\ref{pre-superlinear}(iii) indicates that $c_1$ and $1/c_2$ should be small enough to ensure the acceptance of Newton iterate $\wz^{k+1}$. 
	In fact, we can give a more specific lower bound for $c_2$. 
	By the $(1/\sigma_f)$-smoothness of $f^*$, the upper bound in \eqref{bound-H*} will be $\sup_{ H^* \in \whp^2 h(\bfz^*)}  \bfla_{\max} (H^*_{\S^*}) \leq \|A\|^2/ \sigma_f = \ell_h := \eta_4/2$.
	Then we can set $c_2 \geq 3 \ell_h + 1$. 
	Moreover, an upper bound for $c_1$ can be computed in a special case that $f$ is $\ell_f$-smooth and $A_{\S^*:}$ has full row rank. In this case, $f^*$ is $(1/\ell_f)$-strongly convex, and thus the lower bound in \eqref{bound-H*} will be $	\inf_{ H^* \in \whp^2 h(\bfz^*)}  \bfla_{\min} (H^*_{\S^*}) \geq \sqrt{\bfla_{\min}(A_{\S^*:}A^\top_{\S^*:})}/\ell_f:= 2\eta_3$. Then we can set $c_1 \leq \sqrt{\bfla_{\min}(A_{\S^*:}A^\top_{\S^*:})}/(8\ell_f)$.
	
	
	
	We are ready to state the local superlinear convergence of SGSN.
	
	\begin{theorem} \label{Thm-Local}
		Suppose that Assumptions \ref{seq-bound} and \ref{semismooth+sosc} hold. If we set $c_1 \leq \eta_3/4$ and $c_2 \geq \ell_h + \eta_4 + 1$, then the sequence $\{\bfz^k\}_{k \in \mbN}$ converges to $\bfz^*$ with local superliner rate, i.e. when $k$ is sufficiently large, we have
		\begin{align*}
			\| \bfz^{k+1} - \bfz^* \| \leq o( \| \bfz^k - \bfz^* \| ) .
		\end{align*}
	\end{theorem}
	\bp
	Since we have proved \eqref{superlinear} and $
	\bfz^{k+1} = \wz^{k+1}$ when $k$ is sufficiently large, we have
	\begin{align} \label{superlinear1}
		\| \bfz^{k+1} - \bfz^* \| \leq o( \| \bfv^{k} - \bfz^* \| ).
	\end{align}
	It suffices to find the relationship between $\| \bfz^k - \bfz^* \|$ and $\| \bfv^k - \bfz^* \|$. Using \eqref{pg-step}, $\nabla_{T_k} h(\bfz^*) = 0$, and Lipschitz continuity of $h$, we have
	\begin{align*}
		\| \bfv^k - \bfz^* \| =& \| (\bfv^k - \bfz^*)_{T_k} \| = \|  (\bfz^k - \bfz^*)_{T_k} -\tau( \nabla h( \bfz^k ) - \nabla h(\bfz^*) )_{T_k} \| \\
		\leq& (\ell_h \tau + 1) \| \bfz^k - \bfz^* \|.
	\end{align*}
	Combining this with \eqref{superlinear1}, we can derive our desired conclusion.
	\ep 
	
	\begin{remark} \label{rem-convergence-rate}
		If the semismoothness of $\nabla f^*$ in Assumption \ref{semismooth+sosc} is strengthened to strong semismoothness (see e.g. \cite{facchinei2003finite}), then similar to the classical result in \cite{qi1993nonsmooth}, the convergence rate in Thm.~\ref{Thm-Local} can be improved to be quadratic. 
	\end{remark}
	
	As emphasized in Introduction, various types of semismooth Newton methods
		have been recently developed for nonsmooth equations arising from optimization, see e.g.,
		\cite{zhou2021newton,kanzow2021globalized,khanh2024coderivative,khanh2024globally,li2018highly,zhang2024zero}.
		We make detailed comments below on three such methods: Classical semismooth Newton method \cite{qi1993nonsmooth}, Coderivative-based Newton method \cite{khanh2024coderivative}, and Subspace Newton method \cite{zhou2021newton}.
		The comparisons will highlight the unique features of our method and will be also
		useful when conducting numerical comparison in the next section. 
	
		\begin{remark}[\textbf{Motivation from classical semismooth Newton method}] 
			In general, classical semismooth Newton method serves a guidance for constructing
			new methods for nonsmooth equations. Let us consider the proximal gradient equation:
			\be \label{Eq-Semismooth}
			0 = G(\bfz) := \bfz - \Prox_{\tau g(\cdot)}( \bfz - \tau \nabla h(\bfz)).
			\ee
			This equation is well-defined at any P-stationary point when $\tau$ is sufficiently small, see Lem.~\ref{property-g}(vi).
			We assume that $G(\cdot)$ is locally Lipschitzian
			and semismooth. A semismooth/generalized Newton method
			gives the following update at $\bfz^k$:
			\be \label{SemismoothNewton}
			\left\{
			\begin{array}{rl}
				(A^k + \gamma_k I) \bfd^k &= - G(\bfz^k), \quad A^k \in \partial_C G(\bfz^k)\\
				\bfz^{k+1} &= \bfz^k + \bfd^k,
			\end{array} 
			\right .
			\ee
			where $\partial_C G(\cdot)$ is the generalized Jacobian in the sense of Clarke \cite{clarke1990optimization} and $\gamma_k \ge 0$ is a regularization parameter.
			Various generalized Newton methods mainly differ in three aspects: 
			(i) using different generalized Jacobian instead of $\partial_C G(\cdot)$, 
			(ii) proposing approximation property of $G(\cdot)$ in terms of the proposed Jacobian, and 
			(iii) developing globalization strategy. 
			For (i), we used the generalized Jacobian $\whp^2 h (\bfz)$ defined in \eqref{hat-partial}.
			For (ii), we assumed the semismoothness of $\nabla f^*(\cdot)$. 
			We make more comments on Point (iii).
			
			It is often a nontrivial task to globalize the Newton method in \eqref{SemismoothNewton}.
			A common strategy is to use Newton direction $\bfd^k$ if certain conditions hold and to use a gradient descent direction otherwise.
			Our method roughly follows this strategy.
			Firstly, we work with the set-valued mapping $\Prox_{\tau g}(\cdot)$ 
			and select an element $\bfv^k$ from the set of proximal gradients
			defined in \eqref{pg-step1}.
			Secondly, according to \cite[Prop.~2.3]{khanh2024coderivative}, 
			$\Prox_{\tau g} (\cdot)$ is Lipschitz continuous and single-valued around $\bfz^* - \tau\nabla h(\bfz^*)$ when $\tau$ is sufficiently small and $\bfz^*$ is a P-stationary point.
			This result is useful because 
			locally the Newton equation \eqref{Nt-eq} in SGSN can be derived from \eqref{SemismoothNewton}
			using the index set $T_k$ and the associated properties established in Lem.~\ref{property-g}.
			In our SGSN, $\tau$ can take any value in the range $(0, 1/\ell_h)$.
			Thirdly, this is probably the most difficult part and it is
			to  decide when to use the gradient
			direction $\bfv^k$ or Newton direction $\bfd^k$. 
			Our selection is decided by the two conditions \eqref{c1} and \eqref{c2}, 
			which are carefully designed for the Newton step to be eventually accepted
			and for the convergence to hold for a large class 
			of P{\L}K functions $F(\bfz) = h(\bfz) + g(\bfz)$.
			
			In summary, our SGSN algorithm is motivated by the classical Semismooth Newton method.
			But its convergence does not follow from existing convergence theory due to
			some innovative strategies proposed, especially on \eqref{c1} and \eqref{c2}.
	\end{remark}

	\begin{remark}[\textbf{Comparison with a coderivative-based Newton method}]
			When $h \in {\rm C}^2$ (twice continuously differentiable) and $g$ is prox-regular, 
			the coderivative-based Newton (CBN) method \cite[Alg.~4]{khanh2024coderivative} can be applied for Problem \eqref{01dco-uc}.
			It first implements a proximal gradient step \eqref{pg-step} to compute $\bfv^k$, and  followed by a Newton step to compute direction $\bfd^k$:
			\begin{align} \label{coderivative-newton}
				- \bfs^k - H^k \bfd^k \in \partial^2 g( \bfv^k, \bfr^k ) (\bfd^k),
			\end{align}
			where $H^k = \nabla^2 h(\bfv^k)$, $\bfs^k:= \nabla h(\bfv^k) - \nabla h(\bfz^k) + (\bfz^k - \bfv^k)/\tau$, $\bfr^k:= \bfs^k - \nabla h(\bfv^k)$ and the second-order subdifferential (see e.g. \cite[Def.~3.17]{mordukhovich2018variational}) of $g$ at $\bfv$ relative to $\bfr \in \partial g(\bfv)$ is defined by
			\begin{align} \label{second-order-subdiff}
				\partial^2 g(\bfv, \bfr)(\bfd) = \left\{ \bfq \in \mbR^m \left|  q_i \in \left\{  \begin{aligned}
					& \{ t~|~td_i = 0 \}, && \mbox{if} \ v_i = r_i = 0, \\
					& \mbR, && \mbox{if} \ v_i = 0, r_i \neq 0, d_i = 0, \\
					& \{0\}, && \mbox{if} \ v_i > 0, \\
					& \emptyset, && \mbox{otherwise}.
				\end{aligned} \right.  \right. \right\}
			\end{align}
			using the notation $T_k = \S(\bfv^k)$ and $\bfv^k_{\oT_k} = 0$ from \eqref{Tk-support}, we can compute a solution of \eqref{coderivative-newton} by solving the following linear equation:
			\begin{align*}
				H^k_{T_k} \bfd^k_{T_k} = -\bfs^k_{T_k} \ \mbox{and} \ \bfd^k_{\oT_k} = 0.
			\end{align*}
			Ignoring the regularization term $\gamma_k I$ in \eqref{Nt-eq},
			CBN and SGSN appear similar. But they are very different in a few aspects.
			
			{\bf (A) Globalization strategies are different.}
			Although both methods share the same coefficient matrix $H^k_{T_k}$, the right-hand side vectors are different, $-\bfs^k_{T_k}$ vs $-\nabla_{T_k} h(\bfv^k)$.
			Moreover, CBN adopts a line search technique based on a sufficient descent on a forward-back envelope \cite[Def.~5.1]{khanh2024coderivative} for globalization.
			In contrast, SGSN uses \eqref{step-beta} and Conditions \eqref{c1} or \eqref{c2}.
			
			
			{\bf (B) Assumptions on global convergence are different.}
			CBN requires $h \in C^2$ and $g$ being prox-regular for its global convergence
			under the additional assumption that an accumulation point is isolated.
			In contrast, SGSN assume P{\L}K property of the objective function $F$, allowing the accumulation points to be non-isolated.

			{\bf (C) Assumptions on superlinear convergence are different.}
			Under the premise of global convergence to a {\rm P}-stationary point $\bfz^*$, the local superlinear convergence of CBN further requires 
			\begin{itemize}
				\item[(i)] $g$ is continuously prox-regular (see e.g. \cite[Def.~3.27]{mordukhovich2018variational}) at $\bfz^*$ for $-\nabla h(\bfz^*)$.
				
				\item[(ii)] $\nabla^2 h$ is strictly differentiable (see e.g. \cite[Def.~9.17]{RockWets98}) at $\bfz^*$.
				
				\item[(iii)] $\bfz^*$ is a tilt-stable local minimizer (see e.g. \cite[Def.~1.1]{poliquin1998tilt}) of \eqref{01dco-uc}.   
			\end{itemize}
			Similar to the proof of \cite[Thm.~7.3]{khanh2024coderivative}, (i) holds if and only if $[\bfz^* + \nabla h(\bfz^*)]_i \neq 0$ for all $i \in [m]$, and then through the use of representation \eqref{second-order-subdiff}, \cite[Prop.~1.121]{mordukhovich2006variational} and \cite[Thm.~1.3(b)]{poliquin1998tilt}, one can derive that (iii) is equivalent to positive definiteness of $H^*_{\S^*}$, where $H^*:= \nabla^2 h(\bfz^*)$. Meanwhile, the local superlinear convergence of SGSN is guaranteed by Assumption \ref{semismooth+sosc}. Particularly, in the case that $h$ is ${\rm C}^2$-smooth, this assumption becomes positive definiteness of $H^*_{\S^*}$. Overall, convergence of SGSN is established under weaker assumptions. 
	\end{remark}

	\begin{remark}[\textbf{Comparison with a subspace Newton method}] \label{Remark-CBN}
			\cite{zhou2021newton} proposed a Newton method for $\ell_0$-regularized optimization (NL0R). It also implements Newton step on a subspace determined by certain index sets. Moreover, it enjoys global convergence with local quadratic rate under suitable conditions. 
			However, both methods are significantly different.
			
			\begin{itemize}
				\item[(i)] Problem types: NL0R solves the minimization of a ${\rm C}^2$-smooth function with $\ell_0$ regularization, while SGSN is designed for a more challenging problem \eqref{01dco-uc}, where $h$ is ${\rm SC}^1$ and $g$ comprises a $\ell_0$ regularizer and a nonnegative constraint. 
				NL0R cannot be directly applied to \eqref{01dco-uc}.
				The nonsmoothness of $\nabla h$ and the constraint in $g$ demand more sophisticated algorithmic design to ensure both feasibility and convergence of the iterates. 
				
				\item[(ii)] Globalization strategy: NL0R adopts Armijo line search on Newton direction. In comparison, SGSN combines proximal gradient and semismooth Newton steps with the aid of conditions \eqref{c1} and \eqref{c2}. This helps us utilize P{\L}K property to globalize SGSN. 
				
				\item[(iii)] Convergence: The whole sequence convergent of NL0R requires existence of an isolated accumulation point, whereas SGSN adopts the P{\L}K property. For local quadratic convergence to a solution $\bfz^*$, NL0R requires that the smooth part in the objective function has positive definite and Lipschitz continuous Hessian around $\bfz^*$. In comparison, SGSN requires weaker assumptions that $h$ is strongly semismooth and its generalized Hessian at $\bfz^*$ is positive definite on a subspace indexed by the support set $\S^*$.
			\end{itemize}
	\end{remark}

	\section{Numerical Experiments} \label{numerical-experiments}
	
	In this section, we demonstrate the performance of SGSN on AUC maximization and sparse multi-label classification. All the experiments are implemented on Matlab 2022a by a laptop with 32GB memory and Intel CORE i7 2.6 GHz CPU.
	
	\subsection{AUC Maximization} \label{auc}
	The area under the receiver operating characteristic curve (AUC) \cite{hanley1982meaning} is a widely used evaluation metric in imbalanced classification and anomaly detection.  Let $X^+ = [\bfx^+_1, \cdots, \bfx^+_{q_+}]^\top \in \mbR^{q_+ \times n}$ and $X^- = [\bfx^-_1, \cdots, \bfx^-_{q_-}]^\top \in \mbR^{q_- \times n}$ be the positive and negative sample matrices respectively, then the AUC associated with parameter vector $\bfx \in \mbR^n$ can be represented as 
	\begin{align*}
		\texttt{AUC}( \bfx ) := \frac{1}{q_+q_-}  \sum_{i = 1}^{q_+} \sum_{j = 1}^{q_-} \bfone_{(0,\infty)} \left(  \langle \bfx^+_i, \bfx \rangle - \langle \bfx^-_j, \bfx \rangle \right).
	\end{align*}
	After an opportune scaling of $\bfx$ (a trick commonly used in support vector machines),
	each term $(\langle \bfx^+_i, \bfx \rangle - \langle \bfx^-_j, \bfx \rangle)$ is replaced
	by $(\langle \bfx^+_i, \bfx \rangle - \langle \bfx^-_j, \bfx \rangle -1)$, see e.g. \cite{gao2013one}.
	This scaling operation also removes the trivial solution $\bfx^*=0$.

	Therefore, AUC can be maximized by solving composite optimization \eqref{01cop} with the following setting:
	\begin{align} \label{auc-primal}
		f(\bfx) = \frac{1}{2} \| \bfx \|^2,~ A = (\bfe_{q_+} \otimes I_{q_-}) X^- -  ( I_{q_+} \otimes \bfe_{q_-} ) X^+ \in \mbR^{m \times n},~\bfb = \bfe_{m},
	\end{align}
	where $\otimes$ is the Kroneker product and $m := q_+q_-$. Then the conjugate function $f^*(\bfv) = \| \bfv \|^2/2$, and thereby the associated dual problem \eqref{01dco-uc} can be written as 
	\begin{align} \label{auc-dual}
		\min_{\bfz \in \mbR^m}\; F(\bfz) = \frac{1}{2} \| A^\top \bfz \|^2 - \langle \bfb, \bfz \rangle + g(\bfz).
	\end{align} 
	In this case, $f$ is semi-algebraic, $f^*$ is strongly convex and $\nabla^2 f^*$ is Lipschitz continuous. When SGSN is applied to solving \eqref{auc-dual}, provided that the generated sequence is bounded, then it must converge to a P-stationary point $\bfz^*$ of \eqref{auc-dual} according to Thm.~\ref{Thm-Global}. Additionally, if $A_{\S^*:}$ is a full row-rank matrix, then the sequence converges with local quadratic rate by Remark \ref{rem-convergence-rate}.
	
	In the experiment of AUC maximization, we select four other algorithms for comparison. 
	NM01 \cite{zhou2021quadratic} is a subspace Newton-type method designed for solving the primal Problem \eqref{01cop} with a twice continuously differentiable $f$. 
	CBN \cite[Alg.~4]{khanh2024coderivative} 
		was already commented in Remark~\ref{Remark-CBN}.
	PRSVM \cite{chapelle2010efficient} is a semismooth Newton method which optimizes AUC with a squared hinge surrogate loss function. 
	OPAUC \cite{gao2013one} is a one-pass gradient descent method for maximizing AUC with a least square surrogate loss function.  
	
	In the setting of \eqref{auc-primal} and \eqref{auc-dual}, we compute the Lipschitzian constant of $\nabla h$ by $\ell_h = \| A \|^2 = q_+\| X^- \|^2 + q_-\| X^+ \|^2 - 2 \langle \bfe_{q_-}^\top X^-, \bfe_{q_+}^\top X^+ \rangle$. Then for SGSN, we set $\mu = \tau = 1/(2 \ell_h )$, $\gamma = 10^{-1}$, $c_1 = 1/(3\ell_h)$ and $c_2 = 3\ell_h$. 
	In our own implementation of CBN \cite[Alg.~4]{khanh2024coderivative}
		with parameter fine-tuning, 
		we particularly set the proximal gradient step length $\tau = 1/(2\ell_h)$, the descent quantity coefficient $\sigma = 1/(20\ell_h)$, and the scaling factor of line search $\beta = 0.1$. The parameter setting of other algorithms will be given in the comparison part.

	\subsubsection{Convergence Test}
	
	In this subsection, our main goal is to examine the convergence performance of SGSN. The following numerical example generates various types of simulated datasets for convergence tests.
	
	\begin{example} \label{ex1}
		We first generate vectors $\bfmu_1, \bfmu_2 \in \mathbb{R}^n$ with independent and identically distributed (i.i.d.) elements being selected from standard normal distribution $N(0,1)$. Each element of diagonal covariance matrices $\Lambda_1, \Lambda_2 \in \mathbb{R}^{n\times n}$ is also i.i.d., following a folded normal distribution. Then positive and negative samples are normally distributed with $N(\bfmu_1,\Lambda_1)$ and $N(\bfmu_2,\Lambda_2)$ respectively.
		Let $q$ be the total number of samples and $p \in (0,1)$ be the proportion of positive samples, then $q_+ = \lceil pq \rceil$ and $q_- = q - q_+$. Given noise ratio $r \in (0,1)$, $\lfloor rq_+ \rfloor$ positive and negative samples are chosen to be marked with reverse labels respectively.  
	\end{example}
	
	We use the following violation of dual optimality ($\VDO$) to measure the closeness between a dual variable $\bfz$ and a P-stationary point of \eqref{01dco-uc}
	\begin{align*}
		\VDO : = {\rm dist} ( \bfz,~ \proxg ( \bfz - \tau \nabla h(\bfz) ) ) / \tau.
	\end{align*}
	We also record the value of objective function $F(\bfz)$ and computational time (\texttt{TIME}) of the algorithm.
	\begin{figure}[h]
		\subfigure{
			\begin{minipage}[t]{0.5\linewidth}
				\centering
				\includegraphics[width=3in]{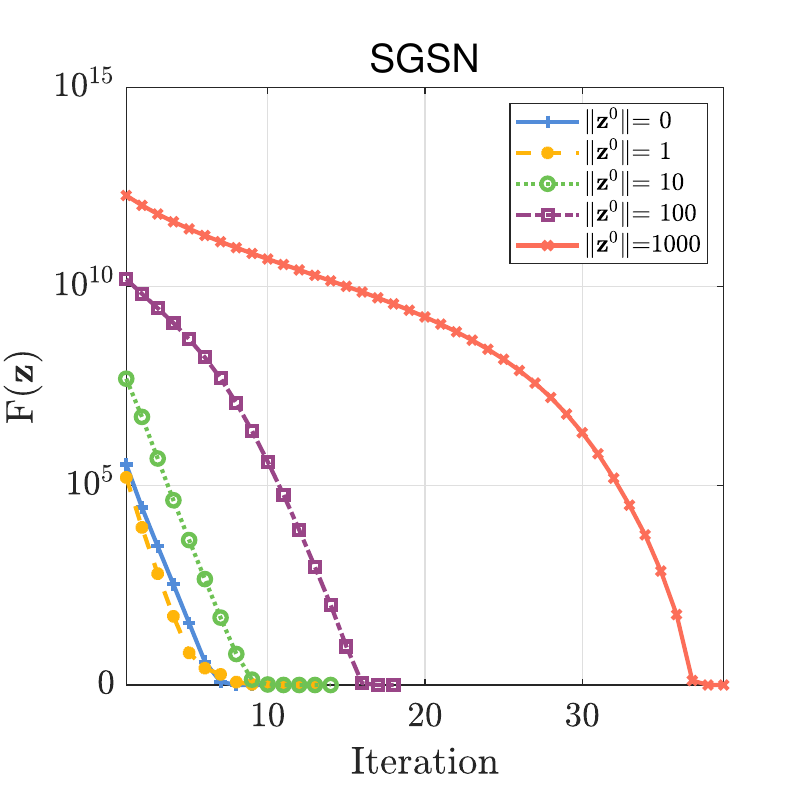}
			\end{minipage}%
		}
		\subfigure{
			\begin{minipage}[t]{0.5\linewidth}
				\centering
				\includegraphics[width=3in]{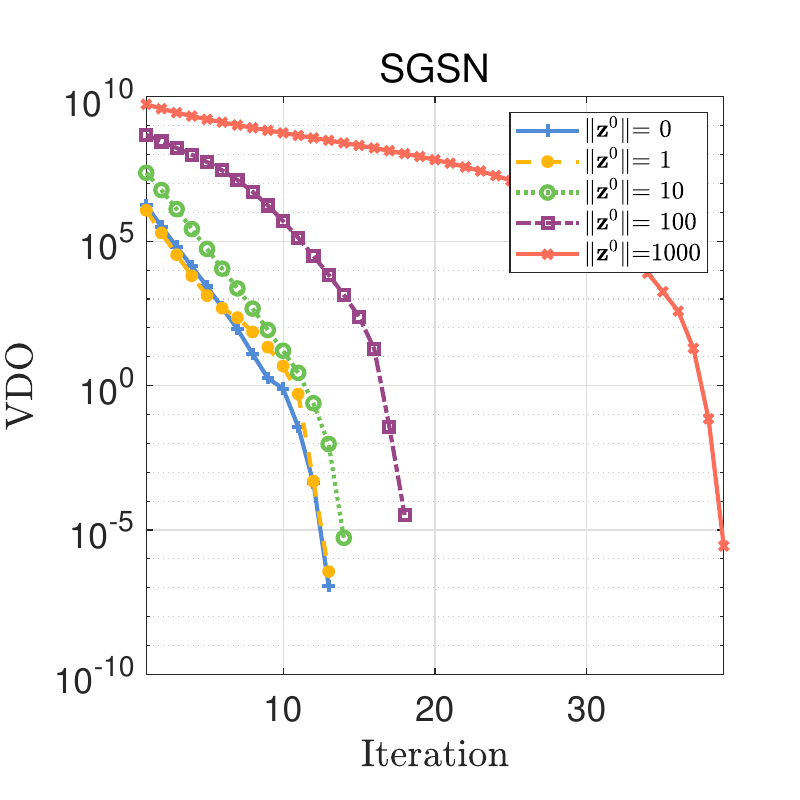}
			\end{minipage}%
		}
		\caption{Convergence test for SGSN with different initial points. 
			\label{fig-z0}}
		{}
	\end{figure}

	$\bullet$ \textbf{Test 0: On change of $\bfz^0$.} This test aims to observe the basic convergence property of SGSN with different initial points $\bfz^0$. It is conducted on a simulated dataset with $q = 1000$, $n = 100$, $p = 0.05$ and $r = 0$. To produce different initial points for SGSN, we first generate a vector $\overline{\bfz}^0 \in \mbR^{n}$ with entries being i.i.d. samples from the uniform distribution $U(0,1)$. 
	Then we compute $\bfz^0 = \rho \overline{\bfz}^0/\| \overline{\bfz}^0 \|$ with $\rho \in \{ 0,1,10,100,1000 \}$. SGSN starts with different $\bfz^0$ and the corresponding \texttt{VDO} and $F(\bfz)$ are recorded in Figure \ref{fig-z0}. The left panel in Figure \ref{fig-z0} demonstrates the sufficient descent of $\{ F(\bfz^k) \}_{k \in \mbN}$, which has been proved in Prop.~\ref{Prop-SGSN-property}. The right panel shows the descent of \VDO, with a faster rate during the last few iterations. This phenomenon illustrates the global and local quadratic convergence of $\{ \bfz^k \}_{k \in \mbN}$ given in Thms.~\ref{Thm-Global} and \ref{Thm-Local}.

	Next, we test SGSN with $\bfz^0 = 0$ on simulated datasets with various $q$ (total number of samples), $n$ (number of features), $p$ (proportion of positive samples), and $r$ (noise ratio). For comparison, we select NM01 and CBN to solve the primal problem \eqref{01cop} and the dual problem \eqref{01dco-uc} respectively. Their convergence performance can be measured by $\VDO$ and violation of primal optimality ($\VPO$) respectively, where $\VPO$ is defined as follows:
	\begin{align*}
		\VPO := \max \left\{ \| \nabla f(\bfx) + A^\top \bfz \|, ~{\rm dist} \left( A \bfx + \bfb, ~{\rm Prox}_{\xi \indf(\cdot) } ( A \bfx + \bfb + \xi \bfz ) \right) / \xi \right\},
	\end{align*}
	where $\xi$ is a proximal parameter adjusted as the defaulted setting of NM01. 
	
	\begin{figure}[h]
		\subfigure{
			\begin{minipage}[t]{0.5\linewidth}
				\centering
				\includegraphics[width=3in]{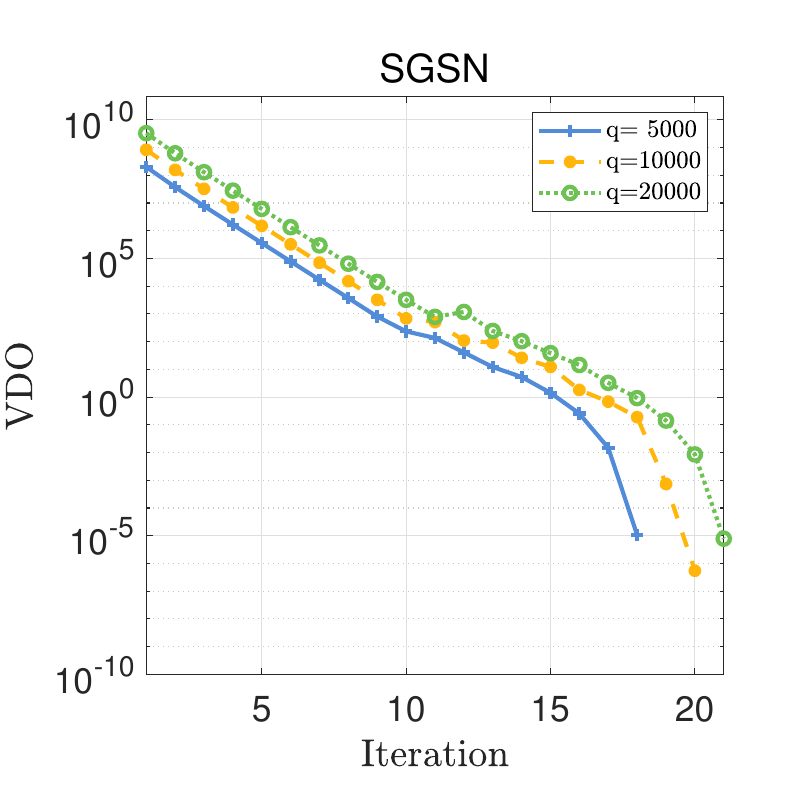}
			\end{minipage}%
			\begin{minipage}[t]{0.5\linewidth}
				\centering
				\includegraphics[width=3in]{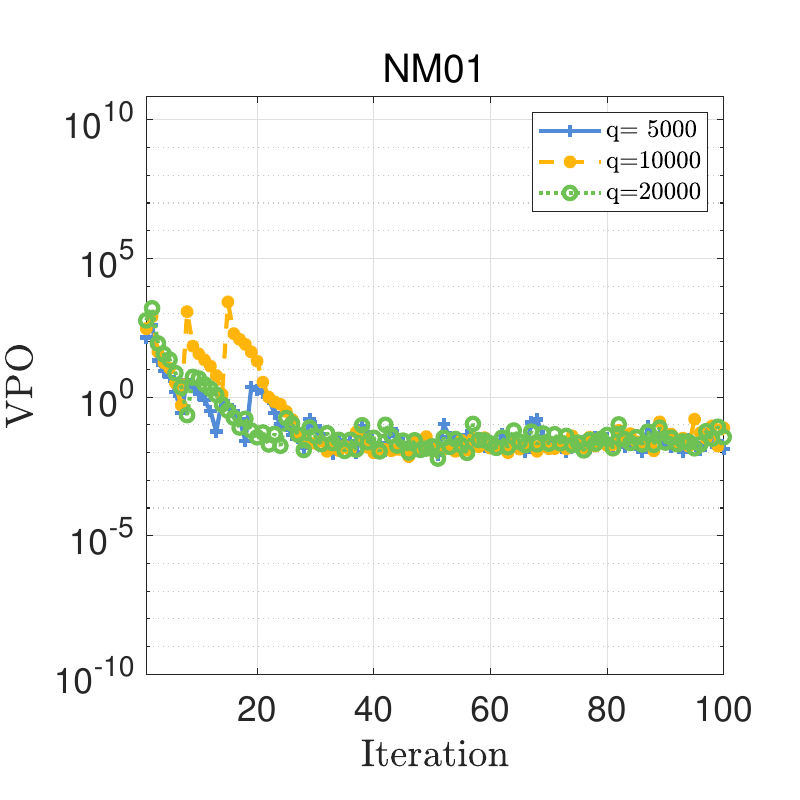}
		\end{minipage}}
		\subfigure{
			\begin{minipage}[t]{0.5\linewidth}
				\centering
				\includegraphics[width=3in]{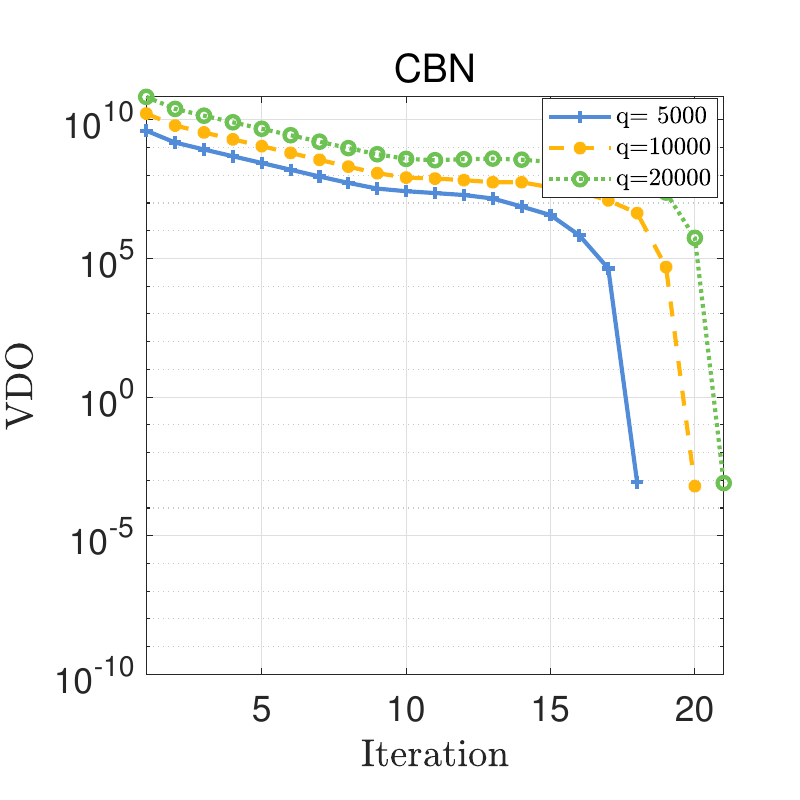}
			\end{minipage}%
			\begin{minipage}[t]{0.5\linewidth}
				\centering
				\includegraphics[width=3in]{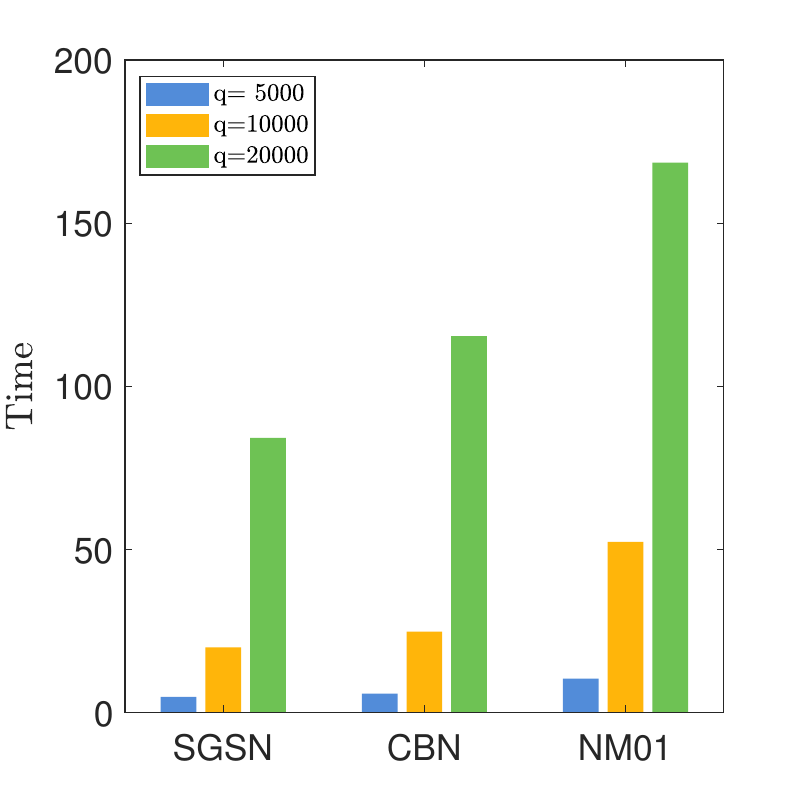}
		\end{minipage}}
		\caption{Convergence comparison on the simulated data with change of $q$. 
			\label{fig-q}}
		{}
	\end{figure}
	
	$\bullet$ \textbf{Test 1: On change of $q$.} In this test, we generate datasets with $q \in \{ 5000, 10000, 20000 \}$, $n = 100$, $p = 0.2$ and $r = 0$. In this setting, $m$ (the number of rows in matrix $A$) is taken from $\{ 4\times 10^6, 1.6\times 10^7,6.4\times 10^7 \}$. SGSN and CBN demonstrate better convergence performance when dealing with these large-scale matrices. From Figure \ref{fig-q}, we can also see that SGSN converges faster at the initial stage compared with CBN. At the final few steps, when the iterate is close to a P-stationary point of \eqref{01dco-uc}, $\VDO$ of SGSN rapidly drops below $10^{-4}$. 
		In comparison, the $\VPO$ of NM01 decreases relatively faster during the first several iterations and then stagnates around $10^{-2}$. In this test, SGSN spends the shortest \TIME to achieve the smallest $\VDO$.
	
	\begin{figure}[h]
		\subfigure{
			\begin{minipage}[t]{0.5\linewidth}
				\centering
				\includegraphics[width=3in]{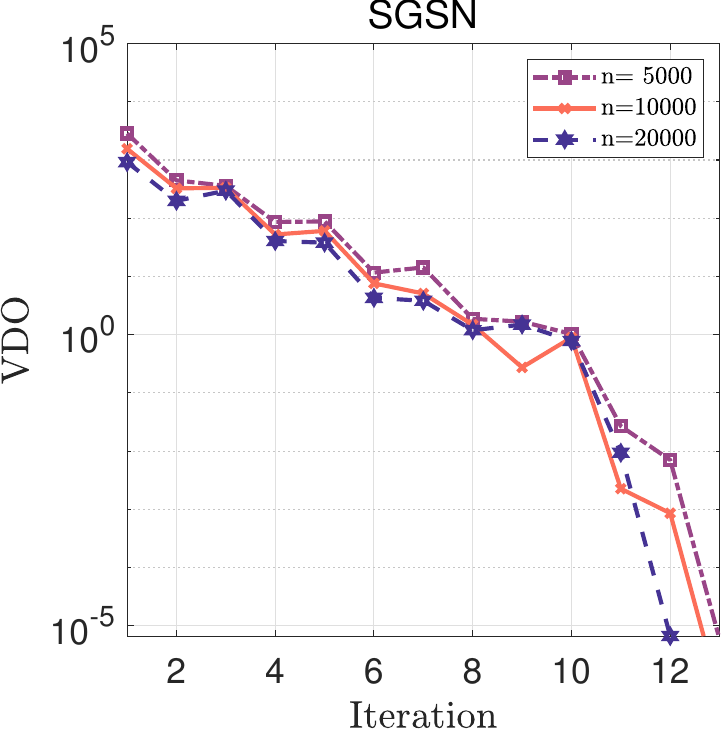}
			\end{minipage}%
			\begin{minipage}[t]{0.5\linewidth}
				\centering
				\includegraphics[width=3in]{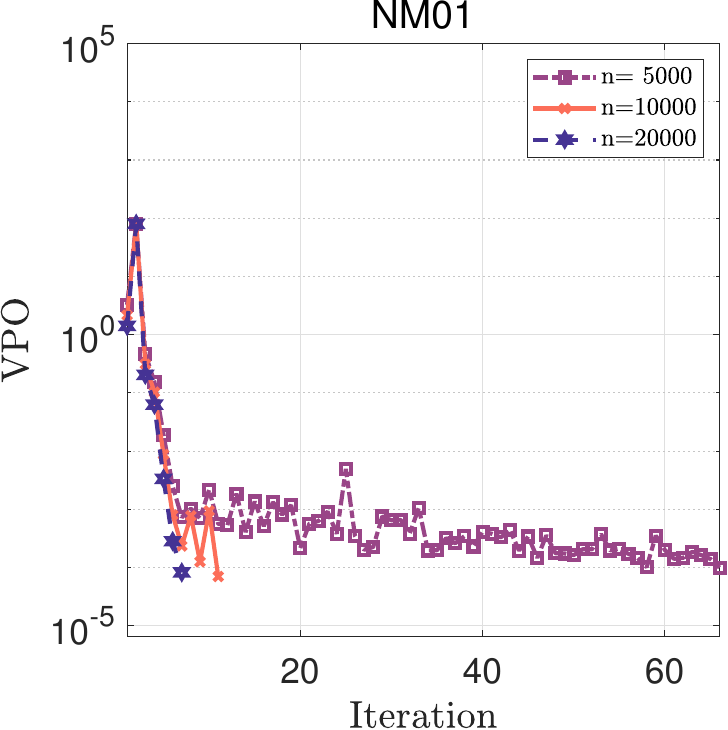}
		\end{minipage}}
		\subfigure{
			\begin{minipage}[t]{0.5\linewidth}
				\centering
				\includegraphics[width=3in]{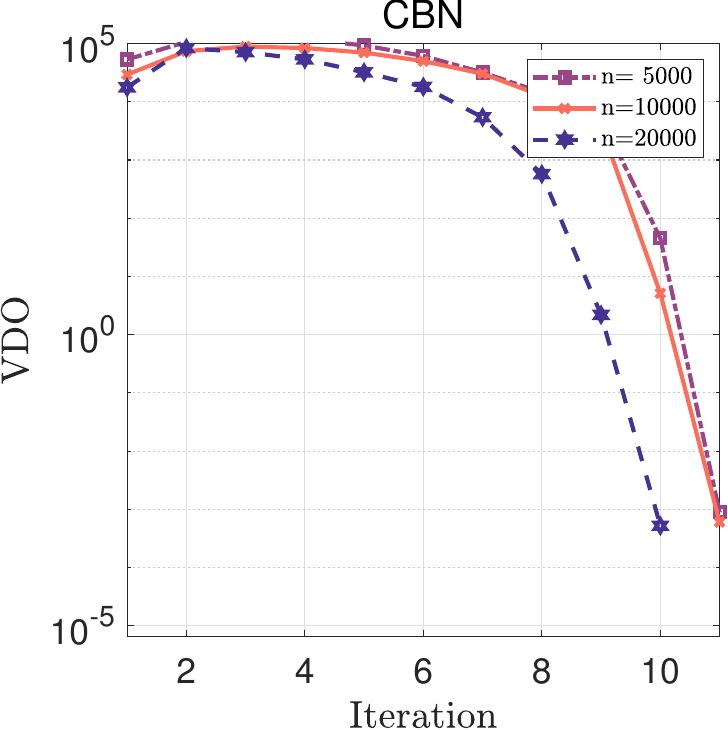}
			\end{minipage}%
			\begin{minipage}[t]{0.5\linewidth}
				\centering
				\includegraphics[width=3in]{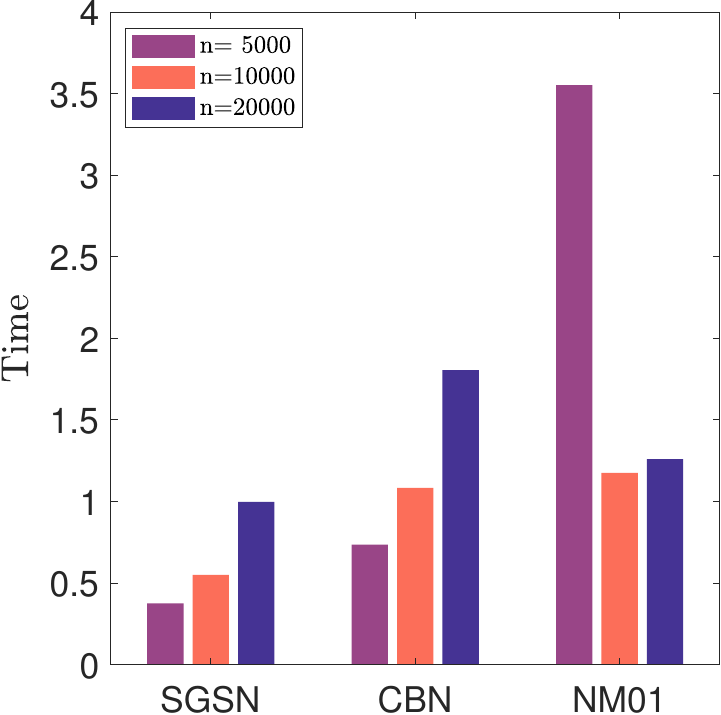}
		\end{minipage}}
		\caption{Convergence comparison on the simulated data with change of $n$. 
			\label{fig-n}}
		{}
	\end{figure}
	
	$\bullet$ \textbf{Test 2: On change of $n$.} We take $n \in \{ 5000, 10000, 20000 \}$, $q = 1000$, $p = 0.2$ and $r = 0$. Figure \ref{fig-n} shows that all the algorithms 
		converge very fast 
		in the cases of $n \geq 10000$. 
	When $n = 5000$, NM01 requires more iterations                                                       for $\VPO$ to get below $10^{-4}$. 
	
	\begin{figure}[h]
		\subfigure{
			\begin{minipage}[t]{0.5\linewidth}
				\centering
				\includegraphics[width=3in]{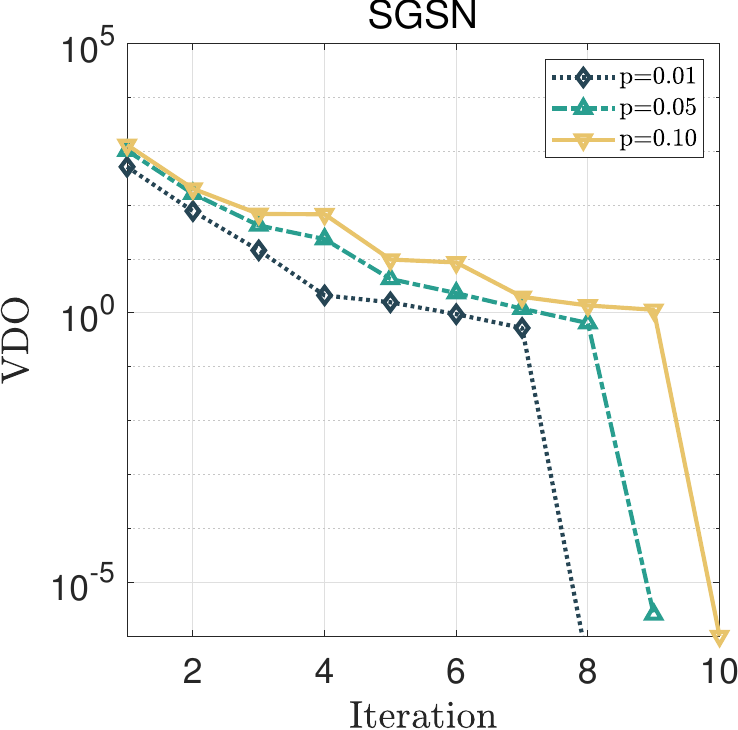}
			\end{minipage}%
			\begin{minipage}[t]{0.5\linewidth}
				\centering
				\includegraphics[width=3in]{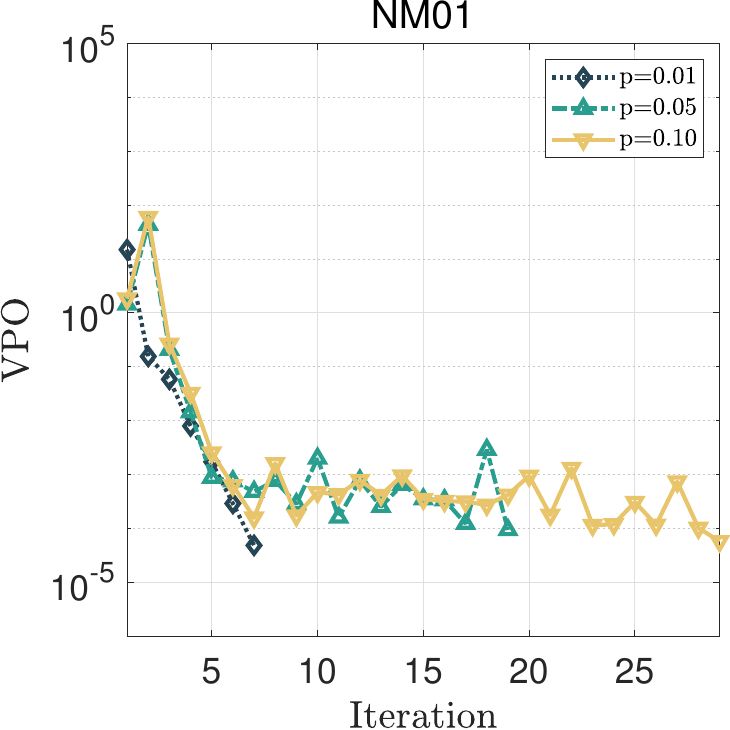}
		\end{minipage}}
		\subfigure{
			\begin{minipage}[t]{0.5\linewidth}
				\centering
				\includegraphics[width=3in]{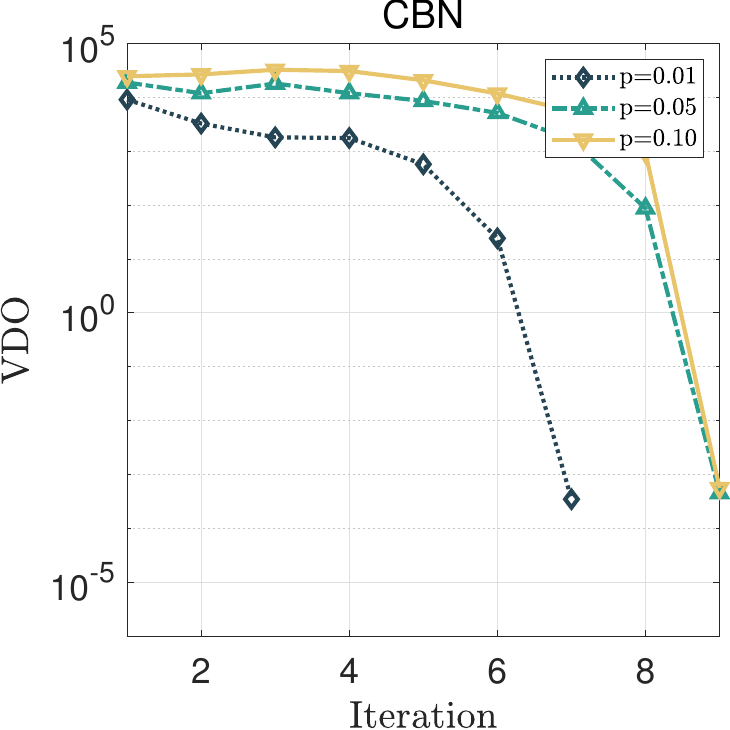}
			\end{minipage}%
			\begin{minipage}[t]{0.5\linewidth}
				\centering
				\includegraphics[width=3in]{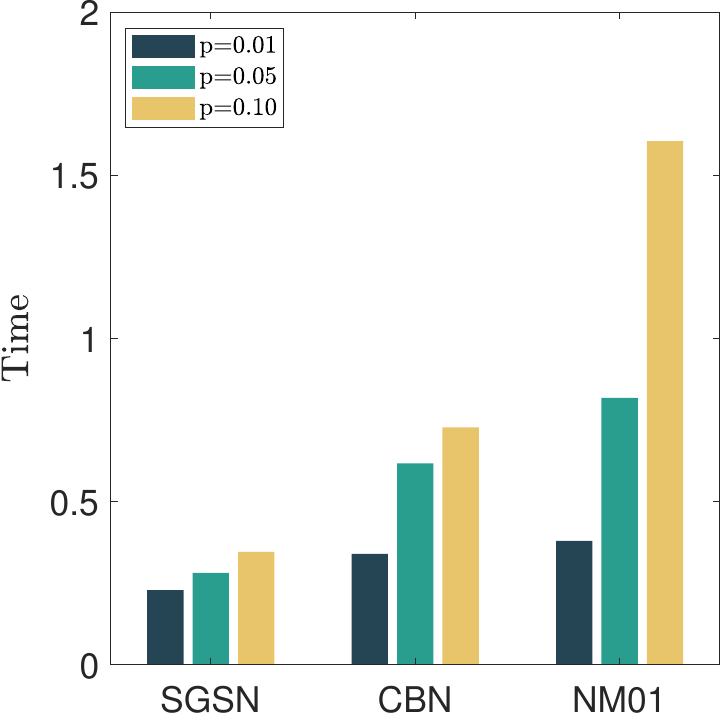}
		\end{minipage}}
		\caption{Convergence comparison on the simulated data with change of $p$. 
			\label{fig-p}}
		{}
	\end{figure}
	
	$\bullet$ \textbf{Test 3: On change of $p$.} We alter $p \in \{ 0.01, 0.05, 0.1 \}$, and fix $q = 1000$, $n = 10000$ and $r = 0$. 
	In all the cases, SGSN and CBN demonstrate significant decrease on $\VDO$ along with iteration. 
	We also observe from Figure \ref{fig-p} that NM01 spends more $\TIME$ and requires more iterations when $p = 0.05$ and $0.1$.
	
	\begin{figure}[h]
		\subfigure{
			\begin{minipage}[t]{0.5\linewidth}
				\centering
				\includegraphics[width=3in]{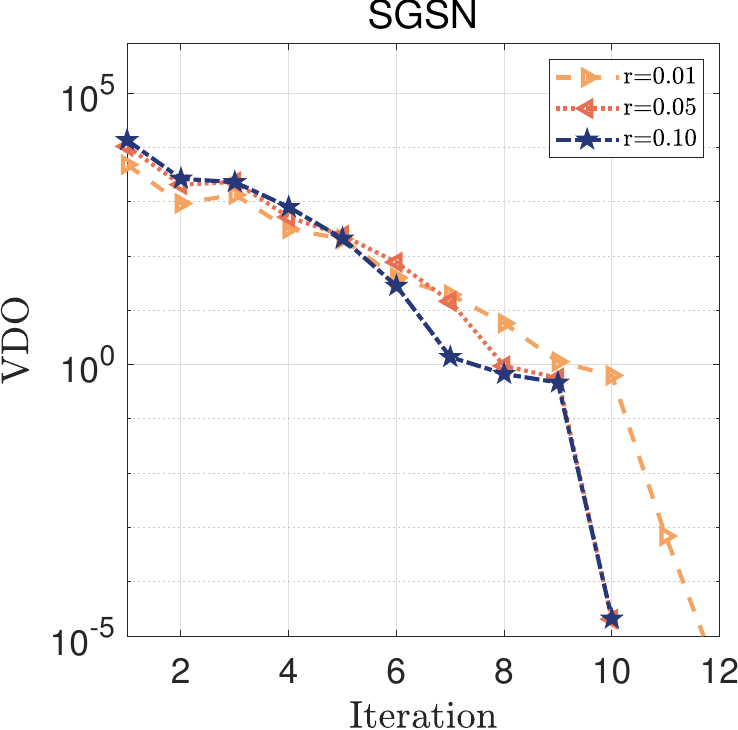}
			\end{minipage}%
			\begin{minipage}[t]{0.5\linewidth}
				\centering
				\includegraphics[width=3in]{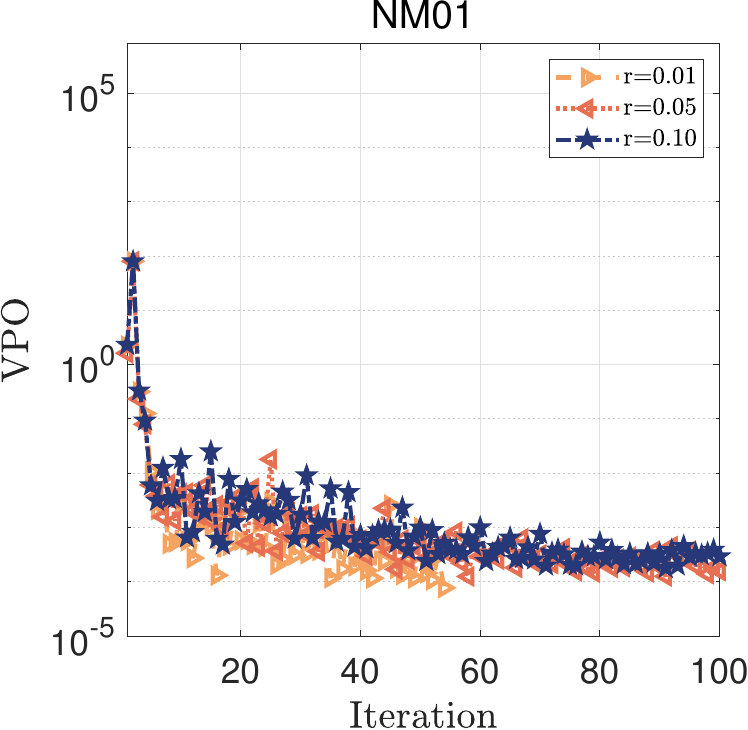}
		\end{minipage}}
		\subfigure{
			\begin{minipage}[t]{0.5\linewidth}
				\centering
				\includegraphics[width=3in]{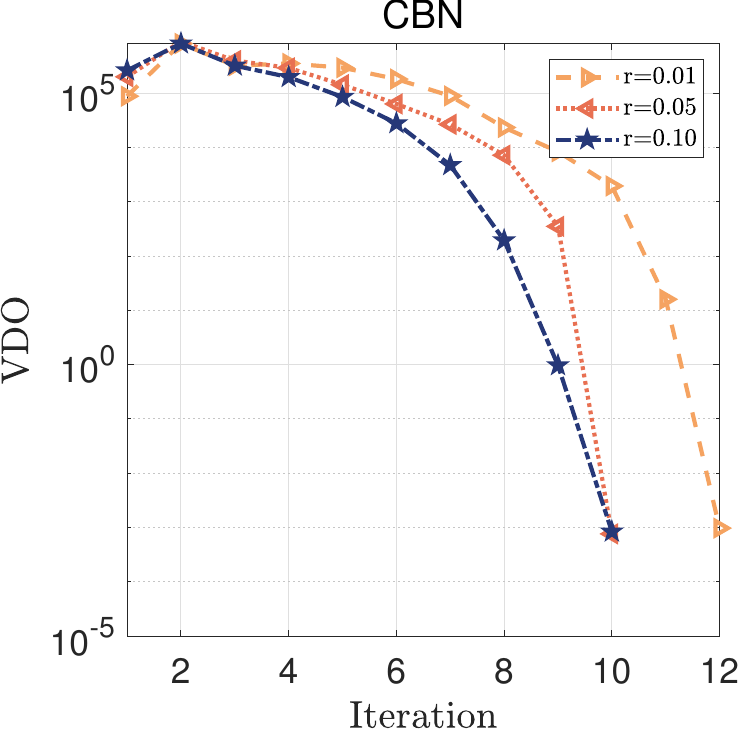}
			\end{minipage}
			\begin{minipage}[t]{0.5\linewidth}
				\centering
				\includegraphics[width=3in]{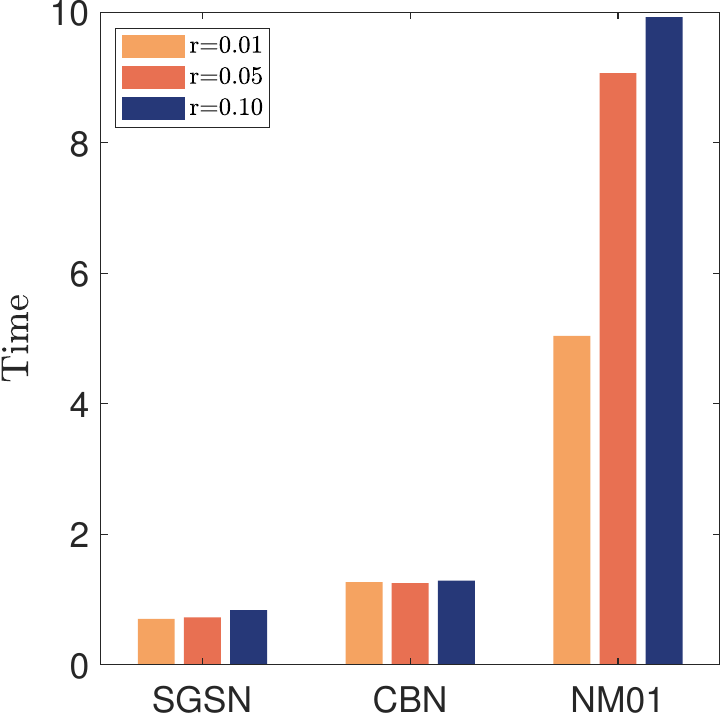}
		\end{minipage}}
		\caption{Convergence comparison on the simulated data with change of $r$. 
			\label{fig-r}}
		{}
	\end{figure}
	
	$\bullet$ \textbf{Test 4: On change of $r$.} Finally, we test cases of $r \in \{ 0.01, 0.05, 0.1 \}$, $q = 1000$, $n = 10000$ and $p = 0.2$.  As shown in Figure \ref{fig-r} again, compared with the other two algorithms, SGSN has the smallest $\VDO$ with the shortest $\TIME$. The $\VPO$ of NM01 decreases faster over the first several iterations, and then fluctuates around $10^{-3}$. 
	
	\subsubsection{Experiments on Real Dataset}
	
	In this part, we test SGSN and the three selected algorithms on 
	the real datasets given in Table \ref{tab-auc-dataset}. 
	For each dataset, five-fold cross-validation is conducted and the average value of $\AUC$ and $\TIME$ for each algorithm are recorded. 
	As for the parameters of the comparison algorithms, NM01 adopts the default setting. 
	The weighted parameter for the loss function in PRSVM and OPAUC is chosen from $\{ 10^{-3}, 10^{-2}, \cdots, 10^3 \}$ with the highest $\AUC$. 
	Moreover, we stop SGSN and CBN when $\texttt{VDO}_k / \texttt{VDO}_1 \leq 10^{-3} $ or $| \texttt{VDO}_k - \texttt{VDO}_{k -1}  | < 10^{-3}$, where $\VDO_k$ denotes the violation of dual optimality at the $k$-th iterate. NM01 is terminated by a similar rule, replacing the above $\VDO$ with $\VPO$.
	
	\begin{example}
		The details of real $\mbox{datasets}^{\mbox{\rm \scriptsize 1-4}}$\footnotetext[1]{https://jundongl.github.io/scikit-feature/}\footnotetext[2]{https://www.openml.org/}\footnotetext[3]{https://www.refine.bio/}\footnotetext[4]{https://www.csie.ntu.edu.tw/~cjlin/libsvmtools/datasets/} for AUC maximization are summarized in Table \ref{tab-auc-dataset}. All the datasets are preprocessed by feature-wise scaling to
		the interval  $[-1,1]$.
		\begin{table}[htb]
			\centering
			\caption{Real datasets for AUC maximization.} \label{tab-auc-dataset}
			\begin{tabular}{ccccccc}
				\hline
				Abbreviation & Dataset           & Domain     & $q_-$ & $q_+$ & $ m$ & $n$   \\ \hline
				\texttt{all} & Allaml            & Text       & 47    & 25    & 1175                  & 7129  \\
				\texttt{aou} & AP\_Ovary\_Uterus & Biology    & 124   & 198   & 24552                 & 10936 \\
				\texttt{bas} & Basehock          & Text       & 999   & 994   & 993006                & 4862  \\
				\texttt{col} & Colon\_Cancer     & Biology    & 40    & 22    & 880                   & 2000  \\
				\texttt{duk} & Duke              & Biology    & 21    & 23    & 483                   & 7129  \\
				\texttt{gli} & Gli85             & Biology    & 26    & 59    & 1534                  & 22283 \\
				\texttt{gs2} & Gse2280           & Biology    & 5     & 22    & 110                   & 11868 \\
				\texttt{gs7} & Gse7670           & Biology    & 27    & 27    & 729                   & 11868 \\
				\texttt{leu} & Leukemia          & Biology    & 47    & 25    & 1175                  & 7070  \\
				\texttt{mad} & Madelon           & Artificial & 1000  & 1000  & 1000000               & 500   \\
				\texttt{ove} & Ova\_Endometrium  & Biology    & 61    & 1484  & 90524                 & 10936 \\
				\texttt{ovk} & Ova\_Kidney       & Biology    & 260   & 1285  & 334100                & 10936 \\
				\texttt{ovo} & Ova\_Ovary        & Biology    & 198   & 1347  & 266706                & 10936 \\
				\texttt{pcm} & Pcmac             & Text       & 961   & 982   & 943702                & 3289  \\
				\texttt{pro} & Prostate\_Ge       & Biology    & 50    & 52    & 2600                  & 5966  \\
				\texttt{rel} & Relathe           & Text       & 648   & 779   & 504792                & 4322  \\
				\texttt{smk} & Smk\_Can187       & Biology    & 90    & 97    & 8730                  & 19993 \\
				\texttt{spl} & Splice            & Biology    & 483   & 517   & 249711                & 60    \\ \hline
			\end{tabular}
		\end{table}
	\end{example} 
	
	The comparison results are summarized in Table \ref{res-auc-real}. We can observe that SGSN has a similar $\AUC$ but shorter \texttt{TIME} than NM01 or CBN because it is more effective in reducing the violation of optimality.
	Although the other two algorithms are also competitive on $\AUC$, 
	their $\TIME$ are longer than that of SGSN on most of the datasets. 
	For example, on the dataset \texttt{ovk}, our SGSN only requires about $10\%$ of the $\TIME$ of PRSVM to achieve a higher $\AUC$.

	\begin{sidewaystable}[]
		\centering 
		\caption{Numerical results on real datasets for AUC maximization.} \label{res-auc-real}
			{\begin{tabular}{c|ccccc|ccccc}
					\hline
					& \multicolumn{5}{c|}{\texttt{AUC}~(\%)~$\uparrow$} & \multicolumn{5}{c}{\texttt{TIME}~(sec)~$\downarrow$} \\
					\hline
					& SGSN  & NM01  & CBN   & PRSVM & OPAUC & SGSN     & NM01     & \multicolumn{1}{c}{CBN} & PRSVM    & OPAUC    \\ \hline
					\texttt{all} & 99.00 & 99.00 & 99.00 & 98.00 & 98.50 & 2.677e-2 & 3.285e-2 & 3.980e-2                & 1.438e-1 & 1.961e+1 \\
					\texttt{aou} & 94.48 & 94.75 & 94.46 & 94.18 & 94.16 & 2.773e-1 & 3.585e-1 & 4.819e-1                & 2.821e+0 & 1.950e+2 \\
					\texttt{bas} & 98.82 & 97.82 & 98.46 & 98.73 & 98.76 & 1.868e+0 & 1.017e+1 & 3.077e+0                & 4.701e+0 & 2.335e+2 \\
					\texttt{col} & 87.92 & 85.63 & 87.92 & 86.04 & 85.83 & 2.061e-2 & 1.085e-1 & 3.083e-2                & 7.766e-2 & 3.407e+0 \\
					\texttt{duk} & 89.54 & 90.79 & 89.54 & 89.54 & 84.54 & 3.082e-2 & 4.377e-2 & 4.797e-2                & 1.517e-1 & 3.031e+1 \\
					\texttt{gli} & 95.76 & 95.76 & 95.76 & 94.18 & 93.33 & 8.698e-2 & 1.100e-1 & 1.422e-1                & 4.148e-1 & 2.701e+2 \\
					\texttt{gs2} & 61.67 & 61.67 & 61.67 & 61.67 & 56.67 & 1.556e-2 & 1.513e-2 & 2.031e-2                & 6.896e-2 & 1.999e+1 \\
					\texttt{gs7} & 98.40 & 98.40 & 98.40 & 98.40 & 97.60 & 4.164e-2 & 3.406e-2 & 7.484e-2                & 1.333e-1 & 3.964e+1 \\
					\texttt{leu} & 99.00 & 99.00 & 99.00 & 98.50 & 65.00 & 2.344e-2 & 2.690e-2 & 4.478e-2                & 1.665e-1 & 1.951e+1 \\
					\texttt{mad} & 55.68 & 53.66 & 55.41 & 54.33 & 58.26 & 5.120e+0 & 1.300e+1 & 1.134e+1                & 3.154e+0 & 9.349e+0 \\
					\texttt{ove} & 95.81 & 94.52 & 95.21 & 95.42 & 93.69 & 7.186e-1 & 2.379e+0 & 1.492e+0                & 2.012e+1 & 1.209e+3 \\
					\texttt{ovk} & 99.16 & 98.95 & 99.11 & 99.09 & 98.26 & 1.251e+0 & 2.151e+0 & 2.314e+0                & 1.844e+1 & 1.288e+3 \\
					\texttt{ovo} & 93.30 & 91.74 & 93.32 & 93.16 & 91.52 & 1.517e+0 & 2.367e+0 & 2.375e+0                & 4.341e+1 & 1.263e+3 \\
					\texttt{pcm} & 91.11 & 88.73 & 90.56 & 90.40 & 90.26 & 1.669e+0 & 2.782e+1 & 5.701e+0                & 8.364e+0 & 1.468e+2 \\
					\texttt{pro} & 96.00 & 95.60 & 96.00 & 95.80 & 90.33 & 3.279e-2 & 3.787e-1 & 5.935e-2                & 1.196e-1 & 2.401e+1 \\
					\texttt{rel} & 93.28 & 90.88 & 92.14 & 92.57 & 92.44 & 1.336e+0 & 4.266e+0 & 6.050e+0                & 3.166e+0 & 1.456e+2 \\
					\texttt{smk} & 78.76 & 77.93 & 77.81 & 78.04 & 76.49 & 1.217e+0 & 1.791e+0 & 4.312e-1                & 9.226e-1 & 4.819e+2 \\
					\texttt{spl} & 88.28 & 88.39 & 88.23 & 88.23 & 87.19 & 1.909e-1 & 5.908e-1 & 5.429e-1                & 1.597e-1 & 1.845e-1 \\ \hline
			\end{tabular}}
			{\\``$\uparrow$" means that a larger metric value corresponds to better algorithm performance, while ``$\downarrow$" has the opposite meaning.}
		\end{sidewaystable}
		
		\subsection{Sparse Multi-Label Classification} 
		For a multi-label classification (MLC) problem with $\ell$ labels, let $\bfc_i \in \mbR^d$, $i \in [q]$ be the feature vector with its last entry $[\bfc_i]_d = 1$ and class label $\bfy_i \in \{ -1,1 \}^\ell$ for $i \in [q]$. Then we denote matrices $\textbf{C}:= [\bfc_1,\cdots,\bfc_q]^\top \in \mbR^{q\times d}$ and $\textbf{Y}:= [\bfy_1,\cdots,\bfy_q]^\top \in \mbR^{q\times \ell}$. Given arbitrary feature vector $\bfc$, a linear binary relevance classifier $\H:\mbR^d \to \mbR^\ell$ returns a label vector by the following rule \cite{boutell2004learning,wu2018unified}:
		\begin{align*}
			\H(\bfc) := \left( {\rm sgn}(\bfc^\top \bfx_1), \cdots,  {\rm sgn}(\bfc^\top \bfx_\ell)\right),
		\end{align*}
		where $\bfx := (\bfx_1;\cdots;\bfx_\ell) \in \mbR^{d\ell}$ is a weighted vector to be 
		determined through optimizing certain loss function,
		and ${\rm sgn} : \mbR \to \mbR$ is the sign function, returning $1$ when the argument variable is positive and $-1$ otherwise. The Hamming loss (\texttt{HL}, see \cite{zhang2013review}) is an important evaluation metric in MLC. Let $\bfy^{(k)} = ( y^{(k)}_1, \cdots, y^{(k)}_q ) \in \mbR^q$ be the $k$-th column of the matrix $\textbf{Y}$, then Hamming loss for classifier $\H$ is defined as follows:
		\begin{align*}
			\texttt{HL} &:= \frac{1}{q\ell} \sum_{i=1}^{q} \sum_{k=1}^{\ell} \bfone_{(0,\infty) } ( -y_i^{(k)} \bfc_i^\top \bfx_k).
		\end{align*} 
		When optimizing $\texttt{HL}$ on the training data, a practical technique is to replace the above $-y_i^{(k)} \bfc_i^\top \bfx_k$ with $-y_i^{(k)} \bfc_i^\top \bfx_k + 1$, see e.g. \cite{wu2018unified}. This helps to avoid trivial solution and improve generalization ability of the classifier. Meanwhile, when the dimension of the parameter $\bfx \in \mbR^{d\ell}$ is large, we can use sparse regularization to alleviate the risk of overfitting and reduce the dimension. Overall, we consider minimizing the Hamming loss with an elastic net regularizer. Denoting $m:= q\ell$ and $n:= d\ell$, this model is a special case of \eqref{01cop} with the following setting:
		\begin{align*}
			f(\bfx) = \frac{1}{2} \| \bfx \|^2 + \lambda_1 \| \bfx \|_1,~
			A := \left[\begin{array}{ccc}
				-( \bfy^{(1)} \bfe^\top_d) \odot \textbf{C} & &\\
				& \ddots &  \\
				&  & -( \bfy^{(\ell)} \bfe^\top_d) \odot \textbf{C} 
			\end{array} \right],~ \bfb = \bfe_m,
		\end{align*}  
		where $\| \cdot \|_1$ is the $\ell_1$ norm and $\odot$ is the Hadamard product.

		We can compute the conjugate function
		\begin{align*}
			f^*(\bfv) = \sum_{i = 1}^{n} \varphi ( v_i ) \quad \hbox{~with~} \ \varphi(t) := \left\{  \begin{aligned}
				& 0, &&\mbox{~if~} |t| \leq \lambda_1, \\
				&( | t | - \lambda_1 )^2/2, &&\mbox{~if~} |t| > \lambda_1,
			\end{aligned} \right.
		\end{align*}
		and thus the corresponding gradient and generalized Jacobian are as follows:
		\begin{align*}
			&\nabla f^*(\bfv) = \left\{ \bfr \in \mbR^n: r_i = \left\{ \begin{aligned}
				& v_i - \lambda_1, &&\mbox{~if~} v_i > \lambda_1 \\
				& 0, &&\mbox{~if~} |v_i| \leq \lambda_1 \\
				& v_i + \lambda_1, &&\mbox{~if~} v_i < -\lambda_1
			\end{aligned} \right. \right\} \\
			& \partial^2 f^* (\bfv) = \left\{ \mbox{Diag}(\bfr): \bfr \in \mbR^n, r_i \in \left\{ \begin{aligned}
				& 1, &&\mbox{~if~} |v_i| > \lambda_1 \\
				& [0,1], &&\mbox{~if~} |v_i| = \lambda_1 \\
				& 0, &&\mbox{~if~} |v_i| < \lambda_1
			\end{aligned} \right.   \right\} .
		\end{align*}
		Let $\bfa^{(i)}$ be the $i$-th column of $A$, our SGSN will solve the following 
		dual problem: 
		\begin{align}\label{mlc-dual}
			\min_{\bfz \in \mbR^{m}}\ F(\bfz) = \sum_{i = 1}^{n} \varphi(-\langle \bfa^{(i)}, \bfz \rangle) - \langle \bfb,  \bfz \rangle + g(\bfz) .
		\end{align} 
		Next, we discuss the convergence property when solving \eqref{mlc-dual}. 
		One can verify that the elastic net regularizer is semi-algebraic by 
		\cite[Def.~5]{bolte2014proximal}. 
		Assuming the sequence $\{ \bfz^k \}_{k \in \mbN}$ is bounded, the whole sequence converges to a P-stationary point $\bfz^*$ of \eqref{mlc-dual} according to Prop.~\ref{Prop-SGSN-property}. $\nabla f^*$ is semismooth. Denoting $\S^*:=\S(\bfz^*)$ and $\Gamma^*:=\{ i \in [n]: | \langle \bfa^{(i)}, \bfz^* \rangle | > \lambda_1 \}$, if we further assume that $A_{\S^*\Gamma^*}$ has full row rank, then for any $H \in \widehat{\partial}^2 h( \bfz^* )$, we have $\bfd^\top H_{\S^*} \bfd = \| A^\top_{\S^*\Gamma^*} \bfd \|^2 >  0$ for any $\bfd \in \mbR^{| \S^* |}\backslash \{ 0 \}$. This means Assumption \ref{semismooth+sosc} holds and $\{ \bfz^k \}_{k \in \mbN}$ converges to $\bfz^*$ with local superlinear rate by Thm.~\ref{Thm-Local}. 
		
		In the following numerical experiments, we set $\bfz^0 = 0$, $\mu = 10^{-5}$, $\gamma = 10^{-2}$, $c_1 = 10^{-4}$ and $c_2 = 10^8$ for SGSN. 
		The proximal parameter is adaptively chosen as $\tau_k = 0.1^j$, where $j$ is the smallest natural number such that $F(\bfv^k) - F(\bfz^k) \geq 10^{-4}\| \bfv^k - \bfz^k \|^2$. The regularization parameter $\lambda_1$ is selected from $\{ 2^{-6}, 2^{-5}, \cdots, 2^{6} \}$. Other setting for SGSN is the same as that in the last subsection. 
		We select three comparison algorithms for the MLC problem with $\ell_1$ norm or group sparse regularizer: LLFS from \cite{huang2016learning}, MDFS from \cite{zhang2019manifold} and RFS from \cite{nie2010efficient}. 
		We adjust the sparse regularization parameters of the three algorithms to achieve better performance. They are selected from $\{ 10^{-6}, 10^{-5}, \cdots, 10^{6} \}$ 
		for RFS and MDFS, and $\{ 2^{-6}, 2^{-5}, \cdots, 2^{6} \}$ for LLFS. 
		Other parameters of the three algorithms are set the same as their default values. 
		We use the Hamming loss ($\texttt{HL}$), computational time (\texttt{TIME}) and number of nonzero elements (\texttt{NNE}) of the solution to evaluate
		the performance of the algorithms.   
		
		\subsubsection{Experiments on Simulated Dataset}
		
		To observe how the change of dimension influence the performance, 
		we test all four algorithms on simulated datasets with various $q$, $d$ and $\ell$ generated as follows.
		
		\begin{example} \label{ex-mlc-simu}
			We first generate a matrix $\overline{\textbf{C}} \in \mbR^{q \times (d-1)}$ with i.i.d elements being chosen from standard normal distribution $N(0,1)$, and then the feature matrix is computed by $\textbf{C} = [ \overline{\textbf{C}}, \bfe_q ]$. Next, we generate a matrix $W \in \mbR^{d \times \ell}$ with each element chosen from the uniform distribution $U(-1,1)$, and then the label matrix is computed by $\textbf{Y} = {\rm Sgn}(\textbf{C}W)$, where ${\rm Sgn}(\cdot)$ is the sign function for a matrix. Finally, $90\%$ of samples are randomly selected as training data and the rest is regarded as testing data. 
		\end{example}
		
		\begin{table}[htb]
			\centering
			\caption{Numerical results on Example \ref{ex-mlc-simu} with different $q$, $d$ and $\ell$.} \label{tab2}
				\setlength{\tabcolsep}{0.8mm}
				{\begin{tabular}{ccccccccccccc}
						\hline
						\multicolumn{1}{c|}{}       & \multicolumn{4}{c|}{\texttt{HL} $\downarrow$}       & \multicolumn{4}{c|}{\texttt{TIME} (sec) $\downarrow$}     & \multicolumn{4}{c}{\texttt{NNE} $\downarrow$} \\ \cline{2-13} 
						\multicolumn{1}{c|}{}       & SGSN  & RFS   & LLFS  & \multicolumn{1}{c|}{MDFS}  & SGSN  & RFS   & LLFS & \multicolumn{1}{c|}{MDFS}                       & SGSN      & RFS       & LLFS      & MDFS      \\ \hline
						\multicolumn{1}{c|}{$q$}                         & \multicolumn{12}{c}{$d = 6000, \ell = 3$}                                                                                                               \\ \hline
						\multicolumn{1}{c|}{1000}   & 0.424 & 0.424 & 0.445 & \multicolumn{1}{c|}{0.453} & 1.186 & 4.908 & 22.10  & \multicolumn{1}{c|}{23.47} & 2874      & 7293      & 2750      & 2970      \\
						\multicolumn{1}{c|}{2000}   & 0.397 & 0.419 & 0.394 & \multicolumn{1}{c|}{0.434} & 2.295 & 19.82 & 32.94 & \multicolumn{1}{c|}{42.41} & 4625      & 12870     & 5312      & 4716      \\
						\multicolumn{1}{c|}{3000}   & 0.368 & 0.423 & 0.373 & \multicolumn{1}{c|}{0.410} & 4.371 & 49.11 & 42.54 & \multicolumn{1}{c|}{75.31} & 6097      & 15174     & 7620      & 6192      \\
						\multicolumn{1}{c|}{4000}   & 0.339 & 0.428 & 0.340 & \multicolumn{1}{c|}{0.402} & 7.818 & 93.03 & 52.79 & \multicolumn{1}{c|}{126.8} & 7470      & 16235     & 9679      & 7560      \\
						\multicolumn{1}{c|}{5000}   & 0.310 & 0.431 & 0.321 & \multicolumn{1}{c|}{0.390} & 11.21 & 154.0 & 64.08 & \multicolumn{1}{c|}{200.0} & 8524      & 16797     & 12726     & 8586      \\
						\multicolumn{1}{c|}{6000}   & 0.300 & 0.428 & 0.319 & \multicolumn{1}{c|}{0.395} & 14.15 & 239.6 & 86.30  & \multicolumn{1}{c|}{295.3} & 9303      & 17161     & 17046     & 9378      \\ \hline
						\multicolumn{1}{c|}{$d$}    & \multicolumn{12}{c}{$q = 2500, \ell = 3$}                                                                                                               \\ \hline
						\multicolumn{1}{c|}{2500}   & 0.304 & 0.438 & 0.331 & \multicolumn{1}{c|}{0.376} & 2.809 & 19.37 & 8.118 & \multicolumn{1}{c|}{20.14} & 4205      & 6362      & 5593      & 4237      \\
						\multicolumn{1}{c|}{4000}   & 0.346 & 0.420 & 0.354 & \multicolumn{1}{c|}{0.406} & 2.880  & 25.36 & 16.23 & \multicolumn{1}{c|}{33.98} & 4942      & 9772      & 4988      & 5004      \\
						\multicolumn{1}{c|}{5500}   & 0.375 & 0.422 & 0.370 & \multicolumn{1}{c|}{0.419} & 3.235 & 30.94 & 29.96 & \multicolumn{1}{c|}{49.35} & 5290      & 13152     & 5476      & 5362      \\
						\multicolumn{1}{c|}{7000}   & 0.401 & 0.426 & 0.404 & \multicolumn{1}{c|}{0.437} & 2.857 & 36.96 & 50.58 & \multicolumn{1}{c|}{69.50}  & 5416      & 16499     & 5755      & 5544      \\
						\multicolumn{1}{c|}{8500}   & 0.406 & 0.431 & 0.401 & \multicolumn{1}{c|}{0.437} & 3.311 & 41.69 & 77.83 & \multicolumn{1}{c|}{95.52} & 5830      & 19755     & 6035      & 5967      \\
						\multicolumn{1}{c|}{10000}  & 0.407 & 0.433 & 0.418 & \multicolumn{1}{c|}{0.451} & 3.761 & 47.29 & 106.8 & \multicolumn{1}{c|}{128.6} & 6147      & 23134     & 6322      & 6330      \\ \hline
						\multicolumn{1}{c|}{$\ell$} & \multicolumn{12}{c}{$q = 1000, d = 2000$}                                                                                                               \\ \hline
						\multicolumn{1}{c|}{5}      & 0.370 & 0.462 & 0.372 & \multicolumn{1}{c|}{0.403} & 1.176 & 2.186 & 4.604 & \multicolumn{1}{c|}{2.756} & 4207      & 5364      & 4732      & 4270      \\
						\multicolumn{1}{c|}{7}      & 0.372 & 0.468 & 0.378 & \multicolumn{1}{c|}{0.415} & 1.698 & 2.252 & 5.854 & \multicolumn{1}{c|}{2.453} & 5872      & 8014      & 6455      & 5950      \\
						\multicolumn{1}{c|}{9}      & 0.373 & 0.450 & 0.374 & \multicolumn{1}{c|}{0.398} & 2.117 & 2.26  & 7.887 & \multicolumn{1}{c|}{2.428} & 7528      & 10855     & 8030      & 7614      \\
						\multicolumn{1}{c|}{11}     & 0.371 & 0.441 & 0.375 & \multicolumn{1}{c|}{0.403} & 2.479 & 2.294 & 9.458 & \multicolumn{1}{c|}{2.069} & 9166      & 13749     & 9832      & 9284      \\
						\multicolumn{1}{c|}{13}     & 0.369 & 0.463 & 0.379 & \multicolumn{1}{c|}{0.407} & 3.143 & 2.389 & 10.82 & \multicolumn{1}{c|}{2.191} & 10872     & 16795     & 11601     & 10972     \\
						\multicolumn{1}{c|}{15}     & 0.364 & 0.467 & 0.374 & \multicolumn{1}{c|}{0.403} & 3.503 & 2.378 & 12.78 & \multicolumn{1}{c|}{2.092} & 12560     & 19877     & 13194     & 12720     \\ \hline
				\end{tabular}}
				{``$\downarrow$" means that a smaller metric value corresponds to better algorithm performance.}
			\end{table}
			
			$\bullet$ \textbf{On change of $q$}. We fix $d = 6000$, $\ell = 3$, and
			select $q \in \{ 1000, 2000, \cdots, 6000 \}$. 
			We can observe from Table \ref{tab2} that SGSN takes the shortest \texttt{TIME} in all the cases. 
			In particular, SGSN is at least $6$ times faster than
			the second fastest algorithm LLFS. 
			SGSN also finds the solution with the smallest \texttt{NNE} when $q > 1000$. \texttt{TIME} and \texttt{NNE} increases when $q$ rises, whereas \texttt{HL} shows the opposite trend. 
			
			$\bullet$ \textbf{On change of $d$}. We fix $q = 2500$, $\ell = 3$, and select $d \in \{ 2500, 4000, \cdots, 10000 \}$. 
			In this test, SGSN achieves the best performance in all the cases. When feature dimension $d$ grows, SGSN has a relatively small increase on \texttt{TIME}. 
			This indicates that SGSN may be more advantageous when applied to large-scale datasets.  
			
			$\bullet$ \textbf{On change of $\ell$}. We fix $q = 1000$, $d = 2000$,
			and select $\ell \in \{ 5,7,\cdots,15 \}$. From Table \ref{tab2}, we can observe that SGSN has the best \texttt{HL} and \texttt{NNE}. Its \texttt{{TIME}} is also the shortest when $\ell \leq 9$, and slightly larger than that of RFS and MDFS when $\ell \geq 11$.
			
			\subsubsection{Experiments on Real Data} In this part, 
			we test all the algorithms on real MLC problems listed in Table \ref{mlc_realdata}.

			\begin{example}
				The real datasets in Table \ref{mlc_realdata} are available on the multi-label classification dataset repository\footnote[5]{https://www.uco.es/kdis/mllresources/}. For each dataset, the partition in $67\%$ train and $33\%$ test is performed by following the Iterative Stratification method proposed by \cite{sechidis2011stratification}. 
				\begin{table}[htb]
					\centering
					\caption{Real multi-label classification datasets with higher dimension.} \label{mlc_realdata}	
						{\begin{tabular}{ccccccc}
								\hline
								Abbreviation & Dataset             & Domain  & $q$   & $d$   & $\ell$ & $q\times d \times \ell$ \\ \hline
								\texttt{bbc} & 3s-bbc1000          & Text    & 352   & 1000  & 6      & 2.11e+06                \\
								\texttt{int} & 3s-inter3000        & Text    & 169   & 3000  & 6      & 3.04e+06                \\
								\texttt{eug} & EukaryoteGO         & Biology & 7766  & 12689 & 22     & 2.17e+09                \\
								\texttt{gng} & GnegativeGO         & Biology & 1392  & 1717  & 8      & 1.91e+07                \\
								\texttt{gnp} & GnegativePseAAC     & Biology & 1392  & 440   & 8      & 4.90e+06                \\
								\texttt{gpg} & GpositiveGO         & Biology & 519   & 912   & 4      & 1.89e+06                \\
								\texttt{gpp} & GpositivePseAAC     & Biology & 519   & 440   & 4      & 9.13e+05                \\
								\texttt{hup} & HumanPseAAC             & Biology & 3106  & 440  & 14     & 1.91e+07                \\
								\texttt{ima} & Image               & Image   & 2000  & 294   & 5      & 2.94e+06                \\
								\texttt{med} & Medical             & Text    & 978   & 1449  & 45     & 6.38e+07                \\
								\texttt{plg} & PlantGO             & Biology & 978   & 3091  & 12     & 3.63e+07                \\
								\texttt{plp} & PlantPseAAC         & Biology & 978   & 440   & 12     & 5.16e+06                \\
								\texttt{sch} & Stackex$\_$chess    & Text    & 1675  & 585   & 227    & 2.22e+08                \\
								\texttt{sco} & Stackex$\_$coffee   & Text    & 225   & 1763  & 123    & 4.88e+07                \\
								\texttt{vig} & VirusGO             & Biology & 207   & 749   & 6      & 9.30e+05                \\
								\texttt{vip} & VirusPseAAC         & Biology & 207   & 440   & 6      & 5.46e+05                \\
								\texttt{yar} & Yahoo$\_$Arts       & Text    & 7484  & 23146 & 26     & 4.50e+09                \\
								\texttt{yed} & Yahoo$\_$Education  & Text    & 12030 & 27534 & 33     & 1.09e+10                \\
								\texttt{yrr} & Yahoo$\_$Recreation & Text    & 12830 & 30324 & 22     & 8.56e+09                \\
								\texttt{yrf} & Yahoo$\_$Reference  & Text    & 8027  & 39679 & 33     & 1.05e+10                \\ \hline
						\end{tabular}}{}
					\end{table}
				\end{example}
				
				
				\begin{sidewaystable}[]
					\centering
					\caption{Numerical results of all algorithms on real multi-label classification datasets.} \label{tab4}
						{\begin{tabular}{c|cccc|cccc|cccc}
								\hline
								& \multicolumn{4}{c|}{\texttt{HL} $\downarrow$}                               & \multicolumn{4}{c|}{\texttt{TIME} (sec) $\downarrow$}         & \multicolumn{4}{c}{\texttt{NNE} $\downarrow$} \\ \hline
								& SGSN     & RSF      & LLSF     & MFDS                         & SGSN     & RSF       & LLSF     & MFDS     & SGSN  & RSF     & LLSF   & MFDS  \\ \hline
								\texttt{bbc} & 1.949e-1 & 1.949e-1 & 1.994e-1 & \multicolumn{1}{c|}{2.426e-1} & 3.856e-2 & 1.364e-1  & 1.541e-1 & 8.103e-2 & 13   & 4224    & 58     & 5394  \\
								\texttt{int} & 1.871e-1 & 1.871e-1 & 1.871e-1 & \multicolumn{1}{c|}{1.988e-1} & 6.921e-2 & 1.169e-1  & 1.674e+0 & 4.506e-1 & 10    & 804     & 109    & 8814  \\
								\texttt{gng} & 1.410e-2 & 4.935e-2 & 1.437e-2 & \multicolumn{1}{c|}{6.291e-2} & 1.704e-1 & 2.078e+0  & 1.042e+0 & 3.584e+0 & 181   & 11160   & 79     & 72    \\
								\texttt{gnp} & 7.674e-2 & 8.026e-2 & 8.080e-2 & \multicolumn{1}{c|}{1.020e-1} & 9.605e-1 & 1.431e+0  & 8.913e-2 & 1.207e+0 & 120   & 3520    & 976    & 3512  \\
								\texttt{gpg} & 3.343e-2 & 8.866e-2 & 3.052e-2 & \multicolumn{1}{c|}{5.087e-2} & 5.279e-2 & 2.027e-1  & 1.170e-1 & 4.915e-1 & 77    & 2664    & 49     & 24    \\
								\texttt{gpp} & 1.599e-1 & 2.137e-1 & 1.817e-1 & \multicolumn{1}{c|}{2.965e-1} & 4.249e-1 & 1.650e-1  & 1.841e-2 & 2.311e-1 & 20    & 1760    & 490    & 56    \\
								\texttt{hup} & 8.221e-2 & 8.242e-2 & 8.235e-2 & \multicolumn{1}{c|}{9.286e-2} & 1.585e+0 & 9.824e+0  & 1.634e-1 & 8.181e+0 & 123   & 378     & 224    & 6146  \\
								\texttt{ima} & 2.192e-1 & 2.473e-1 & 2.473e-1 & \multicolumn{1}{c|}{2.834e-1} & 8.593e-1 & 1.026e+1  & 3.932e-2 & 3.695e+0 & 129   & 105     & 156    & 1175  \\
								\texttt{med} & 1.024e-2 & 1.268e-2 & 1.031e-2 & \multicolumn{1}{c|}{5.726e-2} & 7.167e-1 & 1.030e+0  & 1.001e+0 & 2.196e+0 & 2132  & 53505   & 1809   & 4545  \\
								\texttt{plg} & 3.643e-2 & 8.658e-2 & 4.990e-2 & \multicolumn{1}{c|}{8.832e-2} & 4.637e-1 & 1.503e+0  & 3.650e+0 & 5.607e+0 & 1515  & 4408    & 2584   & 3336  \\
								\texttt{plp} & 8.458e-2 & 8.658e-2 & 8.658e-2 & \multicolumn{1}{c|}{1.742e-1} & 2.042e+0 & 7.156e-1  & 4.372e-2 & 5.203e-1 & 25    & 192     & 52     & 5268  \\
								\texttt{sch} & 9.958e-3 & 1.049e-2 & 9.640e-3 & \multicolumn{1}{c|}{1.315e-2} & 2.187e+0 & 3.015e+0  & 1.774e+0 & 5.188e+0 & 115   & 21565   & 1910   & 227   \\
								\texttt{sco} & 1.570e-2 & 1.581e-2 & 1.558e-2 & \multicolumn{1}{c|}{1.581e-2} & 6.203e-1 & 2.562e-1  & 1.562e+0 & 1.426e+0 & 25    & 13644   & 123    & 1845  \\
								\texttt{vig} & 2.412e-2 & 5.482e-2 & 3.289e-2 & \multicolumn{1}{c|}{8.772e-2} & 7.366e-2 & 3.952e-2  & 6.714e-2 & 1.815e-1 & 288   & 2232    & 229    & 180   \\
								\texttt{vip} & 1.842e-1 & 1.864e-1 & 1.908e-1 & \multicolumn{1}{c|}{1.930e-1} & 1.866e-1 & 3.388e-2  & 2.749e-2 & 3.306e-1 & 25    & 114     & 50     & 1314  \\
								\texttt{eug} & 2.015e-2 & 2.883e-2 & 2.022e-2 & \multicolumn{1}{c|}{8.562e-2} & 1.496e+1 & 3.214e+2  & 1.131e+2 & 7.526e+2 & 2784  & 239316  & 1558   & 13948 \\
								\texttt{yar} & 5.850e-2 & 6.469e-2 & 6.159e-2 & \multicolumn{1}{c|}{6.205e-2} & 2.724e+1 & 4.590e+02 & 2.242e+3 & 1.133e+3 & 1046  & 4420    & 36491  & 2418  \\
								\texttt{yed} & 3.980e-2 & 4.470e-2 & 4.070e-2 & \multicolumn{1}{c|}{4.722e-2} & 4.426e+1 & 1.631e+03 & 2.761e+3 & 1.531e+3 & 680   & 13728   & 34962  & 891   \\
								\texttt{yrr} & 5.173e-2 & 6.568e-2 & 6.357e-2 & \multicolumn{1}{c|}{7.984e-2} & 4.367e+1 & 1.965e+03 & 2.703e+3 & 1.568e+3 & 1825  & 9108    & 20036  & 2684  \\
								\texttt{yrf} & 2.628e-2 & 3.385e-2 & 2.671e-2 & \multicolumn{1}{c|}{3.281e-2} & 3.479e+1 & 8.487e+02 & 2.924e+3 & 1.744e+3 & 1300  & 11451   & 85839  & 2640  \\ \hline
						\end{tabular}}			
						{ \centerline{``$\downarrow$" means that a smaller metric value corresponds to better the algorithm performance.}}
					\end{sidewaystable}

					We report the performance of all the algorithms in terms of \texttt{HL}, \texttt{TIME} and \texttt{NNE}. Particularly, on the large-scale datasets \texttt{yar}, \texttt{yed}, \texttt{yrr}, \texttt{yrf}, we restrict the \texttt{TIME} of LLFS within 3000 seconds, 
					because this algorithm takes longer time to terminate.  
					
					$\bullet$ On \texttt{HL}: We can see from Table \ref{tab4} that SGSN has the best \texttt{HL} on almost all the datasets. LLFS also shows competitive \texttt{HL}, especially on datasets \texttt{gpg} and \texttt{sco}. RFS and MDFS have relatively larger \texttt{HL}.
					
					$\bullet$ On \texttt{TIME}: SGSN shows significant advantage on \texttt{TIME} when solving large-scale dataset, such as \texttt{eug}, \texttt{yar}, \texttt{yed}, \texttt{yrr}, \texttt{yrf}. Particularly, on these datasets, \texttt{TIME} of SGSN is less than 1/10 of other algorithms.  However, on the small dimensional datasets with $d < 600$, LLFS has shorter \texttt{TIME}.
					
					$\bullet$ On \texttt{NNE}: SGSN has the smallest \texttt{NNE} on most datasets in Table \ref{tab4}, especially on the high-dimensional datasets. 
					The smaller \texttt{NNE} is, the faster SGSN runs because the subspace identified becomes
					smaller. 
					This fact also accounts for the shortest \texttt{TIME} of SGSN used on the datasets \texttt{yar}, \texttt{yed}, \texttt{yrr}, \texttt{yrf}. For some low-dimensional datasets, such as \texttt{gng} and \texttt{gpg}, LLFS and MDFS find solutions with smaller \texttt{NNE}.
					
					\section{Conclusions} \label{conclusion}
					
					This paper proposes a new stationary duality approach for
					a class of composite optimization with indicator functions. 
					The dual problem is constructed as $\ell_0$ regularized and nonnegativity constrained sparse optimization. 
					The fundamental one-to-one correspondence between solutions of primal and dual problems is established. 
					To find a dual solution, a 
					subspace gradient semismooth Newton (SGSN) method is developed. 
					It has low per-iteration complexity by exploiting the proximal operator of the sparse regularizer in the dual objective function. 
					Moreover, SGSN enjoys global convergence with a local superlinear (or quadratic) rate under suitable conditions. 
					Extensive experiments on the AUC maximization and multi-label classification problems demonstrate its competitive performance against top solvers for the respective problems.
					
					Although the developed algorithm SGSN explicitly requires the gradient information of the conjugate function $f^*(\cdot)$, 
					the stationary duality result does not rely on the availability of this piece of information.
					In theory, the stationary duality applies to general convex function $f$ as long as the dual problem admits an optimal solution.
					Numerically, this calls for further investigation when the gradient function cannot be
					cheaply computed. 
					For instance, it would be interesting to see if  derivative-free algorithms \cite{conn2009introduction,larson2019derivative}  is a viable choice for
					the dual problem in this situation.
					If successful, it would significantly enlarge the applications of the stationary duality approach.
					We also note that the dual problem provides valuable information for the
					Lagrangian multiplier of the primal problem.
					It would be interesting to investigate whether a primal-dual method is
					possible.
					Those questions are not easy to answer and
					we leave them to our future research.
					
				\end{sloppypar}

				\section*{Acknowledgment}
				This work was supported by the National Key R\&D Program of China (2023YFA1011100), the National Natural Science Foundation of China (12131004), and Research Grants Council of the Hong Kong SAR (PolyU15309223).
				

\bibliographystyle{siamplain}
\bibliography{DCO} 

\end{document}